\documentclass[12pt]{article}
\usepackage{amsmath,amsfonts,amssymb,amsthm}
\newcommand{\1}{{{\bf 1}}}
\newcommand{\id}{{\rm id}}
\newcommand{\Hom}{{\rm Hom}\,}
\newcommand{\End}{{\rm End}\,}
\newcommand{\Res}{{\rm Res}\,}
\newcommand{\Image}{{\rm Im}\,}
\newcommand{\Aut}{{\rm Aut}\,}

\newcommand{\w}{{\omega}}

\newcommand{\Z}{\mathbb{Z}}
\newcommand{\Q}{\mathbb{Q}}
\newcommand{\R}{\mathbb{R}}
\newcommand{\C}{\mathbb{C}}
\newcommand{\N}{\mathbb{N}}
\newcommand{\h}{\mathfrak{h}}

\newcommand{\wt}[1]{{\rm wt}(#1)}

\newtheorem{theorem}{Theorem}[section]
\newtheorem{proposition}[theorem]{Proposition}
\newtheorem{lemma}[theorem]{Lemma}
\newtheorem{corollary}[theorem]{Corollary}

\theoremstyle{definition}
\newtheorem{definition}[theorem]{Definition}

\theoremstyle{remark}
\newtheorem{remark}[theorem]{Remark}

\numberwithin{equation}{section}

\newcommand{\Free}[1]{M(1)^{#1}}
\newcommand{\Fremo}[1]{{M(1,#1)}}
\newcommand{\charge}[1]{V_{L}^{#1}}

\newcommand{\charlam}[1]{V_{#1+L}}

\makeatletter
\@addtoreset{equation}{section}

\makeatother

\newcommand{\Fretw}[1]{M(1)(\theta)^{#1}}
\newcommand{\NO}{\,{\raise0.25em
\hbox{$\mathop{\hphantom{\cdot}}\limits^{_{\circ}}_{^{\circ}}$}}\,}
\newcommand{\fusion}[3]{{\binom{#3}{#1\;#2}}}

\begin{document}

\begin{Large}
\begin{center}
\textbf{Fusion rules for the vertex operator algebras\\
$M(1)^+$ and $V_L^+$}
\end{center}
\end{Large}

\vskip 2ex
\begin{center}
Toshiyuki Abe\footnote{Supported by JSPS Research Fellowships for Young 
Scientists.}
\end{center}
\begin{small}
\begin{center} 
Graduate School of Mathematical Sciences, University of Tokyo\\
3-8-1 Komaba, Tokyo, 153-8914, Japan
\end{center}
\end{small}
\begin{center}
Chongying Dong\footnote{Partially
supported by NSF grants and a research grant from the
Committee on Research, UC Santa Cruz.}\\
\end{center}
\begin{center} 
Department of Mathematics, University of California, Santa Cruz, CA 95064
\end{center}
\begin{center}
and
\end{center}
\begin{center}
Hai-Sheng Li\footnote{Partially
supported by a NSA grant and a grant from Rutgers University Research Council.}
\end{center}
\begin{center} 
Department of Mathematical Sciences, Rutgers University, Camden, NJ 08102, and\\
Department of Mathematics, Harbin Normal University, Harbin, China
\end{center}
\begin{abstract}
The fusion rules for the vertex operator algebras $M(1)^+$
(of any rank) and $V_L^+$ (for any positive definite even lattice $L$) 
are determined completely.
\end{abstract}

\tableofcontents

\section{Introduction}

In this paper we study orbifold vertex operator algebras $\Free{+}$
and $\charge{+}$ for a positive definite even lattice $L.$ The vertex
operator algebra $\charge{+}$ (see \cite{FLM}) 
is the fixed point subalgebra of the lattice
vertex operator algebra $\charge{}$ under the automorphism lifted from
the $-1$ isometry of the lattice and the vertex operator algebra
$\Free{+}$ can be regarded as a subalgebra of $\charge{}.$ The vertex
operator algebra $\charge{+}$ in the case that $L$ is the Leech
lattice was first studied in \cite{FLM} to construct the moonshine
module vertex operator algebra $V^{\natural}$ which is a direct sum of
$\charge{+}$ and an irreducible $\charge{+}$-module in \cite{FLM}.
This construction was extended to some other lattices in \cite{DGM}.

Previously, the vertex operator algebras $\Free{+}$ and $\charge{+}$
have been studied extensively in the literature.  The irreducible
modules for both $\Free{+}$ and $\charge{+}$ have been classified in
\cite{DN1}--\cite{DN3} and \cite{AD}.  If $L$ is of rank 1, the fusion
rules for these vertex operator algebras have been also determined in
\cite{A1} and \cite{A2}.  In this paper we determine the fusion rules
for {\em general} $M(1)^+$ and $V_L^+.$ It turns out that all of the
fusion rules are either $0$ or $1$.

The fusion rules for $M(1)^+$ are obtained in the following way.
First we construct certain untwisted and twisted intertwining
operators which are similar to the untwisted and twisted vertex
operators constructed in Chapters 8 and 9 of \cite{FLM}.  The main
problem is to find the upper bound for each fusion rule.  In order to
achieve this we use a general result about the fusion rules for a
tensor product vertex operator algebra to reduce the problem to the
case when the rank is $1$.  Applying the fusion rules
obtained in \cite{A1} we get the required upper bound. In
particular, the constructed intertwining operators are the only
nonzero intertwining operators up to scalar multiples.

The determination of fusion rules for $V_L^+$ is much more
complicated. The main strategy is to employ the results (on fusion
rules) for $M(1)^+$. (Notice that $M(1)^+$ is a vertex operator
subalgebra of $V_L^+$ and each irreducible $V_L^+$-module is a
completely reducible $M(1)^+$-module.) First, we show that the fusion
rules of certain types are nonzero by exhbiting nonzero intertwining
operators.  Then we prove that the fusion rules for $V_L^+$ are either
$0$ or $1$.  Observe that the intertwining operators constructed in
\cite{DL1} for $V_L$ restrict to nonzero (untwisted) intertwining
operators for $V_L^+.$ We then construct certain (nonzero)
intertwining operators among untwisted and twisted $V_L$-modules and
again restrict to nonzero (twisted) intertwining operators for
$V_L^+.$ The main difficulty is in proving that the constructed
intertwining operators are the all nonzero intertwining operators.
This is achieved by a lengthy calculation involving commutativity and
associativity of vertex operators.

As an application of our main result we show that if $L$ is self dual
and if $V_L^+$ extends to a vertex operator algebra by an irreducible
module from the (unique) twisted $V_L$-module, then the resulted
vertex operator algebra is always holomorphic in the sense that it is
rational and the vertex operator algebra itself is the only
irreducible module.  The moonshine module vertex operator algebra is such an
extension for the Leech lattice and thus it is holomorphic (this result has
been obtained previously in \cite{D3}).  It is expected that the main
result will be useful in the future study of orbifold conformal field
theory for $L$ not self dual.

The organization of the paper is as follows.  Section 2 is 
preliminary; In Section 2.1 we recall definitions of modules for
vertex operator algebras, and in Section 2.2 we review
the notion of intertwining operators and fusion rules and
we also prove that fusion
rules for a tensor product of two vertex operator algebras are equal
to the product of fusion rules for each vertex operator algebra.  In
Section 3.1, we present the construction of
 vertex operator algebras $\Free{+}$ and
$\charge{+}$ and their irreducible modules following \cite{FLM}.  The
classifications of irreducible $\Free{+}$-modules and irreducible
$\charge{+}$-modules given in \cite{DN1}--\cite{DN3} and \cite{AD}
are also stated here.  In Section 3.2 we identify the
contragredient modules of irreducible $\Free{+}$-modules and
$\charge{+}$-modules.  This result is useful to reduce the arguments to
calculate fusion rules.  In Section 4 we determine the fusion rules
for $\Free{+}$ completely.  The nontrivial intertwining operators
among irreducible $\Free{+}$-modules are constructed in Section 4.1,
and it is proved that all of the fusion rules are either $0$ or $1$.

Throughout the paper, $\Z_{\geq0}$ is the set of nonnegative integers.

\section{Preliminaries}\label{CPREL}
\subsection{Vertex operator algebras and modules}
In this section we recall certain basic notions such as the notions of
(weak) twisted module and contragredient module (see 
\cite{FLM}, \cite{FHL}, \cite{DLM3}).

For any vector space $W$ (over $\C$) we set 
\begin{align*}
W[[z,z^{-1}]]
&=\left\{\left.\sum_{n\in\Z}v_nz^{-n-1}\right|\,v_n\in W\,\right\},\\
W((z))&=\left\{\left.\sum_{n\in\Z}v_nz^{-n-1}\right|\,v_n\in W,\,v_n=0
\hbox{ for sufficient small $n$}\,\right\},\\
W\{z\}&=\left\{\left.\sum_{n\in\C}v_nz^{-n-1}\right|\,v_n\in W\,\right\}.
\end{align*}

We first briefly recall the definition of vertex operator algebra (see
\cite{B}, \cite{FLM}).  A vertex operator algebra is a $\Z$-graded
vector space $V=\bigoplus_{n\in \Z} V_{(n)}$ such that $\dim
V_{(n)}<\infty$ for all $n\in \Z$ and such that $V_{(n)}=0$ for $n$
sufficiently small, equipped with a linear map, called the
\textit{vertex operator map},
$$Y(\,\cdot\,,z):V\to(\End V)[[z,z^{-1}]], \quad a\mapsto
Y(a,z)=\sum_{n\in \Z}a_nz^{-n-1}.$$ The vertex operators $Y(a,z)$
satisfy the Jacobi identity.  There are two distinguished vectors; the
{\it vacuum vector} $\1\in V_{(0)}$ and the {\it Virasoro element}
$\w\in V_{(2)}$.  It is assumed that $Y(\1,z)=\id_{V}$ and that the
following Virasoro algebra relations hold for $m,n\in \Z$:
\begin{eqnarray}\label{evirasoro-relations}
[L(m),L(n)]=(m-n)L(m+n)+\frac{1}{12}(m^{3}-m)\delta_{m+n,0}c_{V},
\end{eqnarray}
where $Y(\omega,z)=\sum_{n\in \Z}L(n)z^{-n-2}\,(=\sum_{n\in
\Z}\w_{n}z^{-n-1})$ and $c_{V}$ is a complex scalar, called the
\textit{central charge} of $V$.  It is also assumed that for $n\in
\Z$, the homogeneous subspace $V_{(n)}$ is the eigenspace for $L(0)$
of eigenvalue $n$.  We say that a nonzero vector $v$ of $V_{(n)}$ is a
\textit{homogeneous vector of weight} $n$ and write $\wt v=n$.

Let $V$ be a vertex operator algebra, fixed throughout this section.
An {\em automorphism} of vertex operator algebra $V$ is a linear
isomorphism $g$ of $V$ such that $g(\w)=\w$ and
$gY(a,z)g^{-1}=Y(g(a),z)$ for any $a\in V$.  A simple consequence of
this definition is that $g({\bf 1})={\bf 1}$ and that
$g(V_{(n)})=V_{(n)}$ for $n\in \Z$.  Denote by $\Aut(V)$ the group of
all automorphisms of $V$.  For a subgroup $G<\Aut(V)$ the fixed point
set $V^G=\{a\in V\,|\,g(a)=a\hbox{ for $g\in G$}\}$ is a vertex
operator subalgebra.

Let $g$ be an automorphism of vertex operator algebra $V$ of (finite)
order $T$.  Then $V$ is decomposed into the eigenspaces for $g$:
\[
V=\bigoplus_{r=0}^{T-1}V^{r},\,V^{r}
=\{\,a\in V\,|\,g(a)=e^{-\frac{2 \pi ir}{T}} a\,\}.
\]
\begin{definition}{\rm 
A {\em weak $g$-twisted $V$-module} is a vector space $M$ equipped
with a linear map
\begin{align*}
Y_{M}:V&\to (\End M)\{z\},\\
a&\mapsto Y_{M}(a,z)=\sum_{n\in\Q}a_nz^{-n-1}\quad \left({\rm where }\;a_n\in \End M\right),
\end{align*}
called the {\em vertex operator map}, 
such that the following conditions hold for $0\leq r\leq T-1,\,a\in V^{r},\,b\in V$ and $u\in M$:

(1) $Y_{M}(a,z)v\in z^{-\frac{r}{T}}M((z))$,

(2) $Y_{M}(\1,z)=\id _{M}$,

(3) (the twisted Jacobi identity)
\begin{align*}%\label{Jac1}
\begin{split}
&z_{0}^{-1}\delta\left(\frac{z_{1}-z_{2}}{z_{0}}\right)
Y_{M}(a,z_{1})Y_{M}(b,z_{2})-z_{0}^{-1}\delta\left(
\frac{z_{2}-z_{1}}{-z_{0}}\right)Y_{M}(b,z_{2})Y_{M}(a,z_{1})\\
&\quad=z_{2}^{-1}\left(\frac{z_{1}-z_{0}}{z_{2}}\right)^{-\frac{r}{T}}\delta
\left(\frac{z_{1}-z_{0}}{z_{2}}\right)Y_{M}(Y(a,z_{0})b,z_{2}).
\end{split}
\end{align*}
}
\end{definition}
A weak $g$-twisted $V$-module is denoted by $(M,Y_{M})$, or simply by $M$. 
When $g=1$, a weak $g$-twisted $V$-module is called a \textit{weak $V$-module}. 
A \textit{$g$-twisted weak $V$-submodule} of a $g$-twisted weak module $M$ is 
a subspace $N$ of $M$ such that $a_nN\subset N$ hold for all $a\in V$ and $n\in\Q$. 
If $M$ has no $g$-twisted weak $V$-submodule except $0$ and $M$, $M$ is said to be
 \textit{irreducible}. 

It is known (see \cite{DLM2}) that the operators $L(n)$ for $n\in\Z$ on $M$ 
with $Y_{M}(\w,z)=\sum_{n\in\Z}L(n)z^{-n-2}$ also satisfy 
the Virasoro algebra relations (\ref{evirasoro-relations}).
Moreover, we have the $L(-1)$-derivative property
\begin{equation}\label{DP1} 
Y_{M}(L(-1)a,z)=\frac{d}{dz}Y(a,z)\hbox{ for all }a\in V.
\end{equation}
  
\begin{definition}
{\rm An {\em admissible $g$-twisted $V$-module} is a weak $g$-twisted $V$-module $M$
equipped with a $\frac{1}{T}\N$-grading $M=\bigoplus_{n\in\frac{1}{T}\N} M(n)$ such that 
\begin{align}\label{AD1}
a_mM(n)\subset M(\wt{a}-m-1+n)
\end{align}
for any homogeneous $a\in V$ and for $n\in\frac{1}{T}{\N},\;m\in\Q$.}
\end{definition}
In the case $g=1$, an admissible $g$-twisted $V$-module 
is called an {\em admissible $V$-module}. 
A $g$-twisted weak $V$-submodule $N$ of a $g$-twisted admissible $V$-module 
is called a {\it $g$-twisted admissible $V$-submodule} 
if $N=\bigoplus_{n\in \frac{1}{T}\N} N\cap M(n)$.

A $g$-twisted admissible $V$-module $M$ is said to be \textit{irreducible}
 if $M$ has no nontrivial admissible submodule. 
A $g$-twisted admissible $V$-module $M$ is said to be \textit{completely reducible} 
if $M$ is a direct sum of irreducible admissible submodules. 

\begin{definition}{\rm
The vertex operator algebra $V$ is said to be \textit{$g$-rational} 
if any $g$-twisted admissible $V$-module is completely reducible. 
If $V$ is $\id_{V}$-rational, then $V$ is said to be \textit{rational}.
}
\end{definition} 

\begin{definition}{\rm 
A \textit{$g$-twisted $V$-module} is a weak $g$-twisted $V$-module $M$
which is $\C$-graded by $L(0)$-eigenspace $M=\bigoplus_{\lambda\in\C} M_{(\lambda)}$ 
(where $M_{(\lambda)}=\{u\in M\;|\;L(0)u=\lambda u\}$) 
such that $\dim M_{(\lambda)}<\infty$ for all $\lambda\in \C$ and
such that for any fixed $\lambda\in\C$, 
$M_{(\lambda+n/T)}=0$ for $n\in \Z$ sufficiently small.}
\end{definition}

In the case $g=1$, a $g$-twisted $V$-module is called a \textit{$V$-module}. 
A $V$-module $M$ is said to be \textit{irreducible} 
if $M$ is irreducible as a weak $V$-module.
The vertex operator algebra $V$ is said to be \textit{simple} if
$V$ as a $V$-module is irreducible.

Let $M=\bigoplus_{\lambda\in\C}M_{(\lambda)}$ be a $V$-module. 
Set $M'=\bigoplus_{\lambda\in\C}M_{(\lambda)}^{*}$, 
the restricted dual of $M$.
It was proved in  \cite{FHL} that $M'$ is naturally a $V$-module where 
the vertex operator map, denoted by $Y'$, is defined by the property
\begin{eqnarray}\label{edual}
\langle Y'(a,z)u',v\rangle
=\langle u',Y(e^{zL(1)}(-z^{-2})^{L(0)}a,z^{-1})v\rangle
\end{eqnarray}
for $a\in V,\; u'\in M'$ and $v\in M$.
The $V$-module $M'$ is called the {\em contragredient module} of $M$. 
It was proved therein that if $M$ is irreducible,
then so is $M'$.  A $V$-module $M$ is said to be {\em self-dual} if
$M$ and $M'$ are isomorphic $V$-modules. Then a $V$-module $M$
is self-dual if and only if there exists a nondegenerate invariant
bilinear form on $M$ in the sense that (\ref{edual}) with the obvious
modification holds.
The following result was proved in \cite{liform}:

\begin{lemma}\label{li}
Let $V$ be a simple vertex operator algebra such that $L(1)V_{(1)}\ne
V_{(0)}$. Then $V$ is self-dual.
\end{lemma}

\subsection{Intertwining operators and fusion rules}
We recall the definitions of the notions of intertwining operator 
and fusion rule from \cite{FHL}
and we prove a theorem about fusion rules for a tensor product
vertex operator algebra. 

\begin{definition}\label{define4}
Let $M^1,\,M^2$ and $M^3$ be weak $V$-modules. 
An {\em intertwining operator} $\mathcal{Y}(\,\cdot\,,z)$ of 
type $\fusion{M^1}{M^2}{M^3}$ is a linear map 
\begin{align*}
\mathcal{Y}(\,\cdot\,,z):M^1&\to\Hom(M^2,M^3)\{z\}\\
v^1&\mapsto \mathcal{Y}(v^1,z)=\sum_{n\in\C}v^1_{n}z^{-n-1}\quad \left({\rm where}\;
v^1_{n}\in\Hom(M^2,M^3)\right)
\end{align*}
satisfying the following conditions:

(1) For any $v^1\in M^1,\,v^2\in M^2$ and $\lambda\in\C$, $v^1_{n+\lambda}v^2=0$ 
for $n\in \Z$ sufficiently large.

(2) For any $a\in V,\,v^{1}\in M^1$, 
\begin{align*}%\label{Jac1}
\begin{split}
&z_{0}^{-1}\delta\left(\frac{z_{1}-z_{2}}{z_{0}}\right)Y_{M^3}(a,z_{1})\mathcal{Y}(v^1,z_{2})
-z_{0}^{-1}\delta\left(\frac{z_{2}-z_{1}}{-z_{0}}\right)\mathcal{Y}(v^1,z_{2})Y_{M^2}(a,z_{1})\\
&\quad=z_{2}^{-1}\delta\left(\frac{z_{1}-z_{0}}{z_{2}}\right)\mathcal{Y}(Y_{M^1}(a,z_{0})v^1,z_{2}).
\end{split}
\end{align*}

(3) For $v^1\in M^1$, 
$\frac{d}{dz}\mathcal{Y}(v^1,z)=\mathcal{Y}(L(-1)v^1,z).$
\end{definition}

All of the intertwining operators 
of type $\fusion{M^1}{M^2}{M^3}$ form a vector space,
denoted by $I_V\fusion{M^1}{M^2}{M^3}$. 
The dimension of $I_V\fusion{M^1}{M^2}{M^3}$ is called 
the {\em fusion rule of type $\fusion{M^1}{M^2}{M^3}$ for $V$}.

The following result, which is given \cite{FHL} and \cite{HL}, 
gives the following symmetry:

\begin{proposition}\label{duality}
Let $M,\,N$ and $L$ be $V$-modules.
Then there exist canonical vector space isomorphisms such that
\begin{align*}
I_{V}\fusion{M}{N}{L}\cong I_{V}\fusion{N}{M}{L}\cong I_{V}\fusion{M}{L'}{N'}.
\end{align*}
\end{proposition}

The following proposition can be found in \cite[Proposition 11.9]{DL1}:

\begin{proposition}\label{injectivity}
Let $M^{i}\,(i=1,2,3)$ be $V$-modules.
Suppose that $M^{1}$ and $M^{2}$ are irreducible and 
that $I_{V}\fusion{M^{1}}{M^{2}}{M^{3}}\ne 0$.
Let $\mathcal{Y}(\,\cdot\,,z)$ be any nonzero intertwining operator of type 
$\fusion{M^{1}}{M^{2}}{M^{3}}$. 
Then for any nonzero vectors $u\in M^{1}$ and $v\in M^{2}$, 
$\mathcal{Y}(u,z)v\neq0$.  
\end{proposition}

Assume that $U$ is a vertex operator subalgebra
of $V$ (with the same Virasoro element).
Then every $V$-module is naturally a $U$-module.
Let $M^1,\,M^2,\, M^3$ be $V$-modules and let $N^1$ and $N^2$ be any $U$-submodules of
$M^1$ and $M^2$, respectively. Clearly, any intertwining operator
$\mathcal{Y}(\,\cdot\,,z)$ of type $\fusion{M^1}{M^2}{M^3}$ in the category of
$V$-modules is an intertwining operator of type $\fusion{M^1}{M^2}{M^3}$ 
in the category of $U$-modules. 
Furthermore, the restriction of $\mathcal{Y}(\,\cdot\,,z)$ onto $N^1\otimes N^2$
is an intertwining operator
of type $\fusion{N^1}{N^2}{M^3}$ in the category of $U$-modules. 
Then we have a restriction map
\begin{align*}
I_V\fusion{M^1}{M^2}{M^3}&\to I_U\fusion{N^1}{N^2}{M^3},\\
\mathcal{Y}(\,\cdot\,,z)&\mapsto\mathcal{Y}(\,\cdot\,,z)|N^1\otimes N^2.
\end{align*}
Now, assume that $M^1,\,M^2$ are irreducible $V$-modules and 
$M^3$ is any $V$-module (not necessarily irreducible) and assume that
$N^1$ and $N^2$ are nonzero $U$-modules, e.g.,  irreducible $U$-modules.
It follows immediately from Proposition \ref{injectivity} that 
the restriction map is injective. Therefore we have proved:

\begin{proposition}\label{cfusion-rule-inequality}
Let $V$ be a vertex operator algebra and let $M^1,\,M^2,\, M^3$ be 
$V$-modules among which $M^1$ and $M^2$ are irreducible.
Suppose that $U$ is a vertex operator subalgebra
of $V$ (with the same Virasoro element) and that
$N^1$ and $N^2$ are irreducible $U$-submodules of 
$M^1$ and $M^2$, respectively.
Then the restriction map from $I_V\fusion{M^1}{M^2}{M^3}$ 
to $I_U\fusion{N^1}{N^2}{M^3}$ 
is injective. In particular,
\begin{align}\label{inequality1}
\dim I_V\fusion{M^1}{M^2}{M^3}\leq\dim I_U\fusion{N^1}{N^2}{M^3}.
\end{align}
\end{proposition}

Let $V^{1}$ and $V^{2}$ be vertex operator algebras, 
let $M^{i}\,(i=1,2,3)$ be $V^{1}$-modules and 
let $N^{i}\,(i=1,2,3)$ be $V^{2}$-modules.
For any intertwining operator $\mathcal{Y}_{1}(\,\cdot\,,z)$ of type
$\fusion{M^1}{M^2}{M^3}$ and for any intertwining operator 
$\mathcal{Y}_{2}(\,\cdot\,,z)$ of type
$\fusion{N^1}{N^2}{N^3}$, by using commutativity and rationality one
can prove  that
$\mathcal{Y}_{1}(\,\cdot\,,z)\otimes \mathcal{Y}_{2}(\,\cdot\,,z)$
is an intertwining operator of type 
$\fusion{M^1\otimes N^1}{M^2\otimes N^2}{M^3\otimes M^3}$,
where $(\mathcal{Y}_{1}\otimes \mathcal{Y}_{2})(\,\cdot\,,z)$ is defined by 
\[
(\mathcal{Y}_{1}\otimes \mathcal{Y}_{2})(u^{1}\otimes v^{1},z)u^{2}\otimes v^{2}
=\mathcal{Y}_{1}(u^{1},z)u^{2}\otimes \mathcal{Y}_{2}(v^{1},z)v^{2}
\] 
for $u^{i}\in M^{i}$ and $v^{i}\in N^{i}\,(i=1,\,2)$.  
Then we have a canonical linear map
\begin{align*}
\sigma: I_{V^1}\fusion{M^1}{M^2}{M^3}\otimes I_{V^2}\fusion{N^1}{N^2}{N^3}
&\rightarrow I_{V^1\otimes V^2}\fusion{M^1\otimes N^1}{M^2\otimes N^2}{M^3\otimes M^3}\\
\mathcal{Y}_{1}(\,\cdot\,,z)\otimes \mathcal{Y}_{2}(\,\cdot\,,z)
&\mapsto (\mathcal{Y}_{1}\otimes \mathcal{Y}_{2})(\,\cdot\,,z).
\end{align*}

The following is our main theorem of this section:

\begin{theorem}\label{tensortheorem} 
With the above setting, the linear map $\sigma$ is one-to-one.  
Furthermore, if either
\begin{align*}
\dim I_{V^1}\fusion{M^1}{M^2}{M^3}<\infty\quad\hbox{or}
\quad\dim I_{V^2}\fusion{N^1}{N^2}{N^3}<\infty,
\end{align*}
then $\sigma$ is a linear isomorphism.
\end{theorem}

To prove this theorem we shall need some preparation.
Denote by $\w^{i}$ the Virasoro element of $V^{i}$ for $i=1,2$, and write
\[
Y(\omega^{i},x)=\sum_{n\in \Z}L^{i}(n)x^{-n-2}.
\]

The following proposition is a modification and a
generalization of Proposition 13.18 \cite{DL1}. It can be also 
proved in the same way.

\begin{proposition}\label{pdlgeneral}
Let $V^{1}$ and $V^{2}$ be vertex operator algebras and let
 $W^{i}\,(i=1,2,3)$ be $V^{1}\otimes V^{2}$-modules on which both 
$L^{1}(0)$ and $L^{2}(0)$ act semisimply. Let $\mathcal{Y}(\,\cdot\,,x)$ 
be an intertwining operator of type 
$\fusion{W^1}{W^2}{W^3}$ for $V^{1}\otimes V^{2}$.
Then for any $h\in \C$,
\begin{align*}
x^{-L^{1}(0)}P_{h}\mathcal{Y}(x^{L^{1}(0)}\,\cdot\,,x)x^{L^{1}(0)}\cdot
\end{align*}
is an intertwining operator of type $\fusion{W^1}{W^{2}}{W^{3}(L^{2}(0),h)}$ 
for $V_{1}$-modules, where $W^{3}(L^{2}(0),h)$ is the $L^{2}(0)$-eigenspace 
of $W^{3}$ 
with eigenvalue $h$, which is naturally a (weak) $V_{1}$-module, 
and $P_{h}$ is the projection of $W^{3}$ onto $W^{3}(L^{2}(0),h)$.
\end{proposition}

For a vector space $U$, we say that a formal series $a(x)=\sum_{n\in
\C}a(n)z^{n}\in U\{z\}$ is {\em lower truncated} if $a(n)=0$
for $n$ whose real part is sufficiently small. 
Furthermore, for vector spaces $A$
and $B$, a linear map $g(z)$ from $A$ to $B\{z\}$ is said to be
{\em lower truncated} if $g(z)$ sends every vector in $A$ to a
lower truncated series in $B\{z\}$.

With these notions we formulate the following result, 
which will be very useful in our proof of 
Theorem \ref{tensortheorem}:

\begin{lemma}\label{lsimplefact1}
Let $W=\bigoplus_{h\in\C}W_{(h)}$ be a $\C$-graded vector space 
satisfying the condition
that $\dim W_{(h)}<\infty$ for any $h\in\C$ and that $W_{(h)}=0$ 
for $h$ whose real part is sufficiently small. 
Let $A$ and $B$ be any vector spaces, let $g_{i}(x)\,(i=1,\ldots,r)$ be 
linearly independent lower truncated linear maps from $A$ to $B\{x\}$. 
Suppose that $f_{i}(x)\in W\{x\}\,(i=1,\dots,r)$ are 
lower truncated formal series such that for any $h\in \C$, 
there exists $s\in\C$ such that $P_{h}f_{i}(x)\in x^{s}W_{(h)}$ 
for all $i$, where $P_{h}$ is the projection map of $W$ onto $W_{(h)}$, and such that
\[
f_{1}(x)\otimes g_{1}(x)+\cdots +f_{r}(x)\otimes g_{r}(x)=0
\]
as an element of $\Hom(A,(W\otimes B)\{x\})$.
Then $f_{i}(x)=0$ for all $i$.
\end{lemma}

\begin{proof}
For any $\eta\in W^{*}$, we extend $\eta$ to a linear map from $W\otimes B$ to $B$
by $\eta(w\otimes u)=\eta(w)u$ for $w\in W$ and $u\in B$, and then
 canonically extend it to a linear map from  $W\otimes(B\{x\})$ to $B\{x\}$.
For any $h\in\C,\,\eta\in W^{*}$ and $u\in A$, we see that 
\begin{multline*}
\eta(P_{h}(f_{1}(x)\otimes g_{1}(x)(u)+\cdots +f_{r}(x)\otimes g_{r}(x)(u)))\\=
x^{s}(\eta(w{1})g_{1}(x)(u)+\cdots +\eta(w_{r})g_{r}(x)(u))=0,
\end{multline*}
where we set $P_{h}f_{i}(x)=x^{s}w_{i}$ with $w_{i}\in W_{h}$. 
Since $g_{i}(x)\,(i=1,\ldots,r)$ are linearly independent linear maps from 
$A$ to $B\{x\}$, $\eta(w_{i})=0$ for all $i$.
Thus $w_{i}=0$ for any $h\in\C$ and $i$, that is, $P_{h}f_{i}(x)=0$.
This implies $f_{i}(x)=0$ for all $i$.  
\end{proof}

Now we prove Theorem \ref{tensortheorem}. 
\begin{proof} 
For $h\in \C$, let $P_{h}$ be the projection map of $M^{3}\otimes N^{3}$ 
onto $(M^{3})_{(h)}\otimes N^{3}$.

Suppose that $\mathcal{Y}_{1}^{i}(\,\cdot\,,x)$ for $i=1,\ldots,\,r$ 
are intertwining operators of type $\fusion{M^1}{M^{2}}{M^{3}}$ and 
suppose that $\mathcal{Y}_{2}^{i}(\,\cdot\,,x)$ for $i=1,\ldots,\,r$ 
are linearly independent intertwining operators of type $\fusion{N^1}{N^2}{N^3}$.
Assume that
\begin{align*}
\sum_{i=1}^{r}(\mathcal{Y}_{1}^{i}\otimes \mathcal{Y}_{2}^{i})(\,\cdot\,,x)=0.
\end{align*}
That is,
\begin{eqnarray}\label{ethat-is}
\sum_{i=1}^{r}\mathcal{Y}_{1}^{i}(w^{1},x)w^{2}\otimes \mathcal{Y}_{2}^{i}(v^{1},x)v^{2}=0
\end{eqnarray}
for $w^{j}\in M^{j},\,v^{j}\in N^{j}$ with $j=1,2$.
Write
\begin{eqnarray}
\mathcal{Y}_{1}^{i}(w^{1},x)w^{2}
=\sum_{n\in \C}f^{i}_n(w^{1},w^{2})x^{-n-1}.
\end{eqnarray}
{}From \cite{FHL}, for homogeneous vectors $w^{1},w^{2}$, we have
\begin{align*}
L^{1}(0)f^{i}_n(w^{1},w^{2})
=(\wt{w^{1}}+\wt{w^{2}}-n-1)f^{i}_n(w^{1},w^{2}).
\end{align*}
Then for any $h\in\C$,
\begin{align*}
P_{h}\mathcal{Y}_{1}^{i}(w^{1},x)w^{2}
=f^{i}_{\wt{w^{1}}+\wt{w^{2}}-h-1}(w^1,w^{2})
x^{h-\wt{w^{1}}-\wt{w^{2}}}\in x^{h-\wt{w^{1}}-\wt{w^{2}}}(M^{3})_{(h)}.
\end{align*}
Now it follows immediately from (\ref{ethat-is}) and Lemma \ref{lsimplefact1} that
\begin{align*}
\mathcal{Y}_{1}^{i}(w^{1},x)v^{1}=0\quad\hbox{for $i=1,\ldots,r$.} 
\end{align*}
Thus $\mathcal{Y}_{1}^{i}(\,\cdot\,,x)=0$ for all $i$. 
This proves that $\sigma$ is injective.

Assume $\dim I_{V^2}\fusion{M^2}{N^2}{L^2}<\infty$.
We are going to show that $\sigma$ is also surjective.
Let $\mathcal{Y}(\,\cdot\,,x)$ be any intertwining operator of type 
$\fusion{M^1\otimes N^1}{M^2\otimes N^2}{M^3\otimes N^3}$ for $V^{1}\otimes V^{2}$.
We must prove that $\mathcal{Y}(\,\cdot\,,x)\in \Image\sigma$.

Let $\mathcal{Y}_{2}^{i}(\,\cdot\,,x)\;(i=1,\ldots,\,r)$ be 
a basis of $I_{V^2}\fusion{N^1}{N^2}{N^3}$.
We fix vectors $w^{1}\in M^{1},\,w^{2}\in M^{2}$ arbitrarily.
By Proposition \ref{pdlgeneral}, for $h\in \C$,
\[
x^{-L^{1}(0)}P_{h}\mathcal{Y}(x^{L^{1}(0)}w^{1}\otimes\cdot\,,x)(x^{L^{1}(0)}w^{2}\otimes\cdot\,)
\]
is an intertwining operator of type $\fusion{N^1}{N^2}{(M^3)_{(h)}\otimes N^3}$
for $V_{2}$-modules. (Notice that $(M^3)_{(h)}\otimes N^3$ is the $L^{1}(0)$-eigenspace of
eigenvalue $h$.)
Since $\dim (M^3)_{(h)}<\infty$, we have 
\begin{align*}
I_{V^2}\fusion{N^1}{N^2}{(M^3)_{(h)}\otimes N^3}\cong(M^3)_{(h)}\otimes I_{V^2}\fusion{N^1}{N^2}{N^3}.
\end{align*}
Thus for any $v^{1}\in N^{1}$ and $v^2\in N^2$, we can write 
\begin{align*}
x^{-L_{1}(0)}P_{h}\mathcal{Y}((x^{L_{1}(0)}w^{1})\otimes v^{1},x)(x^{L_{1}(0)}w^{2}\otimes v^{2})
=\sum_{i=1}^{r}f_{i}(w^{1},w^{2},h)\otimes \mathcal{Y}_{2}^{i}(v^{1},x)v^{2}
\end{align*}
for some $f_{i}(w^{1},w^{2},h)\in (M^3)_{(h)}$. 
That is,
\begin{align*}
P_{h}\mathcal{Y}(w^{1}\otimes v^{1},x)(w^{2}\otimes v^{2})
=\sum_{i=1}^{r}x^{h}f_{i}(x^{-L_{1}(0)}w^{1},x^{-L_{1}(0)}w^{2},h)\otimes \mathcal{Y}_{1}^{i}(v^{1},x)v^{2}.
\end{align*}
Then
\begin{align*}
\mathcal{Y}(w^{1}\otimes v^{1},x)(w^{2}\otimes v^{2})
=\sum_{h\in\C}\sum_{i=1}^{r}x^{h}f_{i}(x^{-L_{1}(0)}w^{1},x^{-L_{1}(0)}w^{2},h)\otimes 
\mathcal{Y}_{1}^{i}(v^{1},x)v^{2}
\end{align*}
for any $v^{1}\in N^{1}$ and $v^2\in N^2$.
Now we set
\begin{align*}	
\mathcal{Y}_{1}^{i}(w^{1},x)w^{2}
=\sum_{h\in \C}f_{i}(x^{-L_{1}(0)}w^{1},x^{-L_{1}(0)}w^{2},h)x^{h}.
\end{align*}
Since $M^{3}$ is an ordinary $V$-module, for each $i$,
 $\mathcal{Y}_{1}^{i}(w^{1},x)w^{2}$ 
is a lower truncated element of $M^{3}\{x\}$. 
For example, when $w^{1}\in M^1,\,w^{2}\in M^{2}$ are homogeneous, we have
\begin{align*}
\mathcal{Y}_{1}^{i}(w^{1},x)v^{1}
=\sum_{h\in \C}f_{i}(w^{1},w^{2},h)x^{h-\wt{w^{1}}-\wt{w^{2}}}.
\end{align*}
Then
\begin{align*}
P_{h}\mathcal{Y}_{1}^{i}(w^{1},x)w^{2}\in 
x^{h-\wt{w^{1}}-\wt{w^{2}}}(M^{3})_{(h)}.
\end{align*}
Furthermore, for homogeneous vector $a\in V_{1}$ and for $n\in\Z$, we have
\begin{align*}
P_{h} a_{n}\mathcal{Y}_{1}^{i}(w^{1},x)w^{2}=a_{n}P_{h-\wt{a}+n+1}
\mathcal{Y}_{1}^{i}(w^{1},x)w^{2}
\in x^{h-\wt{a}+n+1-\wt{w^{1}}-\wt{w^{2}}}(M^{3})_{(h)}.
\end{align*}
We are going to prove that $\mathcal{Y}_{1}^{i}(\,\cdot\,,x)$ 
are intertwining operators, 
so that we will have that $\mathcal{Y}(\,\cdot\,,x)\in\Image\sigma$.

Noticing that $L(-1)=L^{1}(-1)\otimes 1+1\otimes L^{2}(-1)$, 
using the $L(-1)$ (resp. $L^{2}(-1)$)-derivative property for 
$\mathcal{Y}(\,\cdot\,,x)$ (resp. $\mathcal{Y}_{2}^{i}(\,\cdot\,,x)$), we get
\begin{align*}
&\sum_{i=1}^{r}\mathcal{Y}_{1}^{i}(L^{1}(-1)w^{1},x)w^{2}\otimes \mathcal{Y}_{2}^{i}(v^{1},x)v^{2}\\
&\quad=\mathcal{Y}(L(-1)(w^{1}\otimes v^{1}),x)(w^{2}\otimes v^{2})
-\sum_{i=1}^{r}\mathcal{Y}_{1}^{i}(w^{1},x)w^{2}\otimes\mathcal{Y}_{2}^{i}(L^{2}(-1)v^{1},x)v^{2}\\
&\quad=\frac{d}{dx}\mathcal{Y}(w^{1}\otimes v^{1},x)(w^{2}\otimes v^{2})
-\sum_{i=1}^{r}\mathcal{Y}_{1}^{i}(w^{1},x)w^{2}\otimes\frac{d}{dx}\mathcal{Y}_{2}^{i}(v^{1},x)v^{2}\\
&\quad=\sum_{i=1}^{r}\left(\frac{d}{dx}\mathcal{Y}_{1}^{i}(w^{1},x)w^{2}\right)\otimes 
\mathcal{Y}_{2}^{i}(v^{1},x)v^{2}.
\end{align*}
Since $\mathcal{Y}_{2}^{i}(\,\cdot\,,x),\,(i=1,\ldots,\,r)$ 
are linearly independent, by Lemma \ref{lsimplefact1} we get
\begin{eqnarray}
\mathcal{Y}_{1}^{i}(L^{1}(-1)w_{1},x)w_{2}
=\frac{d}{dx}\mathcal{Y}_{1}^{i}(w_{1},x)w_{2}
\end{eqnarray}
for any $i=1,\ldots,\,r$ and $w_{j}\in M^{j}\,(j=1,2)$.

Finally, we show that each $\mathcal{Y}_{1}^{i}(\,\cdot\,,x)$ 
satisfies the Jacobi identity.
Let $a\in V^{1},\,w^{1}\in M^{1},\,w^{2}\in M^{2}$.
By linearity we may assume that $a,\,w^{1}$ and $v^{1}$ are homogeneous. 
{}From the Jacobi identity
\begin{align*}
&x_{0}^{-1}\delta\left(\frac{x_{1}-x_{2}}{x_{0}}\right)
Y(a\otimes\1,x_{1})\mathcal{Y}(w^{1}\otimes v^{1},x_{2})(w^{2}\otimes v^{2})\\
&\qquad-x_{0}^{-1}\delta\left(\frac{x_{2}-x_{1}}{-x_{0}}\right)
\mathcal{Y}(w^{1}\otimes v^{1},x_{2})Y(a\otimes\1,x_{1})
(w^{2}\otimes v^{2})\\
&\quad=x_{2}^{-1}\delta\left(\frac{x_{1}-x_{0}}{x_{2}}\right)
\mathcal{Y}(Y(a\otimes\1,x_{0})(w^{1}\otimes v^{1}),x_{2})(w^{2}
\otimes v^{2}),
\end{align*}
we get
\begin{eqnarray}\label{eproof-jacobi-1}
& &\sum_{i=1}^{r}x_{0}^{-1}\delta\left(\frac{x_{1}-x_{2}}{x_{0}}\right)
Y(a,x_{1})\mathcal{Y}_{1}^{i}(w^{1},x_{2})w^{2}\otimes
\mathcal{Y}_{2}^{i}(v^{1},x_{2})v^{2}\nonumber\\
& &\qquad
-\sum_{i=1}^{r}x_{0}^{-1}\delta\left(\frac{x_{2}-x_{1}}{-x_{0}}\right)
\mathcal{Y}_{1}^{i}(w^{1},x_{2})Y(a,x_{1})w^{2}\otimes
\mathcal{Y}_{2}^{i}(v^{1},x_{2})v^{2}\nonumber\\
&=&\sum_{i=1}^{r}x_{2}^{-1}\delta\left(\frac{x_{1}-x_{0}}{x_{2}}\right)
\mathcal{Y}_{1}^{i}(Y(a,x_{0})w^{1},x_{2})w^{2}\otimes
\mathcal{Y}_{2}^{i}(v^{1},x_{2})v^{2}
\end{eqnarray}
for any $v^{j}\in N^{j}\,(j=1,2)$.
For $n\in \Z,\,h\in \C$, we have
\begin{eqnarray}\label{eproof-jacobi-2}
& &\Res_{x_{1}}x_{1}^{n}(x_{1}-x_{2})^{m}
P_{h}Y(a,x_{1})\mathcal{Y}_{1}^{i}(w^{1},x_{2})w^{2}\nonumber\\
& &\quad=\sum_{j=0}^{\infty}\binom{m}{j}(-x_{2})^{j}
P_{h} a_{n+m-j}\mathcal{Y}_{1}^{i}(w^{1},x_{2})v^{1}\nonumber\\
& &\quad=\sum_{j=0}^{\infty}\binom{m}{j}(-x_{2})^{j}a_{n+m-j}
P_{h-\wt{a}+n+m-j+1}\mathcal{Y}_{1}^{i}(w^{1},x_{2})v^{1}\nonumber\\
& &\quad\in x_{2}^{h-\wt{a}-\wt{w^{1}}-\wt{w^{2}}+n+m+1}(M^{3})_{(h)}.
\end{eqnarray}
Similarly, we have
\begin{eqnarray}\label{eproof-jacobi-3}
\Res_{x_{1}}x_{1}^{n}(x_{1}-x_{2})^{m}P_{h}
\mathcal{Y}_{1}^{i}(w^{1},x_{2})Y(a,x_{1})w^{2}
\in x_{2}^{h-\wt{a}-\wt{w^{1}}-\wt{w^{2}}+n+m+1}(M^{3})_{(h)},
\end{eqnarray}
and
\begin{eqnarray}\label{eproof-jacobi-4}
& &\Res_{x_{0}}\Res_{x_{1}}x_{0}^{m}x_{1}^{n}
x_{2}^{-1}\delta\left(\frac{x_{1}-x_{0}}{x_{2}}\right)P_{h}
\mathcal{Y}_{1}^{i}(Y(a,x_{0})w^{1},x_{2})w^{2}\nonumber\\
& &\quad=\Res_{x_{0}}x_{0}^{m}(x_{2}+x_{0})^{n}P_{h}
\mathcal{Y}_{1}^{i}(Y(a,x_{0})w^{1},x_{2})w^{2}\nonumber\\
& &\quad=\sum_{j=0}^{\infty}\binom{n}{j}x_{2}^{n-j}
P_{h}\mathcal{Y}_{1}^{i}(a_{m+j}w^{1},x_{2})v^{1}\nonumber\\
& &\quad\in x_{2}^{h-\wt{a}-\wt{w^{1}}-\wt{w^{2}}+n+m+1}(M^{3})_{(h)}.
\end{eqnarray}
With (\ref{eproof-jacobi-1})--(\ref{eproof-jacobi-4}),
it follows from Lemma \ref{lsimplefact1} that each 
$\mathcal{Y}_{1}^{i}(\,\cdot\,,x)$ satisfies the Jacobi identity. 
Then $\mathcal{Y}_{1}^{i}(\,\cdot\,,x)$ are intertwining operators.
This shows that $\sigma$ is onto, completing the proof.
\end{proof}

\section{Vertex operator algebras $\Free{+}$ and $\charge{+}$}
\subsection{Vertex operator algebras $\Free{+}$ and $\charge{+}$ 
and their modules}\label{ssss}
In this section we review the construction of the vertex operator
algebras $\Free{+}$ and $V_L^+$ associated with a positive definite
even lattice $L$, following \cite{FLM}.

Let $\h$ be a $d$-dimensional vector space equipped with a nondegenerate symmetric bilinear form $(\cdot\,,\cdot)$. 
Consider the Lie algebra $\hat{\h}=\h\otimes\C[t,t^{-1}]\oplus\C C$ 
defined by the commutation relations 
\[
[\beta_{1}\otimes t^{m},\,\beta_{2}\otimes t^{n}]
=m(\beta_1,\beta_2)\delta_{m,-n}C\hbox{ and }[C,\hat{\h}]=0
\] 
for any $\beta_1,\beta_2\in\h,\,m,\,n\in\Z$. 
Set
\[
\hat{\h}^+=\C[t]\otimes\h\oplus\C C,
\]
which is clearly an abelian subalgebra.
For any $\lambda\in\h$, let $\C e^\lambda$ denote the 1-dimensional 
$\hat{\h}^+$-module 
on which $\h\otimes t\C[t]$ acts as zero, $\h\;(=\h\otimes \C t^0)$ acts
according to the character $\lambda$, i.e., 
$he^\lambda=(\lambda,h)e^\lambda$ for $h\in \h$ and 
and $C$ acts as the scalar $1$.
Set
\[
M(1,{\lambda})=U(\hat{\h})\otimes_{U(\hat{\h}^+)}\C e^\lambda\cong S(t^{-1}\C[t^{-1}]\otimes \h),
\]
the induced $\hat{\h}$-module.

For $h\in\h,\,n\in\Z$, we denote by $h(n)$ the corresponding operator 
of $h\otimes t^{n}$ on $M(1,{\lambda})$, and write 
\[
h(z)=\sum_{n\in\Z}h(n)z^{-n-1}.
\]
Define a linear map 
\begin{align}\label{module}
Y(\,\cdot\,,z):\Fremo{0}\to(\End\Fremo{\lambda})[[z,z^{-1}]]
\end{align} 
by
\begin{align*}
Y(v,z)=\NO\frac{1}{(n_{1}-1)!}\left(\frac{d}{dz}\right)^{n_1-1}\beta_1(z)
\cdots\frac{1}{(n_{r}-1)!}\left(\frac{d}{dz}\right)^{n_r-1}\beta_r(z)\NO,
\end{align*}
for the vector $v=\beta_1(-n_1)\cdots\beta_{r}(-n_r)e^{0}$ with $\beta_i\in\h,\,n_i\geq1$, where the normal ordering $\NO\,\cdot\,\NO$ is an operation which reorders the operators so that $\beta(n)\,(\beta\in\h,n< 0$) to be placed to the left of $\beta(n)\,(\beta\in\h,n\geq 0$).

Following \cite{FLM}, we denote $\Free{}=M(1,0)$ and set
\[
\1=e^{0}\in M(1),\quad\w=\frac{1}{2}\sum_{i=1}^{d}h_i(-1)^2e^{0}\in \Free{},
\]
where $\{h_{1},\dots,h_{d}\}$ is an orthonormal basis of $\h$. 
(Note that $\w$ does not depend on the choice of the orthonormal basis.)
Then $(\Free{},\,Y(\,\cdot\,,z),\,\1,\,\w)$ is a simple vertex operator algebra, 
and $(\Fremo{\lambda},Y(\,\cdot\,,z))$ is an irreducible $\Free{}$-module 
for any $\lambda\in\h$ (see \cite{FLM}).

We next recall a construction of the vertex operator algebra $\charge{}$ 
associated to an even lattice and its irreducible modules,  
following \cite{DL1} (see also \cite{FLM} and \cite{D1}).
First we start with a rank $d$ rational lattice $P$ with a positive definite 
symmetric $\Z$-bilinear form $(\cdot\,,\cdot)$. 
We suppose that $L$ is a rank $d$ even sublattice of $P$ such that $(L,P)\subset\Z$. 

Let $q$ be a positive even integer such that $(\lambda,\mu)\in \frac{2}{q} \Z$ 
for all $\lambda,\mu\in P$ and let $\hat{P}$ be
a central extension of $P$ by the cyclic group $\langle\kappa_{q}\rangle$ of order $q\,$:
\[
1\rightarrow \langle \kappa_{q}|{\kappa_{q}}^{q}=1\rangle\rightarrow 
\hat{P}\mathop{\rightarrow}\limits^{-}P\rightarrow 0
\]
with commutator map $c(\cdot\,,\cdot)$ such that 
$c(\alpha,\beta)={\kappa}^{\frac{(\alpha,\beta)}{2}}$ for $\alpha,\,\beta\in L$, 
where $\kappa={\kappa_{q}}^{q/2}$.
It is known that such a central extension exists if $q$ is sufficiently large 
(see Remark 12.18 in \cite{DL1}).
Let $e:P\to\hat{P},\,\lambda\mapsto e_{\lambda}$ be a section such that 
$e_0=1$ and $\epsilon: P\times P\to\langle\kappa_{q}\rangle$ be
the corresponding 2-cocycle, i.e., 
$e_{\lambda}e_{\mu}=\epsilon(\lambda,\mu)e_{\lambda+\mu}$ for any $\lambda,\,\mu\in P$. 
We can assume that $\epsilon$ is bimultiplicative.
Then $\epsilon(\alpha,\beta)\epsilon(\beta,\alpha)=\kappa^{(\alpha,\,\beta)},$
$\epsilon(\alpha+\beta,\gamma)=\epsilon(\alpha,\gamma)\epsilon(\beta,\gamma)$.
We may further assume that 
\[
\epsilon(\alpha,\alpha)=\kappa^{\frac{(\alpha,\alpha)}{2}}
\]
for any $\alpha\in L$.

Denote by $\C[P]=\bigoplus_{\lambda\in P}\C e^{\lambda}$ 
the group algebra. For any subset $M$ of $P$, we write
$\C[M]=\bigoplus_{\lambda\in M}\C e^{\lambda}$.
Then $\C[P]$ becomes a $\hat{P}$-module by the action 
\begin{eqnarray}\label{e-twisted-group-action}
e_{\lambda}e^{\mu}
=\epsilon(\lambda,\mu)e^{\lambda+\mu}\quad\hbox{and}\quad\kappa_{q} 
e^{\mu}=\omega_{q}e^{\mu}
\end{eqnarray}
for $\lambda,\mu\in P$, where $\omega_{q}\in\C^{\times}$ 
is a $q$-th root of unity. 
It is clear that for any $\lambda\in P$, $\C[\lambda+L]$ 
is an $\hat{L}$-module on which $\kappa(={\kappa_{q}}^{q/2})$ 
acts by the scalar $-1$. 

Set $\h=\C\otimes_{\Z}L$ and extend the $\Z$-bilinear form 
$(\cdot\,,\cdot)$ to a $\C$-bilinear form of $\h$.
Then
\[
V_{P}:=M(1)\otimes \C[P]
\]
is endowed with an $\hat{\h}$-module structure such that 
\begin{align*}
h(n)(u\otimes e^{\lambda})=(h(n)u)\otimes e^{\lambda}\quad\hbox{and}\quad
h(0)(u\otimes e^{\lambda})=(h,\lambda)(u\otimes e^{\lambda})
\end{align*}
for $h\in\h,\,n\neq 0,\, \lambda\in P$ and that $C$ acts as the identity. 
We have 
\begin{eqnarray*}
V_{P}\cong\bigoplus_{\lambda\in P}\Fremo{\lambda},
\end{eqnarray*}
as an $\Free{}$-module. 
For any subset $M$ of $P$, we set
$V_{M}=M(1)\otimes \C[M]$, which is an $\Free{}$-submodule of $V_{P}$,
where $\C[M]=\bigoplus_{\lambda\in M}\C e^{\lambda}$.
 
For $\lambda\in P$, we define $Y(e^{\lambda},z)\in (\End V_{P})\{z\}$ by  
\begin{eqnarray}\label{eYoperator-def}
&Y(e^{\lambda},\,z)
=\exp\left(\sum_{n=1}^{\infty}\frac{\lambda(-n)}{n} z^{n}\right)
\exp\left(-\sum_{n=1}^{\infty}\frac{\lambda(n)}{n}z^{-n}\right)e_{\lambda}z^{\lambda},
\end{eqnarray}
where $e_\lambda$ is the left action of $e_{\lambda}\in\hat{P}$ on $\C[P]$ and $z^{\lambda}$ is the operator on $\C[P]$ defined by $z^{\lambda}e^\mu=z^{(\lambda,\mu)}e^\mu$. 
The vertex operator associated to the vector $v=\beta_1(-n_1)\cdots\beta_{r}(-n_r)e^\lambda$ for $\beta_i\in\h,\,n_i\geq1$ and $\lambda\in P$ is defined by 
\begin{align*}
Y(v,z)=\NO\frac{1}{(n_{1}-1)!}\left(\frac{d}{dz}\right)^{n_1-1}\beta_1(z)
\cdots\frac{1}{(n_{r}-1)!}\left(\frac{d}{dz}\right)^{n_r-1}\beta_r(z)Y(e^\lambda,z)\NO,
\end{align*}
where the normal ordering $\NO\,\cdot\,\NO$ is an operation which reorders the operators so that $\beta(n)\,(\beta\in\h,n< 0$) and $e_{\lambda}$ to be placed to the left of $X(n),\,(X\in\h,n\geq 0$) and $z^{\lambda}$.
This defines a linear map 
\begin{align}\label{ffff}
Y(\,\cdot\,,z):V_{P}\to(\End V_{P})\{z\}.
\end{align}
Let $\alpha,\lambda\in P$ be such that $(\alpha,\lambda)\in\Z$.
Then for $u\in\Fremo{\alpha},\,v\in\Fremo{\lambda}$, we have
\begin{align}\label{iiii}
\begin{split}
& z_{0}^{-1}\delta\left({\frac{z_{1}-z_{2}}{z_{0}}}\right)Y(u,z_{1})Y(v,z_{2})-(-1)^{(\alpha,\lambda)}c(\alpha,\lambda)z_{0}^{-1}\delta\left({\frac{z_{2}-z_{1}}{-z_{0}}}\right)Y(v,z_{2})Y(u,z_{1})\\
&\quad=z_{2}^{-1}\delta\left(\frac{z_{1}-z_{0}}{z_{2}}\right)Y(Y(u,z_{0})v,z_{2}).
\end{split}
\end{align}

Set 
\begin{eqnarray*}
L^\circ=\{\,\lambda\in\h\,|\,(\alpha,\lambda)\in\Z\,\},
\end{eqnarray*}
the dual lattice of $L$, and we fix a coset decomposition 
$L^{\circ}=\cup_{i\in L^{\circ}/L}(L+\lambda_i)$ such that $\lambda_0=0$.
In the case $P=L^{\circ}$, we see that $V_{P}=\bigoplus_{i\in L^{\circ}/L}\charlam{\lambda_{i}}$ and that the restriction of $Y(\,\cdot\,,z)$ to $\charge{}$ gives a linear map $\charge{}\to(\End  \charlam{\lambda_{i}})[[z,z^{-1}]]$ for any $i\in L^{\circ}/L$.  
{}From \cite{FLM}, $(\charge{},Y(\,\cdot\,,z),\1,\w)$ 
is a vertex operator algebra 
and $(\charlam{\lambda},Y(\,\cdot\,,z))$ are irreducible $\charge{}$-modules.
Note that $M(1)$ is a vertex operator subalgebra of $V_L$
(with the same vacuum vector and the Virasoro element). 

Now we define a map $\theta$ from $\hat{L^{\circ}}$ to itself by
\[
\theta(\kappa_{q}^{s}e_{\lambda})=\kappa_{q}^{s}e_{-\lambda}
\] 
for any $s\in\Z$ and $\lambda\in L^{\circ}$.
Since the $2$-cocycle $\epsilon$ is bimulticative, $\theta$ 
is in fact an automorphism of $\hat{L^{\circ}}$.
Now we define the action of $\theta$ on $V_{L^{\circ}}$ by  
\begin{align*}
\theta(\beta_{1}(-n_{1})\beta_{2}(-n_{2})\cdots \beta_{k}(-n_{k})
e^{\lambda})=(-1)^{k}\beta_{1}(-n_{1})\beta_{2}(-n_{2})\cdots \beta_{k}(-n_{k})
e^{-\lambda}
\end{align*}
for $\beta_i\in\h,\,n_i\geq1$ and $\lambda\in L^{\circ}$.
Then we see that
\begin{align}\label{oooo}
\theta Y(u,z)v=Y(\theta(u),z)\theta(v)
\end{align}
for any $u,v\in V_{L^{\circ}}$.
In particular, $\theta$ gives an automorphism of $V_L$ 
which induces an automorphism of $M(1)$. 

For any $\theta$-stable subspace $U$ of $V_{L^\circ}$, 
let $U^\pm$ be the $\theta$-eigenspace of 
$U$ (of eigenvalues $\pm 1$).
Then both $(\Free{+},Y(\,\cdot\,,z),\1,\w)$ and $(\charge{+},Y(\,\cdot\,,z),\1,\w)$ 
are simple vertex operator algebras.
We have the following  proposition (see \cite{DM} and \cite{DLM1}):

\begin{proposition}\label{untwisted} 
(1) $\Free{\pm},\,\Fremo{\lambda}$ for $\lambda\in\h-\{0\}$ are irreducible
$M(1)^+$-modules, and $\Fremo{\lambda}\cong\Fremo{-\lambda}.$  

(2) $(\charlam{\lambda_i}+\charlam{-\lambda_i})^{\pm}$ for $i\in L^{\circ}/L$ 
are irreducible $\charge{+}$-modules. 
Moreover if $2\lambda_i\not\in L$ then 
$(\charlam{\lambda_i}+\charlam{-\lambda_i})^{\pm}$,  
$\charlam{\lambda_i}$ and $\charlam{-\lambda_i}$ are isomorphic $V_{L}^{+}$-modules. 
\end{proposition}

Next we recall a construction of $\theta$-twisted modules for
$\Free{}$ and $\charge{}$ following \cite{FLM} and \cite{D2}.  Denote
by $\h[-1]=\h\otimes t^{\frac{1}{2}}\C[t,t^{-1}]\oplus\C C$ the
twisted affinization of $\h$ defined by the commutation relations
\begin{align*}
[\beta_1\otimes t^{m},\beta_2\otimes t^{n}]
=m(\beta_1,\beta_2)\delta_{m,-n}C\quad\hbox{ and }\quad[C,\hat{\h}]=0
\end{align*} 
for any $\beta_1,\beta_2\in\h,\,m,\,n\in\frac{1}{2}+\Z$. 
Set
\begin{eqnarray*}
\Fretw{}=S(t^{-\frac{1}{2}}\C[t^{-1}]\otimes\h)
\end{eqnarray*}
Then $\Fretw{}$ is (up to equivalence)
the unique irreducible $\hat{\h}[-1]$-module such that $C=1$ and
$(\beta\otimes t^n)\cdot 1=0$ if $n>0.$ 
This space is an irreducible $\theta$-twisted $\Free{}$-module
(see \cite{FLM}).

Set $K=\{a^{-1}\theta(a)\,|\,a\in\hat{L}\}$.
For any $\hat{L}/K$-module $T$ such that $\kappa$ acts by the scalar $-1$, we define $\charge{T}=\Fretw{}\otimes T$.
Then there exists a linear map $Y(\,\cdot\,,z):\charge{}\to(\End \charge{T})[[z^{\frac{1}{2}},z^{-\frac{1}{2}}]]$ such that $(\charge{T},\,Y(\,\cdot\,,z))$ becomes a $\theta$-twisted $\charge{}$-module (see \cite{FLM}).  
The cyclic group $\langle\theta\rangle$ acts on $M(1)(\theta)$ and $\charge{T}$ by 
\begin{align*}
\theta(\beta_{1}(-n_{1})\beta_{2}(-n_{2})\cdots \beta_{k}(-n_{k}))=(-1)^{k}\beta_{1}(-n_{1})\beta_{2}(-n_{2})\cdots \beta_{k}(-n_{k})
\end{align*}
and
\begin{align}\label{action4}
\theta(\beta_{1}(-n_{1})\beta_{2}(-n_{2})\cdots \beta_{k}(-n_{k})t)
=(-1)^{k}\beta_{1}(-n_{1})\beta_{2}(-n_{2})\cdots \beta_{k}(-n_{k})t
\end{align}
for $\beta_i\in\h,\,n_i\in \frac{1}{2}+\Z_{\geq0}$ and $t\in T$. 
We denote by $M(1)(\theta)^{\pm}$ and $\charge{T,\pm}$ the $\pm1$-eigenspaces 
for $\theta$ of $M(1)(\theta)$ and $\charge{T}$, respectively.

Following \cite{FLM}, let ${T_{\chi}}$ be the irreducible 
$\hat{L}/K$-module associated to a central character $\chi$ satisfying 
$\chi(\kappa)=-1.$ Then any irreducible $\theta$-twisted $\charge{}$-module 
is isomorphic to $\charge{T_{\chi}}$ for some central character $\chi$ with $\chi(\kappa)=-1$ (see \cite{D2}). 
{}From \cite{DLi}  we have:

\begin{proposition}\label{twisted}
(1) $\Fretw{\pm}$ are irreducible $\Free{+}$-modules.

(2) Let $\chi$ be an central character of $\hat{L}/K$ such that $\chi(\kappa)=-1$, 
and $T_{\chi}$ the irreducible $\hat{L}/K$-module with central character $\chi$. 
Then $\charge{+}$-modules $\charge{T_{\chi},\pm}$ are irreducible.
\end{proposition}

The following classification of the irreducible $\Free{+}$-modules is due to
\cite{DN1} and \cite{DN3}:

\begin{theorem}\label{freeclass} 
The $M(1)^{+}$-modules
\begin{equation}\label{IM1}
\Free{\pm},\Fretw{\pm},\Fremo{\lambda}(\cong \Fremo{-\lambda})
\;\;\mbox{ for }\lambda\in\h-\{0\}
\end{equation} 
are all the irreducible $\Free{+}$-modules (up to equivalence).
\end{theorem}

Furthermore, the following classification of the 
irreducible $\charge{+}$-modules was obtained in
\cite{DN2} and \cite{AD}:

\begin{theorem}\label{chargeclass} 
Let $L$ be a positive-definite even lattice and let $\{\lambda_{i}\}$ be 
a set of representatives of $L^{\circ}/L$.
Then any irreducible $\charge{+}$-module is isomorphic to one of the irreducible modules 
$\charge{\pm},\,\charlam{\lambda_{i}}$ with $2\lambda_{i}\notin L$, 
$\charlam{\lambda_{i}}^{\pm}$ with $2\lambda_{i}\in L$ or $\charge{T_{\chi},\pm}$ 
for a central character $\chi$ of $\hat{L}/K$ with $\chi(\kappa)=-1$.
Furthermore, $\charlam{\lambda_{i}}$ and $\charlam{\lambda_{j}}$ is isomorphic 
if and only if $\lambda_{i}\pm\lambda_{j}\in L$.
\end{theorem}

We refer to the irreducible $\charge{+}$-modules 
$\charge{\pm},\,\charlam{\lambda}\,(2\lambda\notin L)$ 
and $\charlam{\lambda}^{\pm}\,(2\lambda\in L)$ 
as the \textit{irreducible modules of untwisted type} 
and  refer to $\charge{T_{\chi},\pm}$ 
as the \textit{irreducible modules of twisted type}.

\subsection{Contragredient modules}
In this section we identify the contragredient modules 
of the irreducible $\Free{+}$-modules and $\charge{+}$-modules explicitly.

First we have:

\begin{proposition}\label{contfreeboson}
Every irreducible $\Free{+}$-module $W$ is self dual, i.e., $W'\cong W$. 
\end{proposition}

\begin{proof}
First, since $M(1)^{+}$ is simple and $M(1)^{+}_{(1)}=0$, by Lemma
\ref{li} $M(1)^{+}$ is self-dual. Similarly, the vertex operator
algebra $M(1)$ is also
self-dual because $L(1)M(1)_{(1)}=L(1)\h=0$. Then as an $M(1)^{+}$-module
$$M(1)'=(M(1)^{+})'\oplus (M(1)^{-})'\simeq M(1)=M(1)^{+}\oplus M(1)^{-}.$$
Since $M(1)^{+}$ and $M(1)^{-}$ are nonisomorphic irreducible $M(1)^{+}$-modules,
we must have that $M(1)^{-}$ is self-dual.

We claim that for any $\lambda\in \h$, 
$M(1,\lambda)'\simeq M(1,-\lambda)$ as an $M(1)$-module.
Note that the lowest $L(0)$-weight subspace of $M(1,\lambda)$ is
$\C e^{\lambda}$ whose $L(0)$-weight is $(\lambda,\lambda)/2$.
Define a linear functional $\psi\in M(1,\lambda)'$ by
$\psi(e^{\lambda})=1$ and $\psi(u)=0$ for $u\in M(1,\lambda)_{(n)}$
with $n-(\lambda,\lambda)/2\in\Z_{>0}$.
{}From (\ref{edual}) we get $h(0)\psi=-(\lambda,h)\psi$ and $h(n)\psi=0$ for
$h\in \h,\; n\ge 1$. 
Thus $M(1,\lambda)'\simeq M(1,-\lambda)$ as an $\hat{\h}$-module,
since $M(1,\lambda)'$ and $M(1,-\lambda)$ are 
irreducible $\hat{\h}$-modules.
Now that $M(1,\lambda)$ and $M(1,-\lambda)$ are isomorphic
$M(1)^{+}$-modules, we see that $M(1,\lambda)$ as an $M(1)^{+}$-module
is self-dual.

It remains to show that the irreducible $M(1)^{+}$-modules
$M(1)(\theta)^{+}$ and $M(1)(\theta)^{-}$ are self-dual. 
It is known that
the lowest $L(0)$-weights of $M(1)(\theta)^{+}$ and $M(1)(\theta)^{-}$
are $\dim \h/16$ and $1/2+\dim \h/16$, respectively. 
Noticing that any irreducible module and its contragredient module
have the same lowest weight $L(0)$-weight, we see that
$M(1)(\theta)^{\pm }$ must be self-dual.
\end{proof}

Combining Proposition \ref{contfreeboson} with Proposition
\ref{duality} we immediately have:

\begin{proposition}\label{dddd}
Let $M^{i}\,(i=1,2,3)$ be irreducible $\Free{+}$-modules. 
Then the fusion rule of type $\fusion{M^{1}}{M^{2}}{M^{3}}$ 
as a function of $(M^{1},M^{2},M^{3})$ is invariant
under the permutation group of $\{1,2,3\}$.
\end{proposition}

Next we identify the contragredient modules of
the irreducible $\charge{+}$-modules:

\begin{proposition}\label{eeee}
The irreducible $\charge{+}$-modules $\charge{\pm}$ and
$\charlam{\lambda}$ for $\lambda\in L^{\circ}$ with $2\lambda\notin L$ are self dual. 
For any $\lambda\in L^{\circ}$ with $2\lambda\in L$, 
$\charlam{\lambda}^{\pm}$ are self dual if $2(\lambda,\lambda)$ 
is even and $(\charlam{\lambda}^{\pm})'\cong\charlam{\lambda}^{\mp}$ 
if $2(\lambda,\lambda)$ is odd.
Let $\chi$ be a central character of $\hat{L}/K$ such that $\chi(\kappa)=-1$. 
Then the irreducible modules $(\charge{T_{\chi},\pm})'$ are isomorphic
to $\charge{T_{\chi'},\pm}$, 
where $\chi'$ is a central character of $\hat{L}/K$ 
defined by $\chi'(a)=(-1)^{\frac{(\bar{a},\bar{a})}{2}}\chi(a)$ for any $a\in Z(\hat{L}/K)$.
\end{proposition}

\begin{proof}
We first prove that  for $\lambda\in L^{\circ}$
$(\charlam{\lambda})'\cong\charlam{-\lambda}$ 
as a $\charge{}$-module.
Since $\charlam{\lambda}=\oplus_{\alpha\in L}M(1,\lambda+\alpha)$ and since
$(\Fremo{\lambda})'\cong\Fremo{-\lambda}$ as an $\Free{}$-module
(from the proof of Proposition \ref{contfreeboson}),
we have $(\charlam{\lambda})'\cong \oplus_{\alpha\in L} M(1,-\lambda+\alpha)$.
By the classification of irreducible $\charge{}$-modules (see \cite{D1}), 
we must have $(\charlam{\lambda})'\cong\charlam{-\lambda}$.  
Since $\charlam{\lambda}\cong\charlam{-\lambda}$ as
a $\charge{+}$-module 
we see that $\charlam{\lambda}$ as a $\charge{+}$-module is self dual.

Now suppose that $2\lambda\in L$. Then $\lambda+L=-\lambda+L$, 
so that $\charlam{\lambda}'\cong \charlam{\lambda}$.
We have a nondegenerate $V_{L}$-invariant bilinear form 
$\langle\,\cdot\,,\cdot\,\rangle$ 
on $\charlam{\lambda}$. 
{}From the invariance property we have
\[
\langle h(n)u,v\rangle=-\langle u, h(-n)v\rangle
\]
for $h\in \h,\; n\in \Z,\; u,v\in V_{\lambda+L}$, 
noticing that $L(1)h=0$ and $L(0)h=h$.
Thus we get
$\langle e^{\lambda}, e^{-\lambda+\alpha}\rangle =0$ for nonzero $\alpha\in L$.
Since the bilinear form is nondegenerate, we must have that 
$\langle e^{\lambda}, e^{-\lambda}\rangle\neq 0$.
%We may assume that $\langle e^{\lambda},e^{-\lambda}\rangle=1$. 
By (\ref{eYoperator-def}) and (\ref{e-twisted-group-action}) we have
\[
Y(e^{2\lambda},z)e^{-\lambda}
=\epsilon(2\lambda,-\lambda)z^{-2(\lambda,\lambda)}
\exp\left(\sum_{n\ge 1}\frac{2\lambda(-n)}{n}z^{n}\right)e^{\lambda}.
\]
Using this and the invariance property we have
\begin{align*}
\langle Y(e^{2\lambda},z)e^{-\lambda}, e^{-\lambda}\rangle
&=\epsilon(2\lambda,-\lambda)z^{-2(\lambda,\lambda)}
\langle \exp\left(\sum_{n\ge 1}\frac{2\lambda(-n)}{n}z^{n}\right)e^{\lambda},
e^{-\lambda}\rangle\\
&=\epsilon(2\lambda,-\lambda)z^{-2(\lambda,\lambda)}
\langle e^{\lambda},\exp\left(\sum_{n\ge 1}\frac{2\lambda(n)}{n}z^{n}\right)
e^{-\lambda}\rangle\\
&=\epsilon(2\lambda,-\lambda)z^{-2(\lambda,\lambda)}\langle e^{\lambda},
e^{-\lambda}\rangle.
\end{align*}
On the other hand, we have
\begin{align*}
\langle e^{-\lambda}, Y(e^{zL(1)}(-z^{-2})^{L(0)}e^{2\lambda},z^{-1})e^{-\lambda}\rangle
&=\langle e^{-\lambda}, (-1)^{2(\lambda,\lambda)}z^{-4(\lambda,\lambda)}
Y(e^{2\lambda},z^{-1})e^{-\lambda}\rangle\nonumber\\
&=(-1)^{2(\lambda,\lambda)}\epsilon(2\lambda,-\lambda)
z^{-2(\lambda,\lambda)}
\langle e^{-\lambda}, e^{\lambda}\rangle,
\end{align*}
noticing that $L(1)e^{2\lambda}=0$ and 
$L(0)e^{2\lambda}=2(\lambda,\lambda)e^{2\lambda}$,
where $2(\lambda,\lambda)$ is a nonnegative integer. 
By the invariance property we have
\begin{align*}
\langle e^{\lambda},e^{-\lambda}\rangle
=(-1)^{2(\lambda,\lambda)}\langle e^{-\lambda},e^{\lambda}\rangle.
\end{align*}
This shows that 
\[
\langle e^{\lambda}\pm e^{-\lambda},e^{\lambda}\pm(-1)^{2(\lambda,\lambda)}e^{-\lambda}\rangle
=\pm2.
\] 
The irreducibility of $\charlam{\lambda}^{\pm}$ and the $V$-invariance of 
$\langle\,\cdot\,,\cdot\,\rangle$ prove that if $2(\lambda,\lambda)$ is even (resp. odd), 
then $\langle\,\cdot\,,\,\cdot\,\rangle$ gives a nondegenerate invariant bilinear form on $\charlam{\lambda}^{\pm}\times\charlam{\lambda}^{\pm}$ 
(resp. $\charlam{\lambda}^{\pm}\times\charlam{\lambda}^{\mp}$). 
Therefore, $(\charlam{\lambda}^{\pm})'\cong\charlam{\lambda}^{\pm}$ 
if $2(\lambda,\lambda)$ is even and 
$(\charlam{\lambda}^{\pm})'\cong\charlam{\lambda}^{\mp}$ 
if $2(\lambda,\lambda)$ is odd. 

Let $\chi$ be a central character of $\hat{L}/K$ such that $\chi(\kappa)=-1$.
Then $(\charge{T_{\chi}})'$ is a $\theta$-twisted $\charge{}$-module
(see \cite{X}; cf. \cite{FHL}). 
The classification of irreducible $\theta$-twisted modules (see
\cite{D2}) implies that $(\charge{T_{\chi}})'$ is isomorphic to 
$\charge{T_{\chi_{1}}}$ for some central character $\chi_{1}$.
We are going to show that $\chi_{1}=\chi'$, using the same method
that was used for the untwisted modules.

For $\alpha\in L$, we have (\cite[Section 9.1]{FLM}.)
\[
Y(e^{\alpha},z)=2^{-(\alpha,\alpha)} z^{-(\alpha,\alpha)/2}
\exp\left( \sum_{n\in 1/2+\Z_{\geq0}}\frac{\alpha(-n)}{n}z^{n}\right)
\exp\left(-\sum_{n\in 1/2+\Z_{\geq0}}\frac{\alpha(n)}{n}z^{-n}\right)e_{\alpha},
\]
so that
\[
Y(e^{\alpha},z)t=\chi(e^{\alpha}) 2^{-(\alpha,\alpha)} z^{-(\alpha,\alpha)/2}
\exp\left( \sum_{n\in 1/2+\Z_{\geq0}}\frac{\alpha(-n)}{n}z^{n}\right)t
\]
for $t\in T_{\chi}$ and $\alpha\in\bar{R}=\{\,\bar{a}\,|\,a\in Z(\hat{L}/K)\,\}$. Then for any $\alpha\in\bar{R},\,t\in T_{\chi}$
and $t_{1}\in T_{\chi_{1}}$, we have 
\begin{align*}
\langle Y(e^{\alpha},z)t_{1},t\rangle
&=2^{-(\alpha,\alpha)} z^{-(\alpha,\alpha)/2}
\langle t_{1},t\rangle\\
\intertext{and}
\langle t_{1},Y(e^{zL(1)}(-z^{-2})^{L(0)}e^{\alpha},z^{-1})t\rangle&=(-1)^{(\alpha,\alpha)/2}
\chi(e^\alpha) 2^{-(\alpha,\alpha)} z^{-(\alpha,\alpha)/2}
\langle t_{1},t\rangle.
\end{align*}
Therefore, we get  
$
\chi(e^{\alpha})\langle t_{1},t\rangle
=(-1)^{\frac{(\alpha,\alpha)}{2}}\chi_{1}(e^{\alpha})\langle t_{1},t\rangle
$
for any $\alpha\in\bar{R},\,t\in T_{\chi}$ and $t_{1}\in T_{\chi_{1}}$.
This proves $\chi_{1}=\chi'$ and 
$(\charge{T_{\chi}})'\cong\charge{T_{\chi'}}$.
Then it is clear that
$(\charge{T_{\chi},\pm})'\cong\charge{T_{\chi'},\pm}$ 
as a $\charge{+}$-module. 
\end{proof}

\section{Fusion rules for vertex operator algebra $\Free{+}$}
\subsection{Construction of intertwining operators}
In this subsection we prove that the fusion rules of certain types
are not zero for vertex operator algebra $\Free{+}$
by constructing a nonzero intertwining operator. 
This construction of intertwining operator is essentially due to \cite{FLM}.

For any $\lambda,\mu,\nu\in\h$ we call the triple 
$(\lambda,\,\mu,\,\nu)\in\h\times\h\times\h$ an {\em admissible
triple} if $p\lambda+q\mu+r\nu=0$ for some $p,q,r\in\{\pm1\}$.
Clearly, if $(\lambda,\,\mu,\,\nu)$ is admissible,
so is every permutation of $(\lambda,\,\mu,\,\nu)$.
Note that in view of Theorem \ref{freeclass},
$\Fremo{\lambda}$ and $\Fremo{\mu}$ are 
isomorphic $\Free{+}$-modules if and only if $(0,\lambda,\mu)$ 
is an admissible triple.

For $\lambda,\mu\in\h$,
we define a linear map $p_{\lambda}:\Fremo{\mu}\to\Fremo{\lambda+\mu}$ by $p_{\lambda}(u\otimes e^{\mu})=u\otimes e^{\lambda+\mu}$.
The vertex operator associated to the vectors $e^\lambda$ and $v=\beta_1(-n_1)\cdots\beta_{r}(-n_r)e^\lambda$ for $\beta_i\in\h,\,n_i\geq1$ is defined by 
\begin{align}
&\mathcal{Y}_{\lambda,\mu}(e^{\lambda},z)=\exp\left(\sum_{n=1}^{\infty}\frac{\lambda(-n)}{n} z^{n}\right)\exp\left(-\sum_{n=1}^{\infty}\frac{\lambda(n)}{n}z^{-n}\right)p_{\lambda}z^{\lambda},\label{untwistinter1}\\
\begin{split}\label{untwistinter2}
&\mathcal{Y}_{\lambda,\mu}(v,z)\\
&\quad=\NO\left(\frac{1}{(n_1-1)!}\left(\frac{d}{dz}\right)^{n_1-1}\beta_1(z)\right)\cdots\left(\frac{1}{(n_r-1)!}\left(\frac{d}{dz}\right)^{n_r-1}\beta_r(z)\right)\mathcal{Y}_{\lambda,\mu}(e^\lambda,z)\NO,
\end{split}
\end{align}
where $z^{\lambda}$ is the operator on $\C e^{\mu}$ defined by $z^{\lambda}e^\mu=z^{(\lambda,\mu)}e^\mu$, and the normal ordering $\NO\,\cdot\,\NO$ is an operation which reorders the operators so that $\beta(n)\,(\beta\in\h,n< 0$) and $p_{\lambda}$ to be placed to the left of $\beta(n),\,(\beta\in\h,n\geq 0$) and $z^{\lambda}$.

{}From the arguments in \cite[Section 8]{FLM}, we see that the operator 
\begin{align}\label{untwistop}
\mathcal{Y}_{\lambda,\mu}(\,\cdot\,,z):\Fremo{\lambda}\to\Hom(\Fremo{\mu},\Fremo{\lambda+\mu})\{z\}
\end{align} 
satisfies 
\begin{align}\label{ggggg}
\begin{split}
& z_{0}^{-1}\delta\left({\frac{z_{1}-z_{2}}{z_{0}}}\right)\mathcal{Y}_{\lambda,\mu+\nu}(u,z_{1})\mathcal{Y}_{\mu,\nu}(v,z_{2})\\
&-(-1)^{(\lambda,\mu)}z_{0}^{-1}\delta\left({\frac{z_{2}-z_{1}}{-z_{0}}}\right)\mathcal{Y}_{\mu,\lambda+\nu}(v,z_{2})\mathcal{Y}_{\lambda,\nu}(u,z_{1})\\
&\quad=z_{2}^{-1}\delta\left({\frac{z_{1}-z_{0}}{z_{2}}}\right)\mathcal{Y}_{\lambda+\mu,\nu}(\mathcal{Y}_{\lambda,\mu}(u,z_{1})v,z_{2})
\end{split}
\end{align}
for $\lambda,\mu,\nu\in\h$ with $(\lambda,\mu)\in\Z,\,u\in M(1,\lambda)$
and $v\in M(1,\nu).$ 
We also have the $L(-1)$-derivative property $\frac{d}{dz}\mathcal{Y}_{\lambda,\mu}(u,z)=\mathcal{Y}_{\lambda,\mu}(L(-1)u,z).$ 
Noting $\mathcal{Y}_{0,\nu}(\,\cdot\,,z)$ is the vertex operator map of the irreducible $\Free{}$-module $\Fremo{\nu}$, we see that $\mathcal{Y}_{\lambda,\mu}(\,\cdot\,,z)$ is a nonzero intertwining operator 
of type $\fusion{\Fremo{\lambda}}{\Fremo{\mu}}{\Fremo{\lambda+\mu}}$
for $\Free{}$. Consequently, the fusion rule of type 
$\fusion{\Fremo{\lambda}}{\Fremo{\mu}}{\Fremo{\lambda+\mu}}$ for
$\Free{+}$ is not zero.
Since $\Fremo{\nu}\cong\Fremo{-\nu}$ as an $\Free{+}$-module for
any $\nu\in\h$, $\fusion{\Fremo{\lambda}}{\Fremo{\mu}}{\Fremo{-\lambda+\mu}}$
for $\Free{+}$ is not zero. Therefore we have proved:

\begin{proposition}\label{nonzero1}
For any admissible triple $(\lambda,\mu,\nu)$, the fusion rule of type 
$\fusion{\Fremo{\lambda}}{\Fremo{\mu}}{\Fremo{\nu}}$ 
for $\Free{+}$ is nonzero.
\end{proposition}
  
For any $\lambda\in\h$, we define a linear map 
\begin{align}\label{tttt}
\theta:\Fremo{\lambda}\to\Fremo{-\lambda};
\quad\theta(u\otimes e^{\lambda})
=\theta(u)\otimes e^{-\lambda}\quad\hbox{for $u\in\Free{}$}.
\end{align} 
For $h\in \h,\,u\in M(1)$, we have
\begin{eqnarray*}
(\theta\circ h(0)\circ \theta^{-1})(u\otimes e^{\lambda})
=\theta h(0)(\theta^{-1}(u)\otimes e^{-\lambda})
=(h,-\lambda)u\otimes e^{-\lambda}
=-h(0)(u\otimes e^{\lambda})
\end{eqnarray*}
and for $n\neq 0$, we have
\begin{align*}
(\theta\circ h(n)\circ \theta^{-1})(u\otimes e^{\lambda})
=\theta ((h(n)\theta^{-1}(u))\otimes e^{-\lambda})
=(\theta h(n)\theta^{-1}(u))\otimes e^{\lambda}
=-h(n)(u\otimes e^{\lambda}).
\end{align*}
Therefore, we see that $\theta\circ h(z)\circ \theta^{-1}=-h(z)$ for any $h\in\h$. 
Since $\theta\circ p_{\lambda}\circ\theta^{-1}=p_{-\lambda}$ for any $\lambda\in\h$, one has $\theta\circ \mathcal{Y}_{\lambda,-\mu}(e^{\lambda},z)\circ \theta^{-1}= \mathcal{Y}_{-\lambda,\mu}(e^{-\lambda},z)$.
%Since $\h$ generates $M(1)$ as a vertex algebra and $\theta$ is an automorphism of $M(1)$, using induction we have
%\begin{eqnarray}
%\theta( Y(u,x)w)=Y(\theta(u),x)\theta(w)\;\;\;
%\mbox{ for }u\in M(1),\; w\in M(1,\lambda).
%\end{eqnarray}
%This implies that $\theta$ is an $\Free{+}$-module isomorphism.
By using \eqref{untwistinter2} we can prove that the intertwining operator 
$\mathcal{Y}_{\lambda,\mu}(\,\cdot\,,z)$ 
satisfies that 
\begin{align}\label{untwist1}
\theta\mathcal{Y}_{\lambda,-\mu}(u,z)\theta^{-1}(v)
=\mathcal{Y}_{-\lambda,\mu}(\theta(u),z)v
\end{align}
for any $u\in\Fremo{\lambda}$ and $v\in\Fremo{\mu}$. 
By using the isomorphism $\theta$, we define an operator 
\begin{align*}
{}^\theta\mathcal{Y}_{\lambda,\mu}(\,\cdot\,,z):
\Fremo{\lambda}\to\Hom(\Fremo{\mu},\Fremo{-\lambda+\mu})\{z\}
\end{align*}
by
\begin{align*}
{}^\theta\mathcal{Y}_{\lambda\mu}(u,z)v
=\mathcal{Y}_{-\lambda,\mu}(\theta(u),z)v
\end{align*}
for $u\in\Fremo{\lambda}$ and $v\in\Fremo{\mu}$.
Then one can see that ${}^\theta\mathcal{Y}_{\lambda\mu}(\,\cdot\,,z)$ is 
a nonzero intertwining operator of type 
$\fusion{\Fremo{\lambda}}{\Fremo{\mu}}{\Fremo{-\lambda+\mu}}$ 
for $\Free{+}$ by using \eqref{untwist1}.

Now we consider the case $\lambda=\mu=\nu=0$ in Proposition \ref{nonzero1}. 
Since $M(1)$ is simple, $Y(u,z)v\ne 0$ for nonzero vectors
$u,\,v\in\Free{}$ by Proposition \ref{injectivity}. 
Clearly, we have
\begin{align*}
Y(u,z)v\in
\begin{cases}
\Free{+}((z)) &\hbox{\rm if $u\in\Free{\pm}$ and $v\in\Free{\pm}$,}\\
\Free{-}((z)) &\hbox{\rm if $u\in\Free{\pm}$ and $v\in\Free{\mp}$.}
\end{cases}
\end{align*}
The restrictions of $Y(\,\cdot\,,z)$ give nonzero intertwining operators 
of types $\fusion{\Free{\pm}}{\Free{\pm}}{\Free{+}}$ and $\fusion{\Free{\pm}}{\Free{\mp}}{\Free{-}}$. 
Thus we have:

\begin{proposition}\label{nonzero2}
The fusion rules of types $\fusion{\Free{\pm}}{\Free{\pm}}{\Free{+}}$ and 
$\fusion{\Free{\pm}}{\Free{\mp}}{\Free{-}}$ are nonzero.
\end{proposition} 
 
Next we consider the case $\lambda=0$ and $\mu\neq0$ in Proposition \ref{nonzero1}. 
Notice that the vertex operator map
$Y(\,\cdot\,,z)$ of the irreducible $\Free{}$-module $\Fremo{\mu}$ is an intertwining operator. 
Then the restrictions of $Y(\,\cdot\,,z)$ give intertwining operators 
of types $\fusion{\Free{\pm}}{\Fremo{\mu}}{\Fremo{\mu}}$. 
By Proposition \ref{injectivity},  $Y(u,z)v\ne 0$ 
for any nonzero vectors $u\in\Free{}$ and $v\in\Fremo{\mu}$. 
Therefore the following proposition holds:

\begin{proposition}\label{nonzero3}
For any $\mu\in\h-\{0\}$, the fusion rules of types 
$\fusion{\Free{\pm}}{\Fremo{\mu}}{\Fremo{\mu}}$ are nonzero.
\end{proposition} 

We shall discuss the construction of intertwining operators of type 
$\fusion{\Fremo{\lambda}}{\Fretw{\epsilon_{1}}}{\Fretw{\epsilon_{2}}}$ 
for $\lambda\in\h$ and $\epsilon_{i}\in\{\pm\}\,(i=1,2)$. 
Let $\lambda\in\h$. 
Following \cite{FLM}, we define a linear map
\begin{align}\label{twining1}
\mathcal{Y}^{\rm tw}_{\lambda}(\,\cdot\,,z):\Fremo{\lambda}\to(\End\Fretw{})\{z\}
\end{align}
as follows.
First we set   
\begin{align}
&\mathcal{Y}^{\rm tw}_{\lambda}(e^{\lambda},z)\nonumber\\
&\quad=e^{-|\lambda|^{2}\log 2}z^{-\frac{|\lambda|^{2}}{2}}\exp\left(\sum_{n\in\frac{1}{2}+\Z_{\geq0}}\frac{\lambda(-n)}{n} z^{n}\right)\exp\left(-\sum_{n\in\frac{1}{2}+\Z_{\geq0}}\frac{\lambda(n)}{n}z^{-n}\right).\label{twistinter1}
\end{align}
Next we define $W(u,z)$ for $u=\beta_1(-n_1)\cdots\beta_{r}(-n_r)e^\lambda\,(\beta_i\in\h,\,n_i\geq1)$ by
\begin{align}
&W(u,z)\nonumber\\
&\quad=\NO\left(\frac{1}{(n_1-1)!}\left(\frac{d}{dz}\right)^{n_1-1}\beta_1(z)\right)\cdots\left(\frac{1}{(n_r-1)!}\left(\frac{d}{dz}\right)^{n_r-1}\beta_r(z)\right)\mathcal{Y}^{\rm tw}_{\lambda}(e^\alpha,z)\NO,\label{twistinter2}
\end{align}
where the normal ordering $\NO\,\cdot\,\NO$ reorders the operators so that $\beta(n)\,(\beta\in\h,\,n<0)$ to be placed to the left of $\beta(n),\,(\beta\in\h,\,n>0)$.
Now we introduce an operator $\Delta_{z}$ defined by
\begin{align*}
\Delta_{z}=\sum_{i=1}^{d}\sum_{m,n=0}^{\infty}c_{mn}h_{i}(m)h_{i}(n)z^{m+n}
\end{align*}
by using an orthonormal basis $\{h_i\}$ of $\h$ and the coefficients $c_{mn}$ subject to the following formal expansion
\begin{equation*}
\sum_{m,n\geq 0}c_{mn}x^{m}y^{n}=-\log
\left(\frac{(1+x)^{\frac{1}{2}}+(1+y)^{\frac{1}{2}}}{2}\right),
\end{equation*}
Finally we set $\mathcal{Y}^{\rm tw}_{\lambda}(u,z)=W(e^{\Delta_{z}}u,z)$.
Then by using the same arguments in \cite[Chapter 9]{FLM}, we get the following twisted Jacobi identity
\begin{align*}%\label{jacobi1}
\begin{split}
& z_{0}^{-1}\delta\left({\frac{z_{1}-z_{2}}{z_{0}}}\right)Y(a,z_{1})\mathcal{Y}^{\rm tw}_{\lambda}(u,z_{2})-z_{0}^{-1}\delta\left({\frac{z_{2}-z_{1}}{-z_{0}}}\right)\mathcal{Y}^{\rm tw}_{\lambda}(u,z_{2})Y(a,z_{1})\\
&\quad=\frac{1}{2}\sum_{p=0,1}z_{2}^{-1}\delta\left((-1)^{p}{\frac{(z_{1}-z_{0})^{1/2}}{z_{2}^{1/2}}}\right)\mathcal{Y}^{\rm tw}_{\lambda}(Y(\theta^{p}(a),z_{0})u,z_{2})
\end{split}
\end{align*}
for any $a\in\Free{}$ and $u\in\Fremo{\lambda}$ and the $L(-1)$-derivative property $\frac{d}{dz}\mathcal{Y}^{\rm tw}_{\lambda}(u,z)=\mathcal{Y}^{\rm tw}_{\lambda}(L(-1)u,z)$ for $u\in\Fremo{\lambda}$.
These imply that $\mathcal{Y}^{\rm tw}_{\lambda}(\,\cdot\,,z)$ is a nonzero intertwining operator of type $\fusion{\Fremo{\lambda}}{\Fretw{}}{\Fretw{}}$ for $\Free{+}$.
By definition we have 
\begin{align}\label{twistconjugate}
\theta\mathcal{Y}^{\rm tw}_{\lambda}(u,z)\theta^{-1}(v)=\mathcal{Y}^{\rm tw}_{-\lambda}(\theta(u),z)v
\end{align}
for any $u\in\Fremo{\lambda}$ and $v\in\Fretw{}$.

Let $p_{\epsilon}:\Fretw{}\to\Fretw{\epsilon}$ be the canonical projection and $\iota_{\epsilon}:\Fretw{\epsilon}\to\Fretw{}$ the canonical inclusion for $\epsilon\in\{\pm\}$. 
Then for any $\epsilon_{1},\,\epsilon_{2}\in\{\pm\}$, the composition $p_{\epsilon_{2}}\circ\mathcal{Y}^{\rm tw}_{\lambda}(\,\cdot\,,z)\circ\iota_{\epsilon_{1}}$ is an intertwining operator of type $\fusion{\Fremo{\lambda}}{\Fretw{\epsilon_{1}}}{\Fretw{\epsilon_{2}}}$ for $\Free{+}$.  
By direct calculation, one has  
\begin{align*}
\mathcal{Y}^{\rm tw}(e^{\lambda},z)1\equiv e^{-|\lambda|^{2}\log 2}z^{-\frac{|\lambda|^{2}}{2}}\left(1+\lambda\left(-{1}/{2}\right) z^{1/2}\right)\mod{z^{-\frac{|\lambda|^{2}}{2}+1}\Fretw{}[[z^{\frac{1}{2}}]]}
\end{align*}
and  
\begin{align*}
\mathcal{Y}^{\rm tw}(e^{\lambda},z)\lambda(-1/2)\equiv&e^{-|\lambda|^{2}\log 2}z^{-\frac{|\lambda|^{2}}{2}}\left(-|\lambda|^{2}z^{-1/2}+(1-2|\lambda|^2)\lambda\left(-{1}/{2}\right)z^{0}\right.\\
&+2(1-2|\lambda|^2)\lambda\left(-{1}/{2}\right)^2 z^{1/2}+\left.\left(4\lambda\left(-{1}/{2}\right)^3-\frac{2}{3}\lambda\left(-{3}/{2}\right)\right)z\right)\\
&\mod{z^{-\frac{|\lambda|^{2}}{2}+2}\Fretw{}[[z^{\frac{1}{2}}]]}.
\end{align*}
These show that if $\lambda$ is nonzero then the intertwining operator 
$p_{\epsilon_{2}}\circ\mathcal{Y}^{\rm tw}_{\lambda}(\,\cdot\,,z)\circ\iota_{\epsilon_{1}}$ 
is nonzero 
for any $\epsilon_{1},\,\epsilon_{2}\in\{\pm\}$. 
Therefore, the following proposition holds:
\begin{proposition} 
For any $\lambda\in\h-\{0\}$, the fusion rules of types 
$\fusion{\Fremo{\lambda}}{\Fretw{\pm}}{\Fretw{\pm}}$ and 
$\fusion{\Fremo{\lambda}}{\Fretw{\pm}}{\Fretw{\mp}}$ are nonzero. 
\end{proposition}

In the case $\lambda=0$, $\mathcal{Y}^{\rm tw}_{0}(\,\cdot\,,z)$ is
the vertex operator map $Y(\,\cdot\,,z)$ of the $\theta$-twisted
$\Free{}$-module $\Fretw{}$.  In particular, $\mathcal{Y}^{\rm
tw}_{0}(\1,z)=\id$ and $\mathcal{Y}^{\rm tw}_{0}(h(-1)\1,z)=h(z)$ for
any $h\in\h$.  Thus $\mathcal{Y}^{\rm
tw}_{0}(\,\cdot\,,z)\circ\iota_{\epsilon}$ is a nonzero intertwining
operator of type
$\fusion{\Fremo{\lambda}}{\Fretw{\epsilon}}{\Fretw{}}$ for any
$\epsilon\in\{\pm\}$.  By using the conjugation property
\eqref{twistconjugate} we immediately have:

\begin{proposition}\label{twistfree} 
The fusion rules of types $\fusion{\Free{\pm}}{\Fretw{\pm}}{\Fretw{+}},\,
\fusion{\Free{\pm}}{\Fretw{\mp}}{\Fretw{-}}$ are nonzero. 
\end{proposition}

\subsection{Main theorem}
In this section we determine the fusion rules for irreducible
$M(1)^{+}$-modules, generalizing a result of \cite{A1}.

The following result was proved in \cite{A1}:

\begin{theorem}\label{yyyy} 
Let $\h$ be a 1-dimensional vector space equipped with
a symmetric nondegenerate bilinear form $(\cdot,\cdot)$.
For any irreducible $\Free{+}$-modules $M^i\,(i=1,2,3)$, 
the fusion rule of type $\fusion{M^1}{M^2}{M^3}$ is either $0$ or $1$ and 
it is invariant under the permutations of $\{1,2,3\}$. 
The fusion rule of type $\fusion{M^1}{M^2}{M^3}$ is $1$ if and only if 
$M^i\,(i=1,2,3)$ satisfy the following conditions:

\noindent
{\rm(i)} $M^1=\Free{+}$ and $M^2\cong M^3$.

\noindent
{\rm(ii)} $M^1=\Free{-}$ and $(M^2,M^3)$ is one of the following pairs:
\begin{itemize}
\item[] $(\Free{+},\Free{-}),\,(\Free{-},\Free{+})$,
\item[] $(\Fremo{\mu},\Fremo{\nu})$ for $\mu,\nu\in\h-\{0\}$ such that 
$\mu=\pm\nu$,
\item[] $(\Fretw{+},\Fretw{-}),\,(\Fretw{-},\Fretw{+})$.
\end{itemize} 

\noindent
{\rm(iii)} $M^1=\Fremo{\lambda}\,(\lambda\in\h-\{0\})$ and 
$(M^2,M^3)$ is one of the following pairs:
\begin{itemize}
\item[] $(\Free{\pm},\Fremo{\mu}),\,(\Fremo{\mu},\Free{\pm})$ 
for $\mu\in\h-\{0\}$ such that $\lambda=\pm\mu$,
\item[] $(\Fremo{\mu},\Fremo{\nu})$ for $\mu,\nu\in\h-\{0\}$ 
such that $(\lambda,\mu,\nu)$ is an admissible triple,
\item[] $(\Fretw{\pm},\Fretw{\pm}),\,(\Fretw{\pm},\Fretw{\mp})$.
\end{itemize} 

\noindent
{\rm(iv)} $M^1=\Fretw{+}$ and $(M^2,M^3)$ is one of the following pairs:
\begin{itemize}
\item[] $(\Free{\pm},\Fretw{\pm}),\,(\Fretw{\pm},\Free{\pm})$,
\item[] $(\Fremo{\mu},\Fretw{\pm}),\,
(\Fretw{\pm},\Fremo{\mu})\,(\mu\in\h-\{0\})$.
\end{itemize}

\noindent
{\rm(v)} $M^1=\Fretw{-}$ and $(M^2,M^3)$ is one of the following pairs:
\begin{itemize}
\item[] $(\Free{\pm},\Fretw{\mp}),\,(\Fretw{\pm},\Free{\mp})$,
\item[] $(\Fremo{\mu},\Fretw{\pm}),\,(\Fretw{\pm},
\Fremo{\mu})\,(\mu\in\h-\{0\})$.
\end{itemize}
\end{theorem}

This section is devoted to prove the following generalization:

\begin{theorem}\label{ehigher-rank-m1}
Let $\h$ be any finite-dimensional vector space equipped with
a symmetric nondegenerate bilinear form $(\cdot\,,\cdot)$.
Then all the assertions of Theorem \ref{yyyy} hold.
\end{theorem}
 
We write $M_{\h}(1)$ for the vertex operator algebra $M(1)$ 
associated with $\h$ and similarly for the modules.
It is clear that if $\h'$ is a subspace of $\h$ such that 
the bilinear form of $\h$ restricted to $\h'$ is nondegenerate,
then $M_{\h'}(1)^{+}$ is a vertex operator subalgebra 
$M_{\h}(1)^{+}$ (with different Virasoro element if $\h'\neq\h$). 
Furthermore, if $\h=\h_1\oplus\h_2$ such that $(\h_1,\h_2)=0$, 
then the irreducible $M_{\h}(1)^{+}$-modules are decomposed into 
direct sums of 
irreducible $M_{\h_1}(1)^{+}\otimes M_{\h_2}(1)^{+}$-modules as follows:
\begin{align}
M_{\h}(1)^{+}&\cong M_{\h_1}(1)^{+}\otimes M_{\h_2}(1)^{+}\oplus
M_{\h_1}(1)^{-}\otimes M_{\h_2}(1)^{-},\label{decomp1}\\
M_{\h}(1)^{-}&\cong M_{\h_1}(1)^{+}\otimes M_{\h_2}(1)^{-}
\oplus M_{\h_1}(1)^{-}\otimes M_{\h_2}(1)^{+},\label{decomp2}\\
M_{\h}(1,\lambda)&\cong M_{\h_1}(1,\lambda_1)\otimes
M_{\h_2}(1,\lambda_2),\label{decomp3}\\
M_{\h}(1)(\theta)^+&\cong M_{\h_1}(1)(\theta)^+\otimes M_{\h_2}(1)(\theta)^+
\oplus M_{\h_1}(1)(\theta)^-\otimes M_{\h_2}(1)(\theta)^-,\label{decomp4}\\
M_{\h}(1)(\theta)^-&\cong M_{\h_1}(1)(\theta)^+\otimes M_{\h_2}(1)(\theta)^-
\oplus M_{\h_1}(1)(\theta)^-\otimes M_{\h_2}(1)(\theta)^+,\label{decomp5}
\end{align}
where we decompose $\lambda\in\h$ into $\lambda=\lambda_1+\lambda_2$
so that $\lambda_i\in\h_i$.

First we prove the following result:
\begin{proposition}\label{fusionfreeinq}
For any irreducible $\Free{+}$-modules $M,\,N$ and $L$, the
fusion rule of type $\fusion{M}{N}{L}$ is either $0$ or $1$. 
\end{proposition}

\begin{proof}
We shall use induction on $d=\dim \h$. Noticing that
Theorem \ref{ehigher-rank-m1} in the case $d=\dim \h=1$ 
has been proved in \cite{A1} (Theorem \ref{yyyy}), we assume that $d>1$.
Assume that Theorem \ref{ehigher-rank-m1} 
for $M_{\h'}(1)^{+}$ with $\dim\h'<d$ has been proved. 
We decompose $\h$ into a direct sum of mutually orthogonal subspaces 
$\h_1$ and $\h_2$ with $\dim \h_1=1$.
Theorem \ref{ehigher-rank-m1} applies for both $M_{h_1}(1)^{+}$ and
$M_{\h_2}(1)^{+}$.
Recall \eqref{decomp1}--\eqref{decomp5} for the decompositions of
the irreducible $M_{\h}(1)^{+}$-modules into direct sums of 
irreducible $M_{\h_1}(1)^+\otimes M_{\h_2}(1)^+$-modules.
Notice that each of $M,\,N$ and $L$ is isomorphic to one of those
$M_{\h}(1)^{+}$-modules.

Pick up irreducible 
$M_{\h_1}(1)^+\otimes M_{\h_2}(1)^+$-submodules $M^1\otimes M^2$ 
of $M$ and $N^1\otimes N^2$ of $N$, where 
$M^i$ and $N^i$ are irreducible $M_{\h_i}(1)^+$-modules for $i=1,2$.
Decompose $L$ as a direct sum of irreducible
$M_{\h_1}(1)^+\otimes M_{\h_2}(1)^+$-modules:
\[
L\cong\bigoplus_{j}L^1_j\otimes L^2_j,
\]
where $L^i_j$ are irreducible  $M_{\h_i}(1)^+$-modules for $i=1,2$.
By Proposition \ref{cfusion-rule-inequality} and Theorem \ref{tensortheorem} we have
\begin{align}\label{maininequality}
\begin{split}
\dim I_{\Free{+}}\fusion{M}{N}{L}&
\leq \dim I_{M_{\h_1}(1)^+\otimes M_{\h_2}(1)^+}
\fusion{M^1\otimes M^2}{N^1\otimes N^2}{L}\\
&=\sum_{j}\dim I_{M_{\h_1}(1)^+}\fusion{M^1}{N^1}{L^1_j}\cdot
\dim I_{M_{\h_2}(1)^+}\fusion{M^2}{N^2}{L^2_j}.
\end{split}
\end{align}

We take suitable irreducible $M_{\h_1}(1)^+\otimes M_{\h_2}(1)^+$-modules 
$M^1\otimes M^2$ and $N^1\otimes N^2$ from $M$ and $N$ respectively, 
and consider inequality \eqref{maininequality}. 
{}From inductive hypothesis, all the summands in the right hand side of 
\eqref{maininequality} are less than or equal to $1$. 
Furthermore, using Theorem \ref{yyyy} for $M_{\h_1}(1)^+$ 
we see that at most one of summands in the right hand side of 
\eqref{maininequality} is possibly nonzero. 
For example, in the case $M=N=\Free{-}$ and $L=\Free{+}$, we have 
\begin{align*}
&\dim I_{\Free{+}}\fusion{\Free{-}}{\Free{-}}{\Free{+}}\\
&\quad\leq \dim I_{M_{\h_1}(1)^+\otimes M_{\h_2}(1)^+}\fusion{ M_{\h_1}(1)^+\otimes
 M_{\h_2}(1)^+}{M_{\h_1}^-\otimes M_{\h_2}^-}{\Free{+}}\\
&\quad=\dim I_{M_{\h_1}(1)^+}\fusion{\Free{+}}{\Free{-}}{\Free{+}}\cdot
\dim I_{M_{\h_2}(1)^+}\fusion{\Free{-}}{\Free{+}}{\Free{+}}\\
&\qquad+\dim I_{M_{\h_1}(1)^+}\fusion{\Free{+}}{\Free{-}}{\Free{-}}\cdot
\dim I_{M_{\h_2}(1)^+}\fusion{\Free{-}}{\Free{+}}{\Free{-}}\\
&\quad=1.
\end{align*}
Therefore, the right hand side of \eqref{maininequality} is zero or one. 
This proves the proposition.
\end{proof}

Next, we show that fusion rules of certain types for $\Free{+}$
are zero. 

\begin{lemma}\label{freeplmilemma}
The fusion rules of types $\fusion{\Free{+}}{\Free{+}}{\Free{-}}$ 
and $\fusion{\Free{-}}{\Free{-}}{\Free{-}}$ are zero.
\end{lemma}

\begin{proof}
Again we shall use induction on $d=\dim \h$. 
As it was proved in \cite{A1} in the case $\dim\h=1$,
we assume that $\dim\h\geq2$. Take $h\in\h$ such that $(h,h)\neq0$ and 
set $\h_1=\C h,\,\h_2=\h_1^{\perp}$.
Then by using \eqref{maininequality} 
for $M^{1}=M_{\h_1}(1)^+,\,M^{2}=M_{\h_2}^{\pm},\,N^{1}=M_{\h_1}^{\pm}$ and $N^{2}=M_{\h_2}^\pm$ and 
the inductive hypothesis, 
we get $I_{\Free{+}}\fusion{\Free{\pm}}{\Free{\pm}}{\Free{-}}=0$ respectively.
\end{proof}

Using a similar argument we have:

\begin{lemma}\label{freeplmilemma4}
For $\lambda\in\h-\{0\}$, the fusion rules of types 
$\fusion{\Free{\pm}}{\Free{\pm}}{\Fremo{\lambda}}$ and $\fusion{\Free{\pm}}{\Free{\mp}}{\Fremo{\lambda}}$ are zero.
\end{lemma}

We shall need the following simple result in linear
algebra:

\begin{lemma}\label{linear-algebra}
Let $\h$ be a (nonzero) finite-dimensional vector space over $\C$ 
equipped with a nondegenerate symmetric bilinear form
$(\cdot,\cdot)$. Let $S$ be a finite set of nonzero vectors in $\h$.
Then there exists a one-dimensional vector subspace $\h_{1}$ of
$\h$ such that $(\cdot,\cdot)$ is nondegenerate on $\h_{1}$ and
such that $u_{1}\ne 0$ for any $u\in S$, where
$u_{1}$ denotes the orthogonal projection of $u$ into $\h_{1}$.
In particular, for $\lambda,\mu,\nu\in \h$, if the triple
$(\lambda,\mu,\nu)$ is not admissible then 
there exists a one-dimensional vector subspace $\h_{1}$ of
$\h$ such that $(\cdot,\cdot)$ is nondegenerate on $\h_{1}$ and
such that the triple $(\lambda_{1},\mu_{1},\nu_{1})$ is not admissible.
\end{lemma}

\begin{proof} Let $h_{1},\dots, h_{d}$ be an orthonormal basis of $\h$.
Then the bilinear form $(\cdot,\cdot)$ restricted on the $\R$-subspace 
$E=\R h_{1}\oplus \cdots \oplus \R h_{d}$ is positive definite.
For any $u\in \h$, we consider $u$ as a linear
functional on $\h$ through the bilinear form on $\h$. 
If $u\neq 0$, we have $(u,h_{i})\neq 0$ for some $1\leq i\leq d$, so that
$\ker u \cap E$ is a proper $\R$-subspace of $E$.
By a well known fact in linear algebra we have
\[
E\ne \cup_{u\in S}(\ker u\cap E).
\]
Take $h\in E-\cup_{u\in S}(\ker u\cap E)$ and set $\h_{1}=\C h$.
We have $(u,h)\neq 0$ for all $u\in S$.
Then $\h_{1}$ meets our need.

For $\lambda,\mu,\nu\in \h$, set
\[
S=\{a\lambda+b\mu +c\nu\,|\, a,b,c\in \{1,-1\}\}.
\]
We see that the triple $(\lambda,\mu,\nu)$ is not admissible 
if and only if $S$ consists of nonzero vectors.
Then the particular assertion follows immediately.
\end{proof}

Next we prove the following lemma:

\begin{lemma}\label{freeplmilemma2}
(1) For any $\lambda,\mu\in\h-\{0\}$, the fusion rules of types
$\fusion{\Free{\pm}}{\Fremo{\lambda}}{\Fremo{\mu}}$ are zero 
if $(\lambda,\mu,0)$ is not an admissible triple.\\
(2) Let $\lambda,\mu,\nu\in\h-\{0\}$ such that $(\lambda,\mu,\nu)$ 
is not an admissible triple. Then the fusion rule of type 
$\fusion{\Fremo{\lambda}}{\Fremo{\mu}}{\Fremo{\nu}}$ 
for $\Free{+}$ is zero. 
\end{lemma}

\begin{proof}
We also use induction on $\dim\h$. 
As it has been proved (Theorem \ref{yyyy}) in the case $\dim\h=1$,
we assume that $\dim\h\geq2$.
Since $(\lambda,\mu,0)$ is not an admissible triple,
in view of Lemma \ref{linear-algebra}, there exists 
an orthogonal decomposition $\h =\h_{1}\oplus \h_{2}$ such that
$\dim \h_{1}=1$ and $(\lambda_{1},\mu_{1},0)$ is 
not an admissible triple.
Using \eqref{maininequality} and the initial case, we obtain
\begin{align*}
&\dim I_{\Free{+}}\fusion{\Free{\pm}}{\Fremo{\lambda}}{\Fremo{\mu}}\\
&\quad\leq \dim I_{M_{\h_1}(1)^+\otimes M_{\h_2}(1)^+}
\fusion{M_{\h_1}(1)^{\pm}\otimes M_{\h_2}(1)^+}{M_{\h_{1}}(1,\lambda_{1}) 
\otimes M_{\h_{2}}(1,\lambda_{2})}{M_{\h_{1}}(1,\mu_{1})
\otimes M_{\h_{2}}(1,\mu_{2})}\\
&\quad\leq\dim I_{M_{\h_1}(1)^{+}}
\fusion{\Free{\pm}}{\Fremo{\lambda_{1}}}{\Fremo{\mu_{1}}}
\cdot\dim I_{M_{\h_2}(1)^+}
\fusion{\Free{+}}{\Fremo{\lambda_{2}}}{\Fremo{\mu_{2}}}\\
&\quad=0,
\end{align*}
proving the assertion (1). 
{}From this proof the assertion (2) is
also clear.
\end{proof}

We also have:

\begin{lemma}\label{freeplmilemma3}
The fusion rules of types
$\fusion{\Free{-}}{\Free{\pm}}{\Fretw{\pm}},
\,\fusion{\Free{-}}{\Free{\pm}}{\Fretw{\mp}}$
and 
$\fusion{\Free{-}}{\Fretw{\pm}}{\Fretw{\pm}}$ are zero.
\end{lemma}

\begin{proof}
We shall also use induction on $\dim\h$. 
As it was proved in Theorem \ref{yyyy} for rank one, we assume
that $\dim \h>1$.
As we have done before, we decompose $\h=\h_1\oplus \h_2$ (an
orthogonal sum) with $\dim \h_{1}=1$.
For any $\gamma\in \h$, $\gamma$ is decomposed as $\gamma_1+\gamma_2$ 
with $\gamma_{i}\in\h_{i}$ for $i=1,\,2$. 
Using the decomposition \eqref{decomp4}, 
the inequality \eqref{maininequality} and inductive hypothesis, we have 
\begin{align*}
&\dim I_{\Free{+}}\fusion{\Free{-}}{\Fretw{\pm}}{\Fretw{\pm}}\\
&\quad\leq\dim I_{M_{\h_1}(1)^+\otimes M_{\h_2}(1)^+}
\fusion{M_{\h_1}(1)^+\otimes M_{\h_2}(1)^-}{M_{\h_1}(1)(\theta)^+
\otimes M_{\h_{2}}(1)(\theta)^\pm}{\Fretw{\pm}}\\
&\quad\leq\dim I_{M_{\h_1}(1)^+\otimes M_{\h_2}(1)^+}
\fusion{M_{\h_1}(1)^+\otimes M_{\h_2}(1)^-}{M_{\h_{1}}(1)(\theta)^+
\otimes M_{\h_{2}}(1)(\theta)^\pm}{M_{\h_{1}}(1)(\theta)^+\otimes M_{\h_{2}}(1)(\theta)^\pm}\\
&\qquad+\dim I_{M_{\h_1}(1)^+\otimes M_{\h_2}(1)^+}
\fusion{M_{\h_1}(1)^+\otimes M_{\h_2}(1)^-}{M_{\h_{1}}(1)(\theta)^+
\otimes M_{\h_{2}}(1)(\theta)^\pm}{M_{\h_{1}}(1)(\theta)^-\otimes M_{\h_{2}}(1)(\theta)^\mp}\\
&\quad=\dim I_{M_{\h_1}(1)^+}
\fusion{\Free{+}}{\Fretw{+}}{\Fretw{+}}\cdot\dim I_{M_{\h_2}(1)^+}
\fusion{\Free{-}}{\Fretw{\pm}}{\Fretw{\pm}}\\
&\qquad+\dim I_{M_{\h_1}(1)^+}
\fusion{\Free{+}}{\Fretw{+}}{\Fretw{-}}\cdot\dim I_{M_{\h_2}(1)^+}
\fusion{\Free{-}}{\Fretw{\pm}}{\Fretw{\mp}}\\
&\quad=0,
\end{align*}
respectively. Similarly, the fusion rules of types
$\fusion{\Free{-}}{\Free{\pm}}{\Fretw{\pm}},\,
\fusion{\Free{-}}{\Free{\pm}}{\Fretw{\mp}}$ 
are also zero.
\end{proof}

Now we put everything together to prove Theorem \ref{ehigher-rank-m1}.
\begin{proof}
By Propositions \ref{fusionfreeinq}, \ref{contfreeboson} and \ref{duality},
 all the fusion rules among irreducible $\Free{+}$-modules are either 
$0$ or $1$ and 
are stable under the permutation of modules. 

We see that the fusion rule of arbitrary type for $\Free{+}$ coincides 
with one of those 
in Lemmas \ref{freeplmilemma}--\ref{freeplmilemma3} or 
Propositions \ref{nonzero1}--\ref{twistfree} 
after permuting irreducible modules. 
Furthermore, we can show that any type of fusion rule 
indicated in (i)--(v) of Theorem \ref{ehigher-rank-m1} 
agrees with one of that in Propositions \ref{nonzero1}--\ref{twistfree} 
by permuting irreducible modules.
This completes the proof.
\end{proof}

\section{Fusion rules for vertex operator algebra $\charge{+}$}
\subsection{Main theorem}
In this section we state the main result on the fusion rules for irreducible
$\charge{+}$-modules.  To do this we need to introduce a few notations.
First, recall the commutator map $c(\,\cdot\,,\,\cdot\,)$ of
$\hat{L^{\circ}}$. This defines an alternating $\Z$-bilinear form
$c_{0}:L^{\circ}\times L^{\circ}\to\Z/q\Z$ by the property
$c(a,b)=\kappa_{q}^{c_{0}(\bar{a},\bar{b})}$ for $a,b\in
\hat{L}^{\circ}$.  For $\lambda,\mu\in L^{\circ}$, we set
\begin{eqnarray}
\pi_{\lambda,\mu}=e^{(\lambda,\mu)\pi i}\w_{q}^{c_{0}(\mu,\lambda)}.
\end{eqnarray}
Note that $\pi_{\lambda,\alpha}=\pm 1$ 
for any $\alpha\in L$ if $2\lambda\in L$.
Next for a central character $\chi$  
of $\hat{L}/K$ with $\chi(\kappa)=-1$ and $\lambda\in L^{\circ}$ with
$2\lambda\in L$ we set
\begin{align}\label{dddddd}
c_{\chi}(\lambda)
=(-1)^{(\lambda,2\lambda)}\epsilon(\lambda,2\lambda)\chi(e_{2\lambda}).
\end{align}
It is easy to see that $c_{\chi}(\lambda)=\pm1$.  
For any $\lambda\in L^{\circ}$ and a central character 
$\chi$ of $\hat{L}/K$, let $\chi^{(\lambda)}$ be the central character defined by 
$\chi^{(\lambda)}(a)=(-1)^{(\bar{a},\lambda)}\chi(a)$. 
We set $T^{(\lambda)}_{\chi}=T_{\chi^{(\lambda)}}$. 
We call a triple $(\lambda,\,\mu,\,\nu)$ 
an \textit{admissible triple modulo $L$} 
if $p\lambda+q\mu+r\nu\in L$ for some $p,\,q,\,r\in\{\pm1\}$.

\begin{theorem}\label{fusioncharge} 
Let $L$ be a positive-definite even lattice.  
For any irreducible $\charge{+}$-modules $M^i\,(i=1,2,3)$, 
the fusion rule of type 
$\fusion{M^1}{M^2}{M^3}$ is either $0$ or $1$. 
The fusion rule of type $\fusion{M^1}{M^2}{M^3}$ is $1$ 
if and only if $M^i\,(i=1,2,3)$ satisfy the following conditions;

\noindent
{\rm(i)} $M^1=\charlam{\lambda}$ for $\lambda\in L^{\circ}$ 
such that $2\lambda\notin L$ and $(M^2,\,M^3)$ is one of the following pairs:
\begin{itemize}
\item[] $(\charlam{\mu},\charlam{\nu})$ for $\mu,\nu\in L^{\circ}$ 
such that $2\mu,\,2\nu\notin L$ and $(\lambda,\,\mu,\,\nu)$ is 
an admissible triple modulo $L$,

\item[] $(\charlam{\mu}^{\pm},\charlam{\nu}),\,((\charlam{\nu})',
\,(\charlam{\mu}^{\pm})')$ 
for $\mu,\nu\in L^{\circ}$ such that $2\mu\in L$ and 
$(\lambda,\,\mu,\,\nu)$ is an admissible triple modulo $L$,
\item[] 
$(\charge{T_{\chi},\pm},\charge{T_{\chi}^{(\lambda)},\pm}),
\,(\charge{T_{\chi},\pm},\charge{T_{\chi}^{(\lambda)},\mp})$ 
for any irreducible $\hat{L}/K$-module $T_{\chi}$. 
\end{itemize}

\noindent
{\rm(ii)} $M^1=\charlam{\lambda}^+$ for $\lambda\in L^{\circ}$ 
such that $2\lambda\in L$ 
and $(M^2,\,M^3)$ is one of the following pairs:
\begin{itemize}
\item[] $(\charlam{\mu},\charlam{\nu})$ for $\mu,\nu\in L^{\circ}$ 
such that $2\mu\notin L$ and $(\lambda,\,\mu,\,\nu)$ 
is an admissible triple modulo $L$,
\item[] $(\charlam{\mu}^{\pm},\charlam{\nu}^{\pm})$ 
for $\mu,\nu\in L^{\circ}$ such that $2\mu\in
L,\,\pi_{\lambda,2\mu}=1$ 
and $(\lambda,\,\mu,\,\nu)$ is an admissible triple modulo $L$,
\item[] $(\charlam{\mu}^{\pm},\charlam{\nu}^{\mp})$ 
for $\mu,\nu\in L^{\circ}$ such that $2\mu\in
L,\,\pi_{\lambda,2\mu}=-1$ 
and $(\lambda,\,\mu,\,\nu)$ is an admissible triple modulo $L$,
\item[] $(\charge{T_{\chi},\pm},\charge{T_{\chi}^{(\lambda)},\pm}),
\,((\charge{T_{\chi}^{(\lambda)},\pm})',(\charge{T_{\chi},\pm})')$ 
for any irreducible $\hat{L}/K$-module $T_{\chi}$ such that 
$c_{\chi}(\lambda)=1$,
\item[] $(\charge{T_{\chi},\pm},\charge{T_{\chi}^{(\lambda)},\mp}),
\,((\charge{T_{\chi}^{(\lambda)},\pm})',(\charge{T_{\chi},\mp})')$ 
for any irreducible $\hat{L}/K$-module $T_{\chi}$ such that 
$c_{\chi}(\lambda)=-1$.
\end{itemize} 

\noindent
{\rm(iii)} $M^1=\charlam{\lambda}^-$ for $\lambda\in L^{\circ}$ 
such that $2\lambda\in L$  and $(M^2,\,M^3)$ is one of the following pairs:
\begin{itemize}
\item[] $(\charlam{\mu},\charlam{\nu})$ for $\mu,\nu\in L^{\circ}$ 
such that $2\mu\notin L$ and $(\lambda,\,\mu,\,\nu)$ 
is an admissible triple modulo $L$,
\item[] $(\charlam{\mu}^{\pm},\charlam{\nu}^{\mp})$ 
for $\mu,\nu\in L^{\circ}$ such that $2\mu\in
L,\,\pi_{\lambda,2\mu}=1$ 
and $(\lambda,\,\mu,\,\nu)$ is an admissible triple modulo $L$,
\item[] $(\charlam{\mu}^{\pm},\charlam{\nu}^{\pm})$ 
for $\mu,\nu\in L^{\circ}$ such that $2\mu\in
L,\,\pi_{\lambda,2\mu}=-1$ 
and $(\lambda,\,\mu,\,\nu)$ is an admissible triple modulo $L$,
\item[] $(\charge{T_{\chi},\pm},\charge{T_{\chi}^{(\lambda)},\mp}),
\,((\charge{T_{\chi}^{(\lambda)},\mp})',(\charge{T_{\chi},\pm})')$ 
for any irreducible $\hat{L}/K$-module $T_{\chi}$ such that 
$c_{\chi}(\lambda)=1$,
\item[] $(\charge{T_{\chi},\pm},\charge{T_{\chi}^{(\lambda)},\pm}),
\,((\charge{T_{\chi}^{(\lambda)},\pm})',(\charge{T_{\chi},\pm})')$ 
for any irreducible $\hat{L}/K$-module $T_{\chi}$ such that 
$c_{\chi}(\lambda)=-1$.
\end{itemize} 
 
\noindent
{\rm(iv)} $M^1=\charge{T_{\chi},+}$ for an irreducible 
$\hat{L}/K$-module $T_{\chi}$ and $(M^2,\,M^3)$ is one of the following pairs:
\begin{itemize}
\item[] $(\charlam{\lambda},\charge{T_{\chi}^{(\lambda)},\pm}),
\,((\charge{T_{\chi}^{(\lambda)},\pm})',(\charlam{\lambda})')$ 
for $\lambda\in L^{\circ}$ such that $2\lambda\notin L$,
\item[] $(\charlam{\lambda}^{\pm},\charge{T_{\chi}^{(\lambda)},\pm}),
\,((\charge{T_{\chi}^{(\lambda)},\pm})',(\charlam{\lambda}^{\pm})')$ 
for $\lambda\in L^{\circ}$ such that $2\lambda\in L$ and that 
$c_{\chi}(\lambda)=1$,
\item[] $(\charlam{\lambda}^{\pm},\charge{T_{\chi}^{(\lambda)},\mp}),
\,((\charge{T_{\chi}^{(\lambda)},\mp})',(\charlam{\lambda}^{\pm})')$ 
for $\lambda\in L^{\circ}$ such that $2\lambda\in L$ and that 
$c_{\chi}(\lambda)=-1$.
\end{itemize} 

\noindent
{\rm(v)} $M^1=\charge{T_{\chi},-}$ for an irreducible 
$\hat{L}/K$-module $T_{\chi}$ and $(M^2,\,M^3)$ is one of the following pairs:
\begin{itemize}
\item[] $(\charlam{\lambda},\charge{T_{\chi}^{(\lambda)},\pm}),\,((\charge{T_{\chi}^{(\lambda)},\pm})',(\charlam{\lambda})')$ for $\lambda\in L$ such that $2\lambda\notin L^{\circ}$,
\item[] $(\charlam{\lambda}^{\pm},\charge{T_{\chi}^{(\lambda)},\mp}),\,((\charge{T_{\chi}^{(\lambda)},\pm})',(\charlam{\lambda}^{\mp})')$ for $\lambda\in L^{\circ}$ such that $2\lambda\in L$ and that $c_{\chi}(\lambda)=1$,
\item[] $(\charlam{\lambda}^{\pm},\charge{T_{\chi}^{(\lambda)},\pm}),\,((\charge{T_{\chi}^{(\lambda)},\mp})',(\charlam{\lambda}^{\mp})')$ for $\lambda\in L^{\circ}$ such that $2\lambda\in L$ and that $c_{\chi}(\lambda)=-1$.
\end{itemize} 
\end{theorem}

\begin{remark}\label{wwwww}
In the case that the rank of $L$ is one, 
Theorem \ref{fusioncharge} was previously proved in \cite{A2}. 
%But it is not easy to compare Theorem \ref{fusioncharge} 
%with the results in \cite{A2}. 
\end{remark}

We will give a proof of this theorem in Sections \ref{FRAM} and \ref{FRIM}, 
where we deal with the fusion rules for irreducible modules 
of untwisted types and for those of twisted types respectively.

\subsection{Fusion rules among modules of untwisted types}\label{FRAM}
In this section we determine the fusion rules for the irreducible
$\charge{+}$-modules of untwisted types.
We first prove that the fusion rules of certain types 
for irreducible $\charge{+}$-modules of untwisted types are
nonzero by giving nonzero intertwining operators. 
Such intertwining operators come from 
intertwining operators constructed in \cite{DL1}
for irreducible $\charge{}$-modules. 

We recall a construction of intertwining operators 
for irreducible $\charge{}$-modules following \cite{DL1}.
Let $Y(\,\cdot\,,z):V_{L^{\circ}}\to(\End V_{L^{\circ}})\{z\}$ 
be the linear map as in \eqref{ffff} (with $P=L^{\circ}$).
However $Y(\,\cdot\,,z)$ satisfies the $L(-1)$-derivative property, 
the identity \eqref{iiii} implies that $Y(\,\cdot\,,z)$ 
does not give intertwining operators among irreducible $\charge{}$-modules. 
We attach an extra factor to $Y(\,\cdot\,,z)$ to get intertwining operators.  
Let $\lambda\in L^{\circ}$.
We define a linear map $\pi^{(\lambda)}\in\End V_{L^{\circ}}$ 
which acts on $\Fremo{\mu}\,(\mu\in L^{\circ})$ 
as the scalar 
$\pi_{\lambda,\mu}(=e^{(\lambda,\mu)\pi i}\w_{q}^{c_{0}(\mu,\lambda)})$ 
and we then define a linear map 
$\mathcal{Y}_{\lambda}(\,\cdot\,,z):\charlam{\lambda}\to(\End V_{L^{\circ}})\{z\}$ by 
\begin{eqnarray}\label{eYlambda}
\mathcal{Y}_{\lambda}(u,z)v=Y(u,z)\pi^{(\lambda)}(v)
\end{eqnarray}
for any $u\in\Fremo{\lambda}$ and $v\in\Fremo{\mu}$. 
Then the restriction of $\mathcal{Y}_{\lambda}(\,\cdot\,,z)$ 
gives rise to a nonzero intertwining operator of type 
$\fusion{\charlam{\lambda}}{\charlam{\mu}}{\charlam{\lambda+\mu}}$.

Let $\lambda,\mu,\gamma\in L^{\circ}$ (the dual lattice of $L$).
It was proved in \cite[Proposition 12.8]{DL1} that
$I_{\charge{}}\fusion{\charlam{\lambda}}{\charlam{\mu}}
{V_{\gamma+L}}$
is nonzero if and only if $\gamma-\lambda-\mu\in L$ and that
$I_{\charge{}}\fusion{\charlam{\lambda}}{\charlam{\mu}}
{\charlam{\lambda+\mu}}$ is one dimensional.
Thus the fusion of type
$\fusion{\charlam{\lambda}}{\charlam{\mu}}{\charlam{\lambda+\mu}}$ for
$\charge{+}$ is nonzero. 
Using a $\charge{+}$-module isomorphism between $\charlam{\mu}$ and 
$\charlam{-\mu}$, we see that the fusion rule of type 
$\fusion{\charlam{\lambda}}{\charlam{\mu}}{\charlam{\lambda-\mu}}$ 
is also nonzero. Furthermore we have:

\begin{proposition}\label{untw}
For any $\lambda,\,\mu,\,\nu\in L^{\circ}$, the fusion rule of type 
$\fusion{\charlam{\lambda}}{\charlam{\mu}}{\charlam{\nu}}$ 
for $\charge{+}$ is nonzero 
if and only if $(\lambda,\mu,\nu)$ is an admissible triple modulo $L$.
\end{proposition} 

\begin{proof}
Let $(\lambda,\mu,\nu)$ be an admissible triple modulo $L$.
Then $\charlam{\nu}$ is isomorphic to $\charlam{\lambda+\mu}$ or 
$\charlam{\lambda-\nu}$ as a $V_{L}^{+}$-module.
Hence the fusion rule of type 
$\fusion{\charlam{\lambda}}{\charlam{\mu}}{\charlam{\nu}}$ is nonzero.

Conversely, let us assume that the fusion rule of type 
$\fusion{\charlam{\lambda}}{\charlam{\mu}}{\charlam{\nu}}$ is nonzero. 
We take $\lambda,\mu$ to be nonzero if necessary. 
Note that  for any $\gamma\in L^{\circ}$,
$\charlam{\gamma}\cong\bigoplus_{\alpha\in L}\Fremo{\gamma+\alpha}$ 
as an $\Free{+}$-module.
Since $\charlam{\lambda}$ and $\charlam{\mu}$ contain
irreducible $\Free{+}$-modules 
$\Fremo{\lambda}$ and $\Fremo{\mu}$, respectively, 
by Proposition \ref{cfusion-rule-inequality}, the fusion rule of type 
$\fusion{\Fremo{\lambda}}{\Fremo{\mu}}{\charlam{\nu}}$ 
for $\Free{+}$ is nonzero.
By Theorem \ref{ehigher-rank-m1}, $\charlam{\nu}$ must contain 
an irreducible $\Free{+}$-submodule 
isomorphic to $\Fremo{\lambda+\mu}$ or $\Fremo{\lambda-\mu}$. 
Then  $\lambda+\mu\in \nu+L$, or $-\nu+L$, or $\lambda-\mu\in \nu+L$,
or $-\nu+L$. 
This
shows that $(\lambda,\mu,\nu)$ is an admissible triple modulo $L$.
\end{proof}

Furthermore, if $2\lambda\in L$,
by Proposition \ref{cfusion-rule-inequality}, 
we see that the fusion rules of types 
$\fusion{\charlam{\lambda}^{\pm}}{\charlam{\mu}}{\charlam{\lambda+\mu}}$
are not zero. 
Similarly, the fusion rules of types 
$\fusion{\charlam{\lambda}^{\pm}}{\charlam{\mu}}{\charlam{\lambda-\mu}}$ 
are also nonzero. 
Clearly, if one of the fusion rules of types 
$\fusion{\charlam{\lambda}^{\pm}}{\charlam{\mu}}{\charlam{\lambda-\mu}}$ 
is nonzero, the fusion rule of type 
$\fusion{\charlam{\lambda}}{\charlam{\mu}}{\charlam{\nu}}$ 
for $\charge{+}$ is nonzero.
In view of Proposition \ref{untw} we immediately have:

\begin{proposition}\label{qqqq}
For any $\lambda,\mu,\nu\in L^{\circ}$ with $2\lambda\in L$, 
the fusion rules of types 
$\fusion{\charlam{\lambda}^{\pm}}{\charlam{\mu}}{\charlam{\nu}}$ are nonzero 
if and only if $(\lambda,\mu,\nu)$ is an admissible triple modulo $L$.
\end{proposition}

We next prove the following result:

\begin{proposition}\label{bbbb}
Let $M^{1},\,M^{2}$ and $M^{3}$ be irreducible $\charge{+}$-modules 
of untwisted types. 
Suppose that one of $M^{i}\,(i=1,2,3)$ is isomorphic to $\charlam{\lambda}$ 
for $\lambda\in L^{\circ}$ with $2\lambda\notin L$ or $\charge{\pm}$. 
Then the fusion rule of type $\fusion{M^{1}}{M^{2}}{M^{3}}$ is 
either $0$ or $1$.
\end{proposition}

\begin{proof}
%Since $M(1,\alpha)$ and $M(1,-\alpha)$ are isomorphic irreducible $\Free{+}$-modules for any $0\ne \alpha\in L$, we see that 
%\[
%\charge{\pm}\cong\Free{\pm}\oplus\bigoplus_{\alpha\in L-\{0\}}\Fremo{\alpha}
%\]
%as an $\Free{+}$-module,
%where $(\Fremo{\alpha}\oplus \Fremo{-\alpha})^{\pm}\cong\Fremo{\alpha}$.
For $\lambda\in L^{\circ}$, 
the $\charge{+}$-module $\charlam{\lambda}$ is decomposed 
into a direct sum of irreducible $\Free{+}$-modules
as 
\[
\charlam{\lambda}\cong\bigoplus_{\alpha\in L}\Fremo{\lambda+\alpha}.
\]
Moreover, if $2\lambda\in L$, we can take a subset $S_{\lambda}\subset \lambda+L$ 
so that $S_{\lambda}\cap(-S_{\lambda})=\emptyset$ and 
$S_{\lambda}\cup(-S_{\lambda})=\lambda+L\,(=L-\{0\}$ if $\lambda\in L)$, and we have 
\begin{align*}
&\charge{\pm}\cong\Free{\pm}\oplus\bigoplus_{\mu\in S_{\lambda}}\Fremo{\mu}\quad\hbox{if $\lambda\in L$},\\
&\charlam{\lambda}^{+}\cong\charlam{\lambda}^{-}\cong\bigoplus_{\mu\in S_{\lambda}}\Fremo{\mu}\quad\hbox{if $\lambda\notin L$}
\end{align*} 
as $\Free{+}$-modules.
Therefore, the multiplicity of any irreducible $\Free{+}$-module 
in any irreducible $\charge{+}$-module of untwisted type is at most
one and any irreducible $\charge{+}$-module of untwisted type
contains an irreducible $\Free{+}$-submodule isomorphic to
$M(1,\beta)$ with $0\ne \beta\in L^{\circ}$. 

Let $M^{1},\,M^{2}$ and $M^{3}$ be irreducible $\charge{+}$-modules 
of untwisted type. 
{}From the previous paragraph, each $M^{i}$
contains $\Fremo{\lambda_{i}}$
as an irreducible $\Free{+}$-submodule
for some nonzero $\lambda_{i}\in L^{\circ}$ for $i=1,2,3$.
In view of Proposition \ref{cfusion-rule-inequality}, we see 
that the fusion rule of type $\fusion{M^{1}}{M^{2}}{M^{3}}$
for $V_{L}^{+}$-modules is not bigger than that of
type $\fusion{M(1,\lambda_{1})}{M(1,\lambda_{2})}{M^3}$
for $M(1)^{+}$-modules. 
Assume that the fusion rule of type $\fusion{M^{1}}{M^{2}}{M^{3}}$
for $V_{L}^{+}$-modules is not zero. 
{}From Theorem \ref{yyyy} we must have
$a\lambda_{1}+b\lambda_{2}\in \lambda_{3}+L$ for some $a,b\in \{1,-1\}$.
That is, $(\lambda_{1},\lambda_{2},\lambda_{3})$ is
an admissible triple modulo $L$.

By using Propositions \ref{duality} and \ref{eeee} we may assume that 
$M^{1}$ is isomorphic to one of the irreducible modules
$V_{L}^{\pm}$ and $V_{\lambda+L}$ for $\lambda\in L^{\circ}$ with
$2\lambda\notin L$. We divide the proof in the following three cases.

Case 1: $M^{1}=V_{L}^{+}$.
{}From Remark 2.9 of \cite{liform} we have that for any vertex operator
algebra $V$ and for any $V$-modules $W$ and $M$, the fusion rule of
type $\fusion{V}{W}{M}$ equals $\dim \Hom_{V} (W,M)$. It follows from
Schur lemma (see [FHL]) that the fusion rule of
type  $\fusion{V}{W}{M}$ for irreducible $V$-modules $W$ and $M$ 
is either $0$ or $1$.

Case 2: $M^{1}=V_{L}^{-}$. From Theorem \ref{yyyy} (ii),
for any irreducible $M(1)^{+}$-module $W$
the fusion rule of type $\fusion{M(1)^{-}}{M(1,\lambda_{2})}{W}$
for $M(1)^{+}$-modules is $1$ if 
$W\cong M(1,-\lambda_{2})$ and it is zero otherwise. 
We also know that the multiplicity of $M(1,-\lambda_{2})$ in
$M^{3}$ is one. Thus the fusion rule of the type
$\fusion{M^{1}}{M^{2}}{M^{3}}$ is at most $1$. 

Case 3: $M^{1}=V_{\lambda+L}$ for $\lambda\in L^{\circ}$ with
$2\lambda\notin L$.  Because $(\lambda_{1},\lambda_{2},\lambda_{3})$
is an admissible triple modulo $L$, we have that either
$2\lambda_{2}\notin L$ or $2\lambda_{3}\notin L$.  By using
Propositions \ref{duality} and \ref{eeee} we may assume that
$2\lambda_{3}\notin L$.  This implies that $\charlam{\lambda_{3}}$
contains either $\Fremo{\lambda_{1}+\lambda_{2}}$ or
$\Fremo{\lambda_{1}-\lambda_{2}}$, as an $\Free{+}$-submodule with
multiplicity one. In view of Theorem \ref{yyyy} and Proposition
\ref{cfusion-rule-inequality}, the fusion rule of type
$\fusion{M^{1}}{M^{2}}{M^{3}}$ is either $0$ or $1$.
\end{proof}

Let $\lambda,\mu\in L^{\circ}$ such that $2\lambda,\,2\mu\in L$.   
Then we see that $\mathcal{Y}_{\lambda}(\,\cdot\,,z)$ gives rise to
a nonzero intertwining operator of type 
$\fusion{\charlam{\lambda}^{\epsilon_{1}}}{\charlam{\mu}^{\epsilon_{2}}}{\charlam{\lambda+\mu}}$ 
for any $\epsilon_{1},\,\epsilon_{2}\in\{\pm\}$. 
We consider the conjugation $\theta\mathcal{Y}_{\lambda}(\,\cdot\,,z)\theta^{-1}$. 
By definition, we have for any $\beta\in L$ and $v\in\Fremo{\mu+\beta}$,
\begin{align*}
\theta(\pi^{(\lambda)}(\theta^{-1}(v)))
&=e^{(\lambda,-\mu-\beta)\pi i}\w_{q}^{c_{0}(-\mu-\beta,\lambda)}v\\
&=e^{(\lambda,-2\mu-2\beta)\pi i}\w_{q}^{c_{0}(-2\mu-2\beta,\lambda)}
e^{(\lambda,\mu+\beta)\pi i}\w_{q}^{c_{0}(\mu+\beta,\lambda)}v\\
&=e^{(\lambda,-2\beta)\pi i}\w_{q}^{c_{0}(-2\beta,\lambda)}
e^{(\lambda,-2\mu)\pi i}\w_{q}^{c_{0}(-2\mu,\lambda)}\pi^{(\lambda)}(v)\\
&=\pi_{\lambda,-2\mu}\pi^{(\lambda)}(v)\\
&=\pi_{\lambda,2\mu}\pi^{(\lambda)}(v),
\end{align*} 
noticing that
\[
e^{(\lambda,-2\beta)\pi i}=1,\;\;\;\;\w_{q}^{c_{0}(-2\beta,\lambda)}
=\w_{q}^{c_{0}(\beta,-2\lambda)}=(-1)^{(\beta,-2\lambda)}=1.
\]
Using \eqref{oooo} and (\ref{eYlambda}), we get
\begin{align}\label{gggg}
\theta\mathcal{Y}_{\lambda}(u,z)\theta^{-1}(v)
=\pi_{\lambda,2\mu}\mathcal{Y}_{\lambda}(\theta(u),z)v.
\end{align}
Next we prove the following result:

\begin{proposition}\label{qqqqq}
Let $\lambda,\,\mu,\,\nu\in L^{\circ}$ such that $2\lambda,\,2\mu\in L$. 
(1) If $(\lambda,\mu,\nu)$ is not an admissible triple modulo $L$ 
then the fusion rule of type 
$\fusion{\charlam{\lambda}^{\epsilon_{1}}}
{\charlam{\mu}^{\epsilon_{2}}}{\charlam{\nu}^{\epsilon_{3}}}$ 
is zero for any  $\epsilon_{i}\in\{\pm\}\,(i=1,2,3)$.

(2) Let $(\lambda,\mu,\nu)$ be an admissible triple modulo $L$.
Then the fusion rules of types 
$\fusion{\charlam{\lambda}^{\pm}}{\charlam{\mu}^{\pm}}{\charlam{\nu}^{+}}$ 
and 
$\fusion{\charlam{\lambda}^{\pm}}{\charlam{\mu}^{\mp}}{\charlam{\nu}^{-}}$ 
are nonzero if and only if $\pi_{\lambda,2\mu}=1$.
The fusion rules of types 
$\fusion{\charlam{\lambda}^{\pm}}{\charlam{\mu}^{\pm}}{\charlam{\nu}^{-}}$ 
and 
$\fusion{\charlam{\lambda}^{\pm}}{\charlam{\mu}^{\mp}}{\charlam{\nu}^{+}}$ 
are nonzero if and only if $\pi_{\lambda,2\mu}=-1$.
Furthermore, the fusion rules of type 
$\fusion{\charlam{\lambda}^{\epsilon_{1}}}
{\charlam{\mu}^{\epsilon_{2}}}{\charlam{\nu}^{\epsilon_{3}}}$ 
is either $0$ or $1$ for $\epsilon_{i}\in\{\pm\}$. 
\end{proposition}

\begin{proof}
The assertion (1) follows immediately from Propositions \ref{untw}.

We now prove  (2).
By \eqref{gggg} we see that $\mathcal{Y}_{\lambda}(\,\cdot\,,z)$ gives
nonzero intertwining operators of types
$\fusion{\charlam{\lambda}^{\pm}}{\charlam{\mu}^{\pm}}{\charlam{\nu}^{+}}$
$(\fusion{\charlam{\lambda}^{\pm}}{\charlam{\mu}^{\pm}}{\charlam{\nu}^{-}}$
resp.) if $\pi_{\lambda,2\mu}=1$ ($\pi_{\lambda,2\mu}=-1$ reps.), so that
the corresponding fusion rules are nonzero.
 
It is enough to prove that the fusion rule of type 
$\fusion{\charlam{\lambda}^{\epsilon_{1}}}{\charlam{\mu}^{\epsilon_{2}}}{\charlam{\lambda+\mu}}$ 
for $\charge{+}$ is one for any $\epsilon_{1},\epsilon_{2}\in\{\pm\}$. 
We shall demonstrate the proof only for $\epsilon_{1}=\epsilon_{2}=+$. 
The other cases can be proved similarly.

As in the proof of Proposition \ref{bbbb}, for any nonzero 
$\nu\in L^{\circ}$ with $2\nu\in L$, 
we take a subset $S_{\nu}\subset\nu+L$ such that 
$S_{\nu}\cap(-S_{\nu})=\emptyset$ 
and $S_{\nu}\cup(-S_{\nu})=\nu+L\,(L-\{0\}$ if $\nu\in L)$. 
We may assume that $\nu,3\nu\in S_{\nu}$. 
Then we have an $\Free{+}$-isomorphism 
$\phi:\charlam{\nu}^{+}\to\bigoplus_{\gamma\in S_{\nu}}\Fremo{\gamma}\,(\,\charge{+}\to\Free{+}\oplus\bigoplus_{\gamma\in S_{\nu}}\Fremo{\gamma}$ if $\nu\in L)$ 
such that $\phi(u+\theta(u))=u$ for any $\gamma\in S_{\nu}$ and $u\in\Fremo{\gamma}$. 
Set 
\[
\charlam{\nu}^{+}[\gamma]=\Free{+}\otimes (e^{\gamma}+e^{-\gamma})\oplus \Free{-}\otimes (e^{\gamma}-e^{-\gamma})\subset \charlam{\nu}
\]
for $\gamma\in\nu+L$. 
Then $\phi$ gives an $\Free{+}$-isomorphism from  $\charlam{\nu}^{+}[\gamma]$ to $\Fremo{\gamma}$. 

Let $\gamma\in\lambda+L$ and $\delta\in\mu+L$. 
By Theorem \ref{yyyy}, the  vector space $I_{\Free{+}}\fusion{\charlam{\lambda}^{+}[\gamma]}{\charlam{\mu}^{+}[\delta]}{\charlam{\lambda+\mu}}$ is of $4$ dimension and is spanned by $\mathcal{Y}_{i}(\,\cdot\,,z)\,(i=1,2,3,4)$ defined by 
\begin{align*}
\mathcal{Y}_{1}(u,z)v&=\mathcal{Y}_{\gamma,\delta}(\phi(u),z)\phi(v),\\
\mathcal{Y}_{2}(u,z)v&=\mathcal{Y}_{\gamma,-\delta}(\phi(u),z)\theta(\phi(v)),\\
\mathcal{Y}_{3}(u,z)v&=\mathcal{Y}_{-\gamma,\delta}(\theta(\phi(u)),z)\phi(v),\\
\mathcal{Y}_{4}(u,z)v&=\mathcal{Y}_{-\gamma,-\delta}(\theta(\phi(u)),z)\theta(\phi(v))
\end{align*}
for any $u\in \charlam{\lambda}^{+}[\gamma]$ and $v\in\charlam{\mu}^{+}[\delta]$.

For $\alpha\in L$ we set $E^{\alpha}=e^{\alpha}+e^{-\alpha}\in V_L^+.$ 
Then for $u\in\charlam{\lambda}^{+}[\lambda],\,v\in\charlam{\mu}^{+}[\mu]$, there exists a nonnegative integer $k$ such that
\begin{align*}
&(z_{1}-z_{2})^{k}Y(E^{2\mu},z_{1})\mathcal{Y}_{1}(u,z_{2})v\\
&\quad=(z_{1}-z_{2})^{k}Y(E^{2\mu},z_{1})\mathcal{Y}_{\lambda,\mu}(\phi(u),z_{2})\phi(v)\\
&\quad=(z_{1}-z_{2})^{k}\left(\epsilon(2\mu,\lambda+\mu)\mathcal{Y}_{2\mu,\lambda+\mu}(e^{2\mu},z_{1})\mathcal{Y}_{\lambda,\mu}(\phi(u),z_{2})\phi(v)\right.\\
&\qquad\left.+\epsilon(-2\mu,\lambda+\mu)\mathcal{Y}_{-2\mu,\lambda+\mu}(e^{-2\mu},z_{1})\mathcal{Y}_{\lambda,\mu}(\phi(u),z_{2})\phi(v)\right)\\
&\quad=(z_{1}-z_{2})^{k}(-1)^{(2\mu,\lambda)}\left(\epsilon(2\mu,\lambda+\mu)\mathcal{Y}_{\lambda,3\mu}(\phi(u),z_{2})\mathcal{Y}_{2\mu,\mu}(e^{2\mu},z_{1})\phi(v)\right.\\
&\qquad\left.+\epsilon(-2\mu,\lambda+\mu)\mathcal{Y}_{\lambda,-\mu}(\phi(u),z_{2})\mathcal{Y}_{-2\mu,\mu}(e^{-2\mu},z_{1})\phi(v)\right).
\end{align*}
As well, we have 
\begin{align*}
&(z_{1}-z_{2})^{k}Y(E^{2\mu},z_{1})\mathcal{Y}_{2}(u,z_{2})v\\
&\quad=(z_{1}-z_{2})^{k}(-1)^{(2\mu,\lambda)}\left(\epsilon(2\mu,\lambda-\mu)\mathcal{Y}_{\lambda,\mu}(\phi(u),z_{2})\mathcal{Y}_{2\mu,-\mu}(e^{2\mu},z_{1})\theta(\phi(v))\right.\\
&\qquad\left.+\epsilon(-2\mu,\lambda-\mu)\mathcal{Y}_{\lambda,-3\mu}(\phi(u),z_{2})\mathcal{Y}_{-2\mu,-\mu}(e^{-2\mu},z_{1})\theta(\phi(v))\right),\\
&(z_{1}-z_{2})^{k}Y(E^{2\mu},z_{1})\mathcal{Y}_{3}(u,z_{2})v\\
&\quad=(z_{1}-z_{2})^{k}(-1)^{(2\mu,-\lambda)}\left(\epsilon(2\mu,-\lambda+\mu)\mathcal{Y}_{-\lambda,3\mu}(\theta(\phi(u)),z_{2})\mathcal{Y}_{2\mu,\mu}(e^{2\mu},z_{1})\phi(v)\right.\\
&\qquad\left.+\epsilon(-2\mu,-\lambda+\mu)\mathcal{Y}_{-\lambda,-\mu}(\theta(\phi(u)),z_{2})\mathcal{Y}_{-2\mu,\mu}(e^{-2\mu},z_{1})\phi(v)\right),\\
&(z_{1}-z_{2})^{k}Y(E^{2\mu},z_{1})\mathcal{Y}_{4}(u,z_{2})v\\
&\quad=(z_{1}-z_{2})^{k}(-1)^{(2\mu,-\lambda)}\left(\epsilon(2\mu,-\lambda-\mu)\mathcal{Y}_{-\lambda,\mu}(\theta(\phi(u)),z_{2})\mathcal{Y}_{2\mu,-\mu}(e^{2\mu},z_{1})\theta(\phi(v))\right.\\
&\qquad\left.+\epsilon(-2\mu,-\lambda-\mu)\mathcal{Y}_{-\lambda,-3\mu}(\theta(\phi(u)),z_{2})\mathcal{Y}_{-2\mu,-\mu}(e^{-2\mu},z_{1})\theta(\phi(v))\right).
\end{align*}
For simplicity, we set
\begin{align*}
A^{i,j}&=\mathcal{Y}_{(-1)^{i}\lambda,(2+(-1)^{j})\mu}(\theta^{i}(\phi(u)),z_{2})\mathcal{Y}_{2\mu,(-1)^{j}\mu}(e^{2\mu},z_{1})\theta^{j}(\phi(v)),\\
B^{i,j}&=\mathcal{Y}_{(-1)^{i}\lambda,(-2+(-1)^{j})\mu}(\theta^{i}(\phi(u)),z_{2})\mathcal{Y}_{-2\mu,(-1)^{j}\mu}(e^{-2\mu},z_{1})\theta^{j}(\phi(v))
\end{align*}
for $i=0,1$.
Then we see that $A^{i,j}\in\Fremo{(-1)^{i}\lambda+(2+(-1)^{j})\mu}\{z_{1}\}\{z_{2}\}$ and $B^{i,j}\in\Fremo{(-1)^{i}\lambda+(-2+(-1)^{j})\mu}\{z_{1}\}\{z_{2}\}$ and that $A^{i,j}$ and $B^{i,j}$ for $i,j=0,1$ are linearly independent in $\charlam{\lambda+\mu}\{z_{1}\}\{z_{2}\}$.

Now we let $\mathcal{Y}(\,\cdot\,,z)$ be an intertwining operator of type $\fusion{\charlam{\lambda}^{+}}{\charlam{\mu}^{+}}{\charlam{\lambda+\mu}}$.
Then for $\gamma\in \lambda+L,$ $\delta\in \mu +L,$ 
there are $c^{i}_{\gamma,\delta}\in\C$ such that the restriction of $\mathcal{Y}(\,\cdot\,,z)$ to $\charlam{\lambda}^{+}[\gamma]\otimes\charlam{\mu}^{+}[\delta]$ is expressed by $\mathcal{Y}(\,\cdot\,,z)=\sum_{i=1}^{4}c^{i}_{\gamma,\delta}\mathcal{Y}_{i}(\,\cdot\,,z)$. 
Thus, 
\begin{align*}
&(z_{1}-z_{2})^{k}Y(e^{2\mu},z_{1})\mathcal{Y}(u,z_{2})v\\
&=(z_{1}-z_{2})^{k}(-1)^{(2\mu,\lambda)}(\epsilon(2\mu,\lambda+\mu)c_{\lambda,\mu}^{(1)}A^{0,0}+\epsilon(-2\mu,\lambda+\mu)c_{\lambda,\mu}^{(1)}B^{0,0}\\
&\quad+\epsilon(2\mu,\lambda-\mu)c_{\lambda,\mu}^{(2)}A^{0,1}+\epsilon(-2\mu,\lambda-\mu)c_{\lambda,\mu}^{(2)}B^{0,1}\\
&\quad+\epsilon(2\mu,-\lambda+\mu)c_{\lambda,\mu}^{(3)}A^{1,0}+\epsilon(-2\mu,-\lambda+\mu)c_{\lambda,\mu}^{(3)}B^{1,0}\\
&\quad+\epsilon(2\mu,-\lambda-\mu)c_{\lambda,\mu}^{(4)}A^{1,1}\epsilon(-2\mu,-\lambda-\mu)c_{\lambda,\mu}^{(4)}B^{1,1}).
\end{align*}

Since $\mu,\,3\mu\in S_{\mu}$, we see that 
\begin{align*}
\phi(Y(E^{2\mu},z)v)&=Y(e^{2\mu},z)\phi(v)+Y(e^{2\mu},z)\theta(\phi(v))\\
&=\epsilon(2\mu,\mu)\mathcal{Y}_{2\mu,\mu}(e^{2\mu},z)\phi(v)+\epsilon(2\mu,-\mu)\mathcal{Y}_{2\mu,-\mu}(e^{2\mu},z)\theta(\phi(v)).
\end{align*} 
Thus we get 
\begin{align*}
&(z_{1}-z_{2})^{k}\mathcal{Y}(u,z_{2})Y(E^{2\mu},z_{1})v\\
&\quad=(z_{1}-z_{2})^{k}(
  c_{\lambda,3\mu}^{(1)}\epsilon(2\mu,\mu)\mathcal{Y}_{\lambda,3\mu}(\phi(u),z_{2})\mathcal{Y}_{2\mu,\mu}(e^{2\mu},z_{1})\phi(v)\\
&\qquad+c_{\lambda,3\mu}^{(2)}\epsilon(2\mu,\mu)\mathcal{Y}_{\lambda,-3\mu}(\phi(u),z_{2})\mathcal{Y}_{-2\mu,-\mu}(e^{-2\mu},z_{1})\theta(\phi(v))\\
&\qquad+c_{\lambda,3\mu}^{(3)}\epsilon(2\mu,\mu)\mathcal{Y}_{-\lambda,3\mu}(\theta(\phi(u)),z_{2})\mathcal{Y}_{2\mu,\mu}(e^{2\mu},z_{1})\phi(v)\\
&\qquad+c_{\lambda,3\mu}^{(4)}\epsilon(2\mu,\mu)\mathcal{Y}_{-\lambda,-3\mu}(\theta(\phi(u)),z_{2})\mathcal{Y}_{-2\mu,-\mu}(e^{-2\mu},z_{1})\theta(\phi(v))\\
&\qquad+c_{\lambda,\mu}^{(1)}\epsilon(2\mu,-\mu)\mathcal{Y}_{\lambda,-\mu}(\phi(u),z_{2})\mathcal{Y}_{2\mu,-\mu}(e^{2\mu},z_{1})\theta(\phi(v))\\
&\qquad+c_{\lambda,\mu}^{(2)}\epsilon(2\mu,-\mu)\mathcal{Y}_{\lambda,-\mu}(\phi(u),z_{2})\mathcal{Y}_{-2\mu,\mu}(e^{-2\mu},z_{1})\phi(v)\\
&\qquad+c_{\lambda,\mu}^{(3)}\epsilon(2\mu,-\mu)\mathcal{Y}_{-\lambda,-\mu}(\theta(\phi(u)),z_{2})\mathcal{Y}_{2\mu,-\mu}(e^{2\mu},z_{1})\theta(\phi(v))\\
&\qquad+c_{\lambda,\mu}^{(4)}\epsilon(2\mu,-\mu)\mathcal{Y}_{-\lambda,-\mu}(\theta(\phi(u)),z_{2})\mathcal{Y}_{-2\mu,\mu}(e^{-2\mu},z_{1})\phi(v))\\
&\quad=(z_{1}-z_{2})^{k}(
  c_{\lambda,3\mu}^{(1)}\epsilon(2\mu,\mu)A^{0,0}+c_{\lambda,3\mu}^{(2)}\epsilon(2\mu,\mu)B^{0,1}\\
&\qquad+c_{\lambda,3\mu}^{(3)}\epsilon(2\mu,\mu)A^{1,0}+c_{\lambda,3\mu}^{(4)}\epsilon(2\mu,\mu)B^{1,1}\\
&\qquad+c_{\lambda,\mu}^{(1)}\epsilon(2\mu,-\mu)A^{0,1}+c_{\lambda,\mu}^{(2)}\epsilon(2\mu,-\mu)B^{0,0}\\
&\qquad+c_{\lambda,\mu}^{(3)}\epsilon(2\mu,-\mu)A^{1,1}+c_{\lambda,\mu}^{(4)}\epsilon(2\mu,-\mu)B^{1,0}).
\end{align*}  
Since $\mathcal{Y}(\,\cdot\,,z)$ is an intertwining operator for $\charge{+}$, we have $(z_{1}-z_{2})^{k}Y(E^{2\mu},z)\mathcal{Y}(u,z)v=(z_{1}-z_{2})^{k}\mathcal{Y}(u,z)Y(E^{2\mu},z)v$ for sufficiently large integer $k$. 
Therefore, the linearly independence of $A^{i,j}$ and $B^{i,j}$ for $i,j=0,1$ gives the following equations:
\begin{align*}
&(-1)^{(2\mu,\lambda)}\epsilon(2\mu,\lambda+\mu)c_{\lambda,\mu}^{(1)}=c_{\lambda,3\mu}^{(1)}\epsilon(2\mu,\mu),\\
&(-1)^{(2\mu,\lambda)}\epsilon(-2\mu,\lambda+\mu)c_{\lambda,\mu}^{(1)}=c_{\lambda,\mu}^{(2)}\epsilon(2\mu,-\mu),\\
&(-1)^{(2\mu,\lambda)}\epsilon(2\mu,\lambda-\mu)c_{\lambda,\mu}^{(2)}=c_{\lambda,\mu}^{(1)}\epsilon(2\mu,-\mu),\\
&(-1)^{(2\mu,\lambda)}\epsilon(-2\mu,\lambda-\mu)c_{\lambda,\mu}^{(2)}=c_{\lambda,3\mu}^{(2)}\epsilon(2\mu,\mu),\\
&(-1)^{(2\mu,\lambda)}\epsilon(2\mu,-\lambda+\mu)c_{\lambda,\mu}^{(3)}=c_{\lambda,3\mu}^{(3)}\epsilon(2\mu,\mu),\\
&(-1)^{(2\mu,\lambda)}\epsilon(-2\mu,-\lambda+\mu)c_{\lambda,\mu}^{(3)}=c_{\lambda,\mu}^{(4)}\epsilon(2\mu,-\mu),\\
&(-1)^{(2\mu,\lambda)}\epsilon(2\mu,-\lambda-\mu)c_{\lambda,\mu}^{(4)}=c_{\lambda,\mu}^{(3)}\epsilon(2\mu,-\mu),\\
&(-1)^{(2\mu,\lambda)}\epsilon(-2\mu,-\lambda-\mu)c_{\lambda,\mu}^{(4)}=c_{\lambda,3\mu}^{(4)}\epsilon(2\mu,\mu).
\end{align*}
{}From these equations we get
\begin{align}\label{lllp}
c_{\lambda,\mu}^{(2)}
=(-1)^{(2\mu,\lambda)}\epsilon(-2\mu,\lambda)c_{\lambda,\mu}^{(1)},
\quad c_{\lambda,\mu}^{(4)}
=(-1)^{(2\mu,\lambda)}\epsilon(2\mu,\lambda)c_{\lambda,\mu}^{(3)}.
\end{align}

Next we shall apply a similar argument to the associativity
\[
(z_{0}+z_{2})^{k}Y(E^{2\lambda},z_{0}+z_{2})\mathcal{Y}(u,z_{2})v=(z_{2}+z_{0})^{k}\mathcal{Y}(Y(E^{2\lambda},z_{0})u,z_{2})v
\]
for $u\in\charlam{\lambda}[\lambda],\,v\in\charlam{\mu}[\mu]$ and sufficiently large integer $k$.
By using \eqref{ggggg} we have 
\begin{align*}
&(z_{0}+z_{2})^{k}Y(E^{2\lambda},z_{0}+z_{2})\mathcal{Y}_{\lambda,\mu}(\phi(u),z_{2})\phi(v)\\
&\quad=(z_{2}+z_{0})^{k}(\epsilon(2\lambda,\lambda+\mu)\mathcal{Y}_{3\lambda,\mu}(\mathcal{Y}_{2\lambda,\lambda}(e^{2\lambda},z_{0})\phi(u),z_{2})\phi(v)\\
&\qquad+\epsilon(-2\lambda,\lambda+\mu)\mathcal{Y}_{-\lambda,\mu}(\mathcal{Y}_{-2\lambda,\lambda}(e^{-2\lambda},z_{0})\phi(u),z_{2})\phi(v)).
\end{align*}
Similarly, 
\begin{align*}
&(z_{0}+z_{2})^{k}Y(E^{2\lambda},z_{0}+z_{2})\mathcal{Y}_{\lambda,-\mu}(\phi(u),z_{2})\theta(\phi(v))\\
&\quad=(z_{2}+z_{0})^{k}(\epsilon(2\lambda,\lambda-\mu)\mathcal{Y}_{3\lambda,-\mu}(\mathcal{Y}_{2\lambda,\lambda}(e^{2\lambda},z_{0})\phi(u),z_{2})\theta\phi(v))\\
&\qquad+\epsilon(-2\lambda,\lambda-\mu)\mathcal{Y}_{-\lambda,-\mu}(\mathcal{Y}_{-2\lambda,\lambda}(e^{-2\lambda},z_{0})\phi(u),z_{2})\theta(\phi(v))),\\
&(z_{0}+z_{2})^{k}Y(E^{2\lambda},z_{0}+z_{2})\mathcal{Y}_{-\lambda,\mu}(\theta(\phi(u)),z_{2})\phi(v)\\
&\quad=(z_{2}+z_{0})^{k}(\epsilon(2\lambda,-\lambda+\mu)\mathcal{Y}_{\lambda,\mu}(\mathcal{Y}_{2\lambda,-\lambda}(e^{2\lambda},z_{0})\theta(\phi(u)),z_{2})\phi(v)\\
&\qquad+\epsilon(-2\lambda,-\lambda+\mu)\mathcal{Y}_{-3\lambda,\mu}(\mathcal{Y}_{-2\lambda,-\lambda}(e^{-2\lambda},z_{0})\theta(\phi(u)),z_{2})\phi(v)),\\
&(z_{0}+z_{2})^{k}Y(E^{2\lambda},z_{0}+z_{2})\mathcal{Y}_{-\lambda,-\mu}(\theta(\phi(u)),z_{2})\theta(\phi(v))\\
&\quad=(z_{2}+z_{0})^{k}(\epsilon(2\lambda,-\lambda-\mu)\mathcal{Y}_{\lambda,-\mu}(\mathcal{Y}_{2\lambda,-\lambda}(e^{2\lambda},z_{0})\theta(\phi(u)),z_{2})\theta(\phi(v))\\
&\qquad+\epsilon(-2\lambda,-\lambda-\mu)\mathcal{Y}_{-3\lambda,-\mu}(\mathcal{Y}_{-2\lambda,-\lambda}(e^{2\lambda},z_{0})\theta(\phi(u)),z_{2})\theta(\phi(v))).
\end{align*}
Hence,
\begin{align*}
&(z_{0}+z_{2})^{k}Y(E^{2\lambda},z_{0}+z_{2})\mathcal{Y}(u,z_{2})v\\
&\quad=(z_{2}+z_{0})^{k}(c_{\lambda,\mu}^{(1)}\epsilon(2\lambda,\lambda+\mu)\mathcal{Y}_{3\lambda,\mu}(\mathcal{Y}_{2\lambda,\lambda}(e^{2\lambda},z_{0})\phi(u),z_{2})\phi(v)\\
&\qquad+c_{\lambda,\mu}^{(1)}\epsilon(-2\lambda,\lambda+\mu)\mathcal{Y}_{-\lambda,\mu}(\mathcal{Y}_{-2\lambda,\lambda}(e^{-2\lambda},z_{0})\phi(u),z_{2})\phi(v))\\
&\qquad+c_{\lambda,\mu}^{(2)}\epsilon(2\lambda,\lambda-\mu)\mathcal{Y}_{3\lambda,-\mu}(\mathcal{Y}_{2\lambda,\lambda}(e^{2\lambda},z_{0})\phi(u),z_{2})\theta\phi(v))\\
&\qquad+c_{\lambda,\mu}^{(2)}\epsilon(-2\lambda,\lambda-\mu)\mathcal{Y}_{-\lambda,-\mu}(\mathcal{Y}_{-2\lambda,\lambda}(e^{-2\lambda},z_{0})\phi(u),z_{2})\theta(\phi(v))),\\
&\qquad+c_{\lambda,\mu}^{(3)}\epsilon(2\lambda,-\lambda+\mu)\mathcal{Y}_{\lambda,\mu}(\mathcal{Y}_{2\lambda,-\lambda}(e^{2\lambda},z_{0})\theta(\phi(u)),z_{2})\phi(v)\\
&\qquad+c_{\lambda,\mu}^{(3)}\epsilon(-2\lambda,-\lambda+\mu)\mathcal{Y}_{-3\lambda,\mu}(\mathcal{Y}_{-2\lambda,-\lambda}(e^{-2\lambda},z_{0})\theta(\phi(u)),z_{2})\phi(v)),\\
&\qquad+c_{\lambda,\mu}^{(4)}\epsilon(2\lambda,-\lambda-\mu)\mathcal{Y}_{\lambda,-\mu}(\mathcal{Y}_{2\lambda,-\lambda}(e^{2\lambda},z_{0})\theta(\phi(u)),z_{2})\theta(\phi(v))\\
&\qquad+c_{\lambda,\mu}^{(4)}\epsilon(-2\lambda,-\lambda-\mu)\mathcal{Y}_{-3\lambda,-\mu}(\mathcal{Y}_{-2\lambda,-\lambda}(e^{2\lambda},z_{0})\theta(\phi(u)),z_{2})\theta(\phi(v))).
\end{align*}
Now we set
\begin{align*}
C^{i,j}&=\mathcal{Y}_{(2+(-1)^{i})\lambda,(-1)^{j}\mu}(\mathcal{Y}_{2\lambda,(-1)^{i}\lambda}(e^{2\lambda},z_{0})\theta^{i}(\phi(u)),z_{2})\theta^{j}(\phi(v)),\\
D^{i,j}&=\mathcal{Y}_{(-2+(-1)^{i})\lambda,(-1)^{j}\mu}(\mathcal{Y}_{-2\lambda,(-1)^{i}\lambda}(e^{-2\lambda},z_{0})\theta^{i}(\phi(u)),z_{2})\theta^{j}(\phi(v)).
\end{align*}
for $i=0,1$.
Since $C^{i,j}\in\Fremo{(2+(-1)^{i})\lambda,(-1)^{j}\mu}((z_{0}))((z_{2}))\subset \charlam{\lambda+\mu}((z_{0}))((z_{2}))$ and $D^{i,j}\in\Fremo{(-2+(-1)^{i})\lambda,(-1)^{j}\mu}((z_{0}))((z_{2}))\subset \charlam{\lambda+\mu}((z_{0}))((z_{2}))$, $C^{i,j}$ and $D^{i,j}$ for $i=0,1$ are linearly independent in $\charlam{\lambda+\mu}((z_{0}))((z_{2}))$. 
Using $C^{i,j}$ and $D^{i,j}$, we can rewrite the identity above as 
\begin{align*}
&(z_{0}+z_{2})^{k}Y(E^{2\lambda},z_{0}+z_{2})\mathcal{Y}(u,z_{2})v\\
&\quad=(z_{2}+z_{0})^{k}(c_{\lambda,\mu}^{(1)}\epsilon(2\lambda,\lambda+\mu)C^{0,0}+c_{\lambda,\mu}^{(1)}\epsilon(-2\lambda,\lambda+\mu)D^{0,0}\\
&\qquad+c_{\lambda,\mu}^{(2)}\epsilon(2\lambda,\lambda-\mu)C^{0,1}
+c_{\lambda,\mu}^{(2)}\epsilon(-2\lambda,\lambda-\mu)D^{0,1}\\
&\qquad+c_{\lambda,\mu}^{(3)}\epsilon(2\lambda,-\lambda+\mu)C^{1,0}+c_{\lambda,\mu}^{(3)}\epsilon(-2\lambda,-\lambda+\mu)D^{1,0}\\
&\qquad+c_{\lambda,\mu}^{(4)}\epsilon(2\lambda,-\lambda-\mu)C^{1,1}+c_{\lambda,\mu}^{(4)}\epsilon(-2\lambda,-\lambda-\mu)D^{1,1}.
\end{align*}

As before we note that
\[
\phi(Y(E^{2\lambda},z_{0})u)=\epsilon(2\lambda,\lambda)\mathcal{Y}_{2\lambda,\lambda}(e^{2\lambda},z_{0})\phi(u)+\epsilon(2\lambda,-\lambda)\mathcal{Y}_{2\lambda,-\lambda}(e^{2\lambda},z_{0})\theta(\phi(u)).
\] 
So 
\begin{align*}
(z_{0}+z_{2})^{k}&\mathcal{Y}(Y(E^{2\mu},z_{0})u,z_{2})v\\
=&(z_{0}+z_{2})^{k}(
c_{3\lambda,\mu}^{(1)}\epsilon(2\lambda,\lambda)\mathcal{Y}_{3\lambda,\mu}(\mathcal{Y}_{2\lambda,\lambda}(e^{2\lambda},z_{0})\phi(u),z_{2})\phi(v)\\
&+c_{3\lambda,\mu}^{(2)}\epsilon(2\lambda,\lambda)\mathcal{Y}_{3\lambda,-\mu}(\mathcal{Y}_{2\lambda,\lambda}(e^{2\lambda},z_{0})\phi(u),z_{2})\theta(\phi(v))\\
&+c_{3\lambda,\mu}^{(3)}\epsilon(2\lambda,\lambda)\mathcal{Y}_{-3\lambda,\mu}(\mathcal{Y}_{-2\lambda,-\lambda}(e^{-2\lambda},z_{0})\theta(\phi(u)),z_{2})\phi(v)\\
&+c_{3\lambda,\mu}^{(4)}\epsilon(2\lambda,\lambda)\mathcal{Y}_{-3\lambda,-\mu}(\mathcal{Y}_{-2\lambda,-\lambda}(e^{-2\lambda},z_{0})\theta(\phi(u)),z_{2})\theta(\phi(v))\\
&+c_{\lambda,\mu}^{(1)}\epsilon(2\lambda,-\lambda)\mathcal{Y}_{\lambda,\mu}(\mathcal{Y}_{2\lambda,-\lambda}(e^{2\lambda},z_{0})\theta(\phi(u)),z_{2})\phi(v)\\
&+c_{\lambda,\mu}^{(2)}\epsilon(2\lambda,-\lambda)\mathcal{Y}_{\lambda,-\mu}(\mathcal{Y}_{2\lambda,-\lambda}(e^{2\lambda},z_{0})\theta(\phi(u)),z_{2})\theta(\phi(v))\\
&+c_{\lambda,\mu}^{(3)}\epsilon(2\lambda,-\lambda)\mathcal{Y}_{-\lambda,\mu}(\mathcal{Y}_{-2\lambda,\lambda}(e^{-2\lambda},z_{0})\phi(u),z_{2})\phi(v)\\
&+c_{\lambda,\mu}^{(4)}\epsilon(2\lambda,-\lambda)\mathcal{Y}_{-\lambda,-\mu}(\mathcal{Y}_{-2\lambda,\lambda}(e^{-2\lambda},z_{0})\phi(u),z_{2})\theta(\phi(v))\\
=&(z_{1}-z_{2})^{k}(
  c_{3\lambda,\mu}^{(1)}\epsilon(2\lambda,\lambda)C^{0,0}+c_{3\lambda,\mu}^{(2)}\epsilon(2\lambda,\lambda)C^{0,1}\\
&+c_{3\lambda,\mu}^{(3)}\epsilon(2\lambda,\lambda)D^{1,0}+c_{3\lambda,\mu}^{(4)}\epsilon(2\lambda,\lambda)D^{1,1}\\
&+c_{\lambda,\mu}^{(1)}\epsilon(2\lambda,-\lambda)C^{1,0}+c_{\lambda,\mu}^{(2)}\epsilon(2\lambda,-\lambda)C^{1,1}\\
&+c_{\lambda,\mu}^{(3)}\epsilon(2\lambda,-\lambda)D^{0,0}+c_{\lambda,\mu}^{(4)}\epsilon(2\lambda,-\lambda)D^{0,1}).
\end{align*}  
Since $C^{i,j}$ and $D^{i,j}$ for $i,j=0,1$ are linearly independent, the associativity formula implies the equations
\begin{align*}
&c_{\lambda,\mu}^{(1)}\epsilon(2\lambda,\lambda+\mu)=c_{3\lambda,\mu}^{(1)}\epsilon(2\lambda,\lambda),\quad c_{\lambda,\mu}^{(1)}\epsilon(-2\lambda,\lambda+\mu)=c_{\lambda,\mu}^{(3)}\epsilon(2\lambda,-\lambda),\\
&c_{\lambda,\mu}^{(2)}\epsilon(2\lambda,\lambda-\mu)=c_{3\lambda,\mu}^{(2)}\epsilon(2\lambda,\lambda),\quad c_{\lambda,\mu}^{(2)}\epsilon(-2\lambda,\lambda-\mu)=c_{\lambda,\mu}^{(4)}\epsilon(2\lambda,-\lambda),\\
&c_{\lambda,\mu}^{(3)}\epsilon(2\lambda,-\lambda+\mu)=c_{\lambda,\mu}^{(1)}\epsilon(2\lambda,-\lambda),\quad c_{\lambda,\mu}^{(3)}\epsilon(-2\lambda,-\lambda+\mu)=c_{3\lambda,\mu}^{(3)}\epsilon(2\lambda,\lambda),\\
&c_{\lambda,\mu}^{(4)}\epsilon(2\lambda,\lambda-\mu)=c_{\lambda,\mu}^{(2)}\epsilon(2\lambda,-\lambda),\quad c_{\lambda,\mu}^{(4)}\epsilon(-2\lambda,\lambda+\mu)=c_{3\lambda,\mu}^{(4)}\epsilon(2\lambda,\lambda).
\end{align*}
This proves that 
\begin{align}\label{llpp}
c_{\lambda,\mu}^{(3)}=\epsilon(-2\lambda,\mu)c_{\lambda,\mu}^{(1)}=\epsilon(-\lambda,2\mu)c_{\lambda,\mu}^{(1)}.
\end{align} 
Combining \eqref{llpp} with \eqref{lllp} we see that 
\begin{align*}
\mathcal{Y}(u,z)v=&c_{\lambda,\mu}^{(1)}(\mathcal{Y}_{\lambda,\mu}(\phi(u),z_{2})\phi(v)+(-1)^{(2\mu,\lambda)}\epsilon(-2\mu,\lambda)\mathcal{Y}_{\lambda,-\mu}(\phi(u),z_{2})\theta(\phi(v))\\
&+\epsilon(-\lambda,2\mu)\mathcal{Y}_{-\lambda,\mu}(\theta(\phi(u)),z_{2})\phi(v)+(-1)^{(2\mu,\lambda)}\mathcal{Y}_{-\lambda,\mu}(\theta(\phi(u)),z_{2})\theta(\phi(v)).
\end{align*}
Thus the image of $I_{\charge{+}}\fusion{\charlam{\lambda}^{+}}{\charlam{\mu}^{+}}{\charlam{\lambda+\mu}}$ in $I_{\Free{+}}\fusion{\charlam{\lambda}^{+}[\lambda]}{\charlam{\mu}^{+}[\mu]}{\charlam{\lambda+\mu}}$ is spanned by one intertwining operator, in particular, the dimension is one.
This concludes that the fusion rule of type $\fusion{\charlam{\lambda}^{+}}{\charlam{\mu}^{+}}{\charlam{\lambda+\mu}}$ is at most one.  
This completes the proof. 
\end{proof}

In view of Propositions \ref{duality}, \ref{untw}, \ref{qqqq}, \ref{qqqqq} 
and \ref{bbbb} we immediately have:

\begin{proposition}\label{fffffff} 
Let $M^i\,(i=1,2,3)$ be irreducible $\charge{+}$-modules of untwisted type.
Then the fusion rule of type $\fusion{M^1}{M^2}{M^3}$ is either $0$ or $1$. 
The fusion rule of type $\fusion{M^1}{M^2}{M^3}$ is $1$ 
if and only if $M^i\,(i=1,2,3)$ satisfy the following conditions;

\noindent
{\rm(i)} $M^1=\charlam{\lambda}$ for $\lambda\in L^{\circ}$ 
such that $2\lambda\notin L$ and $(M^2,\,M^3)$ 
is one of the following pairs:
\begin{itemize}
\item[] $(\charlam{\mu},\charlam{\nu})$ 
for $\mu,\nu\in L^{\circ}$ such that $2\mu,\,2\nu\notin L$ 
and $(\lambda,\,\mu,\,\nu)$ is an admissible triple modulo $L$,

\item[] $(\charlam{\mu}^{\pm},\charlam{\nu}),\,((\charlam{\nu})',
\,(\charlam{\mu}^{\pm})')$ for $\mu,\nu\in L^{\circ}$ 
such that $2\mu\in L$ and $(\lambda,\,\mu,\,\nu)$ 
is an admissible triple modulo $L$. 
\end{itemize}

\noindent
{\rm(ii)} $M^1=\charlam{\lambda}^+$ for $\lambda\in L^{\circ}$ 
such that $2\lambda\in L$ and $(M^2,\,M^3)$ is one of the following pairs:
\begin{itemize}
\item[] $(\charlam{\mu},\charlam{\nu})$ for $\mu,\nu\in L^{\circ}$ 
such that $2\mu\notin L$ and $(\lambda,\,\mu,\,\nu)$ 
is an admissible triple modulo $L$,
\item[] $(\charlam{\mu}^{\pm},\charlam{\nu}^{\pm})$ 
for $\mu,\nu\in L^{\circ}$ such that $2\mu\in
L,\,\pi_{\lambda,2\mu}=1$ 
and $(\lambda,\,\mu,\,\nu)$ is an admissible triple modulo $L$,
\item[] $(\charlam{\mu}^{\pm},\charlam{\nu}^{\mp})$ 
for $\mu,\nu\in L^{\circ}$ such that 
$2\mu\in L,\,\pi_{\lambda,2\mu}=-1$ and $(\lambda,\,\mu,\,\nu)$ 
is an admissible triple modulo $L$.
\end{itemize} 

\noindent
{\rm(iii)} $M^1=\charlam{\lambda}^-$ for $\lambda\in L^{\circ}$ 
such that $2\lambda\in L$  and $(M^2,\,M^3)$ is one of the following pairs:
\begin{itemize}
\item[] $(\charlam{\mu},\charlam{\nu})$ for $\mu,\nu\in L^{\circ}$ 
such that $2\mu\notin L$ and $(\lambda,\,\mu,\,\nu)$ 
is an admissible triple modulo $L$,
\item[] $(\charlam{\mu}^{\pm},\charlam{\nu}^{\mp})$ 
for $\mu,\nu\in L^{\circ}$ 
such that $2\mu\in L,\,\pi_{\lambda,2\mu}=1$ 
and $(\lambda,\,\mu,\,\nu)$ is an admissible triple modulo $L$,
\item[] $(\charlam{\mu}^{\pm},\charlam{\nu}^{\pm})$ 
for $\mu,\nu\in L^{\circ}$ 
such that $2\mu\in L,\,\pi_{\lambda,2\mu}=-1$ 
and $(\lambda,\,\mu,\,\nu)$ is an admissible triple modulo $L$.
\end{itemize} 
\end{proposition}

\subsection{Fusion rules involving modules of twisted type}\label{FRIM}
We construct nonzero intertwining operators among irreducible $\charge{+}$-modules involving modules of twisted type in this section. We use 
$\chi$ for a central character of $\hat{L}/K$ with $\chi(\kappa)=-1$ and use
$T_{\chi}$ to denote the corresponding irreducible $\hat{L}/K$-module.

Let $\lambda\in L^{\circ}$. 
We define an automorphism $\sigma_{\lambda}$ of $\hat{L}$ by 
\begin{align*}
\sigma_{\lambda}(a)=\kappa^{(\lambda,\bar{a})}a
\end{align*}
for any $a\in \hat{L}$. 
Since $\sigma_{\lambda}(\theta(a))=\theta(\sigma_{\lambda}(a))$, $\sigma_{\lambda}$ stabilizes $K$.
Hence $\sigma_{\lambda}$ induces an automorphism of $\hat{L}/K$.

For any $\hat{L}/K$-module $T$ we denote by $T\circ\sigma_{\lambda}$ the $\hat{L}/K$-module twisted by $\sigma_{\lambda}.$
That is, $T\circ\sigma_{\lambda}=T$ as vector spaces but with
a new action defined by $a.t=\sigma_{\lambda}(a)t$ for $a\in\hat{L}/K$ and $t\in T$. 
Since $T_{\chi}\circ\sigma_{\lambda}$ is also irreducible, there is a unique central character $\chi^{(\lambda)}$ of $\hat{L}/K$ (with $\chi^{(\lambda)}(\kappa)=-1$), such that $T_{\chi}\circ\sigma_{\lambda}\cong T_{\chi^{(\lambda)}}$.
Let $f$ be an $\hat{L}/K$-module isomorphism $T_{\chi}\circ\sigma_{\lambda}\mathop{\rightarrow}\limits^{\sim} T_{\chi^{(\lambda)}}$.
Then $f$ is a linear isomorphism from $T_{\chi}$ to $T_{\chi^{(\lambda)}}$ satisfying 
\begin{align}\label{conjuga1}
f(\sigma_{\lambda}(a)t)=af(t)
\end{align}
for $a\in\hat{L}/K$ and $t\in T_{\chi}$. 

We now fix $\lambda\in L^{\circ}$ and an $\hat{L}/K$-module isomorphism 
$f:T_{\chi}\circ\sigma_{\lambda}\to T_{\chi^{(\lambda)}}$.
For any $\alpha\in L$, we define a linear isomorphism 
$\eta_{\lambda+\alpha}:T_{\chi}\circ\sigma_{\lambda}\to T_{\chi^{(\lambda)}}$ 
by $\eta_{\lambda+\alpha}=\epsilon(-\alpha,\lambda)e_{\alpha}\circ f$, 
where we write $\epsilon(\mu,\nu)=\w_{q}^{\epsilon_{0}(\mu,\mu)}$ 
for $\mu,\nu\in L^{\circ}$ as before. 
Then we have a linear isomorphism 
\begin{align*}
\eta_{\gamma}:T_{\chi}\rightarrow T_{\chi^{(\lambda)}}
\end{align*}
for any $\gamma\in\lambda+L$.

\begin{lemma}\label{iiiii}
For any $\gamma\in\lambda+L$ and $\alpha\in L$, 
\begin{align}
&e_{\alpha}\circ\eta_{\gamma}
=(-1)^{(\alpha,\gamma)}\eta_{\gamma}\circ e_{\alpha}\label{intertwiner1},\\
&e_{\alpha}\circ\eta_{\gamma}=\epsilon(\alpha,\gamma)\eta_{\gamma+\alpha}
=\epsilon(-\alpha,\gamma)\eta_{\gamma-\alpha}.\label{intertwiner2}
\end{align}
\end{lemma}

\begin{proof}
Set $\beta=\gamma-\lambda\in L$. 
Since $e_{\alpha}\circ f=(-1)^{(\alpha,\lambda)}f\circ e_{\alpha}$ and 
$e_{\alpha}e_{\beta}=(-1)^{(\alpha,\beta)}e_{\beta}e_{\alpha}$, 
we have $e_{\alpha}\circ\eta_{\gamma}=(-1)^{(\alpha,\gamma)}\eta_{\gamma}\circ e_{\alpha}$.
This proves \eqref{intertwiner1}.
By definition we have
\begin{align*}
e_{\alpha}\circ\eta_{\gamma}&=\epsilon(-\beta,\lambda)e_{\alpha}\circ e_{\beta}\circ f\\
&=\epsilon(-\beta,\lambda)\epsilon(\alpha,\beta)e_{\alpha+\beta}\circ f\\
&=\epsilon(-\beta,\lambda)\epsilon(\alpha,\beta)\epsilon(\alpha+\beta,\lambda)\eta_{\gamma+\alpha}\\
&=\epsilon(\alpha,\gamma)\eta_{\gamma+\alpha}.
\end{align*}
Thus the first equality in \eqref{intertwiner2} holds.
The second equality in \eqref{intertwiner2} follows from the fact that 
$e_{-\alpha}=\theta(e_{\alpha})=e_{\alpha}$ on $T_{\chi}$.
\end{proof}
\begin{remark}
In the case $L=\Z\alpha$ of rank one, there are two irreducible 
$\hat{L}/K$-modules $T^{1},\,T^{2}$ on which $e_{\alpha}$ acts as $1$ and $-1$ respectively.
Then for any $\lambda=\frac{r}{|\alpha|^{2}}\alpha\in L^{\circ}$, $\eta_{\lambda}$ 
stabilizes $T^{i}$ for $i=1,2$ if $r$ is even and switches $T^{1}$ and $T^{2}$ if $r$ is odd. 
Thus the map $\eta_{\lambda}$ coincides with $\psi_{\lambda}$ in \cite{A2} up to a scalar multiple.
\end{remark}

Let $\lambda\in L^{\circ}$. 
Recall operators $\mathcal{Y}_{\lambda,\mu}(\,\cdot\,,z)$ 
and $\mathcal{Y}^{\rm tw}_{\lambda}(\,\cdot\,,z)$  
defined in \eqref{untwistop} and \eqref{twining1}.
Following the arguments in \cite[Chapter 9]{FLM}, 
we have the  following identity 
for any $\alpha\in L,\,\lambda\in L^{\circ},\,a\in\Fremo{\alpha}$ and 
$u\in\Fremo{\lambda}$
\begin{align}\label{jacobi2}
\begin{split}
&z_{0}^{-1}\delta\left({\frac{z_{1}-z_{2}}{z_{0}}}\right)\mathcal{Y}^{\rm tw}_{\alpha}(a,z_{1})\mathcal{Y}^{\rm tw}_{\lambda}(u,z_{2})\\
&\qquad-(-1)^{(\alpha,\lambda)}z_{0}^{-1}\delta\left({\frac{z_{2}-z_{1}}{-z_{0}}}\right)\mathcal{Y}^{\rm tw}_{\lambda}(u,z_{2})\mathcal{Y}^{\rm tw}_{\alpha}(a,z_{1})\\
&\quad=\frac{1}{2}\sum_{p=0,1}z_{2}^{-1}\delta\left((-1)^{p}{\frac{(z_{1}-z_{0})^{1/2}}{z_{2}^{1/2}}}\right)\mathcal{Y}^{\rm tw}_{\lambda+(-1)^{p}\alpha}(\mathcal{Y}_{(-1)^{p}\alpha,\lambda}(\theta^{p}(a),z_{0})u,z_{2})
\end{split}
\end{align}
and 
\begin{align}\label{deriver2}
\mathcal{Y}^{\rm tw}_{\lambda}(L(-1)u,z)=\frac{d}{dz}\mathcal{Y}^{\rm tw}_{\lambda}(u,z)
\end{align} 
on $\Fretw{}$.

Now we define an operator $\tilde{\mathcal{Y}}^{\rm tw}(u,\cdot\,,z)$ 
from $\charge{T_{\chi}}$ to $\charge{T_{\chi}^{(\lambda)}}$ 
by 
\begin{align}\label{twistop}
\tilde{\mathcal{Y}}^{\rm tw}_{\lambda}(u,z)=\mathcal{Y}_{\lambda+\beta}^{\rm tw}(u,z)\otimes \eta_{\lambda+\beta}
\end{align}
for any $u\in\Fremo{\lambda+\beta}\subset\charlam{\lambda}$.
So we have a linear map 
\begin{align*}
\tilde{\mathcal{Y}}^{\rm tw}_{\lambda}(\,\cdot\,,z):\charlam{\lambda}\rightarrow\Hom(\charge{T_{\chi}},\charge{T^{(\lambda)}_{\chi}})\{z\}.
\end{align*}

We remark that if $\lambda=0$ then $T^{(\lambda)}_{\chi}=T_{\chi},$
$\eta_{\lambda+\alpha}=e_{\alpha}$ and  
$\tilde{\mathcal{Y}}^{\rm tw}_{\lambda}(a,z)=\mathcal{Y}^{\rm tw}_{\alpha}(a,z_{1})\otimes e_{\alpha}$ is exactly
the twisted vertex operator $Y(a,z)$ associated to $a\in\Fremo{\alpha}\subset\charge{}$ which defines the twisted module structure on 
$\charge{T_{\chi}}$ (see [FLM]). 
 
\begin{proposition}\label{twistedinter}
Let $\lambda\in L^{\circ}$ and $\chi$ be a central character of $\hat{L}/K$ with $\chi(\kappa)=-1$. 
Then for any $a\in\charge{}$ and $u\in\charlam{\lambda}$, the identities
\begin{align*}
\begin{split}
&z_{0}^{-1}\delta\left({\frac{z_{1}-z_{2}}{z_{0}}}\right)
Y(a,z_{1})\tilde{\mathcal{Y}}^{\rm tw}(u,z_{2})
-z_{0}^{-1}\delta\left({\frac{z_{2}-z_{1}}{-z_{0}}}\right)
\tilde{\mathcal{Y}}^{\rm tw}(u,z_{2})Y(a,z_{1})\\
&\quad=\frac{1}{2}\sum_{p=0,1}z_{2}^{-1}
\delta\left((-1)^{p}{\frac{(z_{1}-z_{0})^{1/2}}{z_{2}^{1/2}}}\right)
\tilde{\mathcal{Y}}^{\rm tw}(Y(\theta^{p}(a),z_{0})u,z_{2})
\end{split}
\end{align*}
and 
\[
\frac{d}{dz}\tilde{\mathcal{Y}}^{\rm tw}(u,z)=\tilde{\mathcal{Y}}^{\rm tw}(L(-1)u,z)\] 
hold on $\charge{T_{\chi}}$. 
In particular, $\tilde{\mathcal{Y}}^{\rm tw}(\,\cdot\,,z)$ is an intertwining operator 
of type $\fusion{\charlam{\lambda}}{\charge{T_{\chi}}}{\charge{T_{\chi}^{(\lambda)}}}$ 
for $\charge{+}$.
\end{proposition}
\begin{proof} 
By \eqref{intertwiner1}--\eqref{jacobi2} we see that for any $a\in\Fremo{\alpha}\subset\charge{}$ and $u\in\Fremo{\gamma}$ with $\gamma\in\lambda+L$, 
\begin{align*}
&z_{0}^{-1}\delta\left({\frac{z_{1}-z_{2}}{z_{0}}}\right)Y(a,z_{1})\tilde{\mathcal{Y}}^{\rm tw}_{\lambda}(u,z_{2})-z_{0}^{-1}\delta\left({\frac{z_{2}-z_{1}}{-z_{0}}}\right)\tilde{\mathcal{Y}}^{\rm tw}_{\lambda}(u,z_{2})Y(a,z_{1})\\
&\quad=z_{0}^{-1}\delta\left({\frac{z_{1}-z_{2}}{z_{0}}}\right)\left(\mathcal{Y}^{\rm tw}_{\alpha}(a,z_{1})\otimes e_{\alpha}\right)\left(\mathcal{Y}^{\rm tw}_{\gamma}(u,z_{2})\otimes \eta_{\gamma}\right)\\
&\qquad-z_{0}^{-1}\delta\left({\frac{z_{2}-z_{1}}{-z_{0}}}\right)\left(\mathcal{Y}^{\rm tw}_{\gamma}(u,z_{2})\otimes \eta_{\gamma}\right)\left(\mathcal{Y}^{\rm tw}_{\alpha}(a,z_{1})\otimes e_{\alpha}\right)\\
&\quad=z_{0}^{-1}\delta\left({\frac{z_{1}-z_{2}}{z_{0}}}\right)\mathcal{Y}^{\rm tw}_{\alpha}(a,z_{1})\mathcal{Y}^{\rm tw}_{\gamma}(u,z_{2})\otimes \left(e_{\alpha}\circ\eta_{\gamma}\right)\\
&\qquad-z_{0}^{-1}\delta\left({\frac{z_{2}-z_{1}}{-z_{0}}}\right)\mathcal{Y}^{\rm tw}_{\gamma}(u,z_{2})\mathcal{Y}^{\rm tw}_{\alpha}(a,z_{1})\otimes \left(\eta_{\gamma}\circ e_{\alpha}\right)\\
&\quad=\left(z_{0}^{-1}\delta\left({\frac{z_{1}-z_{2}}{z_{0}}}\right)\mathcal{Y}^{\rm tw}_{\alpha}(a,z_{1})\mathcal{Y}^{\rm tw}_{\lambda}(u,z_{2})\right.\\
&\qquad\left.-(-1)^{\langle\alpha,\gamma\rangle}z_{0}^{-1}\delta\left({\frac{z_{2}-z_{1}}{-z_{0}}}\right)\mathcal{Y}^{\rm tw}_{\gamma}(u,z_{2})\mathcal{Y}^{\rm tw}_{\alpha}(a,z_{1})\right)\otimes\left(e_{\alpha}\circ\eta_{\gamma}\right)\\
&\quad=\frac{1}{2}z_{2}^{-1}\delta\left({\frac{(z_{1}-z_{0})^{1/2}}{z_{2}^{1/2}}}\right)\mathcal{Y}^{\rm tw}_{\gamma+\alpha}(\mathcal{Y}_{\alpha,\gamma}(a,z_{0})u,z_{2})\otimes\left(\epsilon(\alpha,\gamma)\eta_{\gamma+\alpha}\right)\\
&\qquad+\frac{1}{2}z_{2}^{-1}\delta\left(-{\frac{(z_{1}-z_{0})^{1/2}}{z_{2}^{1/2}}}\right)\mathcal{Y}^{\rm tw}_{\gamma-\alpha}(\mathcal{Y}_{-\alpha,\gamma}(\theta(a),z_{0})u,z_{2})\otimes\left(\epsilon(-\alpha,\gamma)\eta_{\gamma-\alpha}\right)\\
&\quad=\frac{1}{2}z_{2}^{-1}\delta\left({\frac{(z_{1}-z_{0})^{1/2}}{z_{2}^{1/2}}}\right)\mathcal{Y}^{\rm tw}_{\gamma+\alpha}(Y(a,z_{0})u,z_{2})\otimes\eta_{\gamma+\alpha}\\
&\qquad+\frac{1}{2}z_{2}^{-1}\delta\left(-{\frac{(z_{1}-z_{0})^{1/2}}{z_{2}^{1/2}}}\right)\mathcal{Y}^{\rm tw}_{\gamma-\alpha}(Y(\theta(a),z_{0})u,z_{2})\otimes\eta_{\gamma-\alpha}\\
&\quad=\frac{1}{2}\sum_{p=0,1}z_{2}^{-1}\delta\left((-1)^{p}{\frac{(z_{1}-z_{0})^{1/2}}{z_{2}^{1/2}}}\right)\tilde{\mathcal{Y}}^{\rm tw}_{\lambda}(Y(\theta^{p}(a),z_{0})u,z_{2}).
\end{align*}
It follows from \eqref{deriver2} that $\tilde{\mathcal{Y}}^{\rm tw}(\,\cdot\,,z)$ satisfies the $L(-1)$-derivative property. 
The last assertion is clear.
\end{proof}
 
%We next discuss how to construct intertwining operators of type
%$\fusion{\charlam{\lambda}}{\charge{T_{\chi},\epsilon_{1}}}{\charge{T_{\chi}^{(\lambda)},\epsilon_{2}}}$ for $\charge{+}$.  
We recall the canonical projection $p_{\pm}:\Fretw{}\to\Fretw{\pm}$ and the canonical injection $\iota_{\pm}:\Fretw{\pm}\to\Fretw{}$.   
We then have the projection $p_{\pm}\otimes\id:\charge{T}\to\charge{T,\pm}$ and the injection $\iota_{\pm}\otimes\id:\charge{T,\pm}\to\charge{T}$ for any irreducible $\hat{L}/K$-module on which $\kappa=-1$, noting that $\charge{Y}=\Fretw{}\otimes T$.  
We also write for them by $p_{\pm}$ and $\iota_{\pm}$ respectively. 
Let $\epsilon_{1},\epsilon_{2}\in\{\pm\}$ and $\lambda\in L^\circ$.  
It is clear from the definition that the restriction
 $p_{\epsilon_2}^{\rm tw}\circ\tilde{\mathcal{Y}}^{\rm tw}(\,\cdot\,,z)\circ\iota_{\epsilon_1}$ is a nonzero intertwining operator of type
$\fusion{\charlam{\lambda}}{\charge{T_{\chi},\epsilon_{1}}}{\charge{T_{\chi}^{(\lambda)},\epsilon_{2}}}$ for $\charge{+}$.  Thus
we have:

\begin{proposition}\label{twist1}
For any $\lambda\in L^{\circ}$, the fusion rules of types 
$\fusion{\charlam{\lambda}}{\charge{T_{\chi},\pm}}{\charge{T_{\chi}^{(\lambda)},\pm}}$ 
and $\fusion{\charlam{\lambda}}{\charge{T_{\chi},\pm}}{\charge{T_{\chi}^{(\lambda)},\mp}}$ 
for $\charge{+}$ are nonzero. 
\end{proposition}

We now consider the case $2\lambda\in L$.
Let $\tilde{\mathcal{Y}}^{\rm tw}_{\lambda}(\,\cdot\,,z)$ be the intertwining operator of type 
$\fusion{\charlam{\lambda}}{\charge{T_{\chi}}}{\charge{T_{\chi}^{(\lambda)}}}$ 
defined in \eqref{twistop}.
By the conjugation formula \eqref{twistconjugate}, one has 
\begin{align}\label{conjuga2}
\theta\tilde{\mathcal{Y}}^{\rm tw}_{\lambda}(u,z)\theta^{-1}(v\otimes t)
=(\mathcal{Y}_{-\lambda-\alpha}^{\rm tw}(\theta(u),z)v)\otimes \eta_{\lambda+\alpha}(t)
\end{align}
for $\alpha\in L,\,u\in\Fremo{\lambda+\alpha},\,v\in\Fretw{}$ and $t\in T_{\chi}$.

By \eqref{intertwiner2}
\begin{align*}
\eta_{-\lambda-\alpha}(t)&=\epsilon(2\lambda+\alpha,\lambda)e_{-2\lambda-\alpha}\eta_{\lambda}(t)\\
&=\epsilon(2\lambda+\alpha,\lambda)
\epsilon(\alpha,2\lambda)(-1)^{(\alpha,2\lambda)}e_{-2\lambda}e_{-\alpha}\eta_{\lambda}(t)\\
&=\epsilon(2\lambda,\lambda)\epsilon(\alpha,3\lambda)e_{2\lambda}e_{\alpha}\eta_{\lambda}(t)\\
&=\epsilon(2\lambda,\lambda)\epsilon(\alpha,4\lambda)e_{2\lambda}\eta_{\lambda+\alpha}(t)\\
&=\epsilon(2\lambda,\lambda)e_{2\lambda}\eta_{\lambda+\alpha}(t).
\end{align*}  
Note that $e_{2\lambda}$ is in the center of $\hat{L}$ 
as $(2\lambda,\beta)\in2\Z$ for any $\beta\in L$.
Therefore, $e_{2\lambda}$ acts on $T^{(\lambda)}$ by the scalar 
$\chi^{(\lambda)}(e_{2\lambda})=(-1)^{(\lambda,2\lambda)}\chi(e_{2\lambda})=\chi(e_{2\lambda}).$ 
Hence   
\begin{align}\label{intertwiner3}
\eta_{-\lambda-\alpha}(t)
=\epsilon(2\lambda,\lambda)(-1)^{(\lambda,2\lambda)}\chi(e_{2\lambda})\eta_{\lambda+\alpha}(t)
=c_{\chi}(\lambda)\eta_{\lambda+\alpha}(t),
\end{align} 
where $c_{\chi}(\lambda)$ is the constant defined in \eqref{dddddd}. 
It follows from  \eqref{conjuga2} and \eqref{intertwiner3} that   
\begin{align}\label{conjuga3}
\theta\tilde{\mathcal{Y}}^{\rm tw}_{\lambda}(u,z)\theta^{-1}w
=c_{\chi}(\lambda)^{-1}\tilde{\mathcal{Y}}^{\rm tw}_{\lambda}(\theta(u),z)w
\end{align}
for any $u\in\charge{\lambda_{i}},\,w\in\charge{T,\epsilon_{2}}$.
It is clear that $c_{\chi}(\lambda)$ depends on $\lambda$ up to modulo $L$.

Recall that $c_{\chi}(\lambda)\in\{\pm1\}$.
We have the following proposition:

\begin{proposition}\label{twistporp3}
Let $\chi$ be a central character of $\hat{L}/K$ such that $\chi(\kappa)=-1$. 
For any $\lambda\in L^{\circ}$ with $2\lambda\in L$, the fusion rules of types 
$\fusion{\charlam{\lambda}^{+}}{\charge{T_{\chi},\pm}}{\charge{T_{\chi}^{(\lambda)},\pm}}$ 
and $\fusion{\charlam{\lambda}^{-}}{\charge{T_{\chi},\pm}}{\charge{T_{\chi}^{(\lambda)},\mp}}$ 
are nonzero if $c_{\chi}(\lambda)=1$, and the fusion rules of types 
$\fusion{\charlam{\lambda}^{+}}{\charge{T_{\chi},\pm}}{\charge{T_{\chi}^{(\lambda)},\mp}}$ 
and $\fusion{\charlam{\lambda}^{+}}{\charge{T_{\chi},\pm}}{\charge{T_{\chi}^{(\lambda)},\pm}}$ 
are nonzero if $c_{\chi}(\lambda)=-1$. 
\end{proposition}

We are now in the position to prove that the fusion rules of type 
$\fusion{M^{1}}{M^{2}}{M^{3}}$ is less than $1$ if one of $M^{1},\,M^{2}$ 
and $M^{3}$ is of twisted type.

\begin{proposition}\label{proptwst}
Let $M^{1},\,M^{2}$ and $M^{3}$ be irreducible $\charge{+}$-modules.
The fusion rule of type $\fusion{M^{1}}{M^{2}}{M^{3}}$ is $0$ if one
of $M^{1},\,M^{2}$ and $M^{3}$ is of twisted type and if the others
are of untwisted type or if all of $M^{1},\,M^{2}$ and $M^{3}$ are of
twisted type.
\end{proposition}

\begin{proof}
First we consider the case that one of $M^{1},\,M^{2}$ and $M^{3}$ is
of twisted type and the others are of untwisted type.  In view of
Propositions \ref{duality} and \ref{eeee}, we may assume that $M^{1}$
and $M^{2}$ are modules of untwisted types.  Then there exist
$\lambda,\,\mu\in L^{\circ}$ such that $M^{1}$ and $M^{2}$ contains
irreducible $\Free{+}$-submodules isomorphic to $\Fremo{\lambda}$ and
$\Fremo{\mu}$, respectively.  By Proposition
\ref{cfusion-rule-inequality}, the fusion rule of type
$\fusion{M^{1}}{M^{2}}{M^{3}}$ is less than or equal to the fusion
rule of type $\fusion{\Fremo{\lambda}}{\Fremo{\mu}}{M^{3}}$ for
$\Free{+}$.  Since $M^{3}$ is a module of twisted type and is a direct
sum of irreducible $\Free{+}$ modules isomorphic $\Fretw{+}$ or
$\Fretw{-}$, the fusion rule of type $\fusion{M^{1}}{M^{2}}{M^{3}}$ is
$0$ by Theorem \ref{ehigher-rank-m1}.

Next we consider the case all $M^{1},\,M^{2}$ and $M^{3}$ are of twisted type.
Then each $M^{i}$ is a direct sum of $\Fretw{+}$ or $\Fretw{-}$.
Proposition \ref{cfusion-rule-inequality} and 
Theorem \ref{ehigher-rank-m1} show that the fusion rule of type 
$\fusion{M^{1}}{M^{2}}{M^{3}}$ is $0$.
\end{proof}

Let $\chi_{1}$ and $\chi_{2}$ be central characters of $\hat{L}/K$
such that $\chi_{i}(\kappa)=-1$ and $M$ an irreducible
$\charge{+}$-module of untwisted type.  We shall prove that for any
$\epsilon_{1},\,\epsilon_{2}\in\{\pm\}$, the fusion rule of type
$\fusion{M}{\charge{T_{\chi_{1}},\epsilon_{1}}}{\charge{T_{\chi_{2}},
\epsilon_{2}}}$ is $0$ if $\chi_{2}\neq\chi_{1}^{(\lambda)}$.

Suppose that the fusion rule of type
$\fusion{M}{\charge{T_{\chi_{1}},\epsilon_{1}}}{\charge{T_{\chi_{2}},
\epsilon_{2}}}$ is nonzero, and let $\mathcal{Y}(\,\cdot\,,z)$ be a
nonzero intertwining operator of the corresponding type.  Since $M$ is
an irreducible $\charge{+}$-module of untwisted type, there is an
$\Free{+}$-submodule $W$ isomorphic to $\Fremo{\lambda}$ for some
$\lambda\in L^{\circ}$.  Let $\xi$ be an $\Free{+}$-module isomorphism
from $W$ to $\Fremo{\lambda}$, and define
$\tilde{\mathcal{Y}}(\,\cdot\,,z)$ by
\[
\tilde{\mathcal{Y}}(u,z)v=\mathcal{Y}(\xi^{-1}(u),z)v
\]  
for $u\in\Fremo{\lambda}$ and $v\in\charge{T_{\chi_{1}},\epsilon_{1}}$.
Then $\tilde{\mathcal{Y}}(\,\cdot\,,z)$ is a nonzero intertwining operator of type
$\fusion{\Fremo{\lambda}}{\charge{T_{\chi_{1}},\epsilon_{1}}}{\charge{T_{\chi_{2}},\epsilon_{2}}}$ for $\Free{+}$.  
Since $\charge{T_{\chi_{i}},\epsilon_{i}}\cong\Fretw{\epsilon_{i}}\otimes T_{\chi_{i}}$ as $\Free{+}$-modules, we have the following isomorphism of vector spaces
\begin{align}\label{grouphom}
I_{\Free{+}}\fusion{\Fremo{\lambda}}{\charge{T_{\chi_{1}},\epsilon_{1}}}{\charge{T_{\chi_{2}},\epsilon_{2}}}
\cong I_{\Free{+}}\fusion{\Fremo{\lambda}}{\Fretw{\epsilon_{1}}}{\Fretw{\epsilon_{2}}}\otimes \Hom_{\C}(T_{\chi_{1}},T_{\chi_{2}}).
\end{align}
We recall that $I_{\Free{+}}\fusion{\Fremo{\lambda}}{\Fretw{\epsilon_{1}}}{\Fretw{\epsilon_{2}}}$ is one dimensional and has a basis $p_{\epsilon_{2}}\circ{\mathcal{Y}}^{\rm tw}_{\lambda}(\,\cdot\,,z)\circ \iota_{\epsilon_{1}}$.
By using \eqref{grouphom} we see that there exists $f_{\lambda}\in\Hom(T_{\chi_{1}},T_{\chi_{2}})$ such that 
\begin{align}\label{wwww}
\tilde{\mathcal{Y}}(u,z)(v\otimes t)=\left(p_{\epsilon_{2}}({\mathcal{Y}}_{\lambda}^{\rm tw}(u,z)\iota_{\epsilon_{1}}(v))\right)\otimes f_{\lambda}(t)
\end{align}
for any $u\in\Fremo{\lambda},\,v\in\Fretw{\epsilon_{1}}$ and $t\in T_{\chi_{1}}$.
The vertex operator $Y(a,z)$ associated to $a\in\charge{+}[\alpha]$ acts on $\charge{T_{\chi_{2}},\epsilon_{2}}$ as $Y(a,z)=\left(\mathcal{Y}_{\alpha}^{\rm tw}(b,z)+\mathcal{Y}_{-\alpha}^{\rm tw}(\theta(b),z)\right)\otimes e_{\alpha}$, where $a=b+\theta(b)$ with $b\in\Fremo{\alpha}$. 
Thus we have 
\begin{align}\label{twjacobi1}
\begin{split}
&Y(a,z_{1})\mathcal{Y}(u,z_{2})(v\otimes t)\\
&\quad=\left(p_{\epsilon_{2}}\left(\mathcal{Y}_{\alpha}^{\rm tw}(b,z_{1})+\mathcal{Y}_{-\alpha}^{\rm tw}(\theta(b),z_{1})\right){\mathcal{Y}}_{\lambda}^{\rm tw}(\xi(u),z_{2})\iota_{\epsilon_{1}}(v)\right)\otimes e_{\alpha}f_{\lambda}(t)
\end{split}
\end{align} 
for any $u\in W,\,v\in\Fretw{\epsilon_{1}}$ and $t\in T_{\chi_{1}}$. 
Similarly, we get 
\begin{align}\label{twjacobi2}
\begin{split}
&\mathcal{Y}(u,z_{2})Y(a,z_{1})(v\otimes t)\\
&\quad=\left(p_{\epsilon_{2}}{\mathcal{Y}}_{\lambda}^{\rm tw}(\xi(u),z_{2})\left(\mathcal{Y}_{\alpha}^{\rm tw}(b,z_{1})+\mathcal{Y}_{-\alpha}^{\rm tw}(\theta(b),z_{1})\right)\iota_{\epsilon_{1}}(v)\right)\otimes f_{\lambda}(e_{\alpha}t).
\end{split}
\end{align}

{}From \eqref{jacobi2}, we see that for sufficiently large integer $k$,
\[
(z_{1}-z_{2})^{k}\mathcal{Y}_{\pm\alpha}^{\rm tw}(b,z_{1})
{\mathcal{Y}}_{\lambda}^{\rm tw}(\xi(u),z_{2})v
=(-1)^{(\alpha,\lambda)}(z_{1}-z_{2})^{k}{\mathcal{Y}}_{\lambda}^{\rm tw}(\xi(u),z_{2})
\mathcal{Y}_{\pm\alpha}^{\rm tw}(b,z_{1})v
\]
for $b\in\Fremo{\pm\alpha},\,u\in W$ and $v\in\Fretw{\epsilon_{1}}$, respectively.
Therefore, \eqref{twjacobi1} and \eqref{twjacobi2} shows that 
\begin{align}\label{twjacobi4}
\begin{split}
&(z_{1}-z_{2})^{k}\mathcal{Y}(u,z_{2})Y(a,z_{1})(v\otimes t)\\
&\quad=(z_{1}-z_{2})^{k}\left(p_{\epsilon_{2}}
{\mathcal{Y}}_{\lambda}^{\rm tw}(\xi(u),z_{2})
\left(\mathcal{Y}_{\alpha}^{\rm tw}(b,z_{1})
+\mathcal{Y}_{-\alpha}^{\rm tw}(\theta(b),z_{1})\right)\iota_{\epsilon_{1}}(v)\right)
\otimes f_{\lambda}(e_{\alpha}t)\\
&\quad=(-1)^{(\alpha,\lambda)}(z_{1}-z_{2})^{k}\left(p_{\epsilon_{2}}
\left(\mathcal{Y}_{\alpha}^{\rm tw}(b,z_{1})
+\mathcal{Y}_{-\alpha}^{\rm tw}(\theta(b),z_{1})\right)
\mathcal{Y}_{\lambda}^{\rm tw}(\xi(u),z_{2})\iota_{\epsilon_{1}}(v)\right)
\otimes f_{\lambda}(e_{\alpha}t)\\
&\quad=(-1)^{(\alpha,\lambda)}(z_{1}-z_{2})^{k}Y(a,z_{1})
\mathcal{Y}(u,z_{2})(v\otimes f_{\lambda}^{-1}((e_{\alpha})^{-1}f_{\lambda}(e_{\alpha}t))).
\end{split}
\end{align}

Since $\mathcal{Y}(\,\cdot\,,z)$ is an intertwining operator for $\charge{+}$, we have 
\[
(z_{1}-z_{2})^{k}\mathcal{Y}(u,z_{2})Y(a,z_{1})(v\otimes t)=(z_{1}-z_{2})^{k}Y(a,z_{1})\mathcal{Y}(u,z_{2})(v\otimes t)
\]
for large $k.$ 
Thus \eqref{twjacobi1}, \eqref{twjacobi4} and Proposition \ref{injectivity} imply the identity 
\begin{align}\label{llll}
e_{\alpha}f_{\lambda}(t)
=(-1)^{(\lambda,\alpha)}f_{\lambda}(e_{\alpha}t)
=f_{\lambda}(\sigma_{\lambda}(e_{\alpha})(t))
\end{align}
for any $\alpha\in L$ and $t\in T_{\chi_{1}}$. 
Therefore, 
$f_{\lambda}\in\Hom_{\hat{L}/K}(T_{\chi_{1}}^{(\lambda)},T_{\chi_{2}})$.
Consequently, we see that there exists an injective linear map 
\begin{align}\label{xxxx}
I_{\charge{+}}\fusion{M}{\charge{T_{\chi_{1}},
\epsilon_{1}}}{\charge{T_{\chi_{2}},\epsilon_{2}}}\rightarrow
I_{\Free{+}}\fusion{\Fremo{\lambda}}
{\Fretw{\epsilon_{1}}}{\Fretw{\epsilon_{2}}}\otimes
\Hom_{\hat{L}/K}(T_{\chi_{1}}^{(\lambda)},T_{\chi_{2}}).
\end{align} 
We have
$\dim_\C\Hom_{\hat{L}/K}(T_{\chi_{1}}^{(\lambda)},T_{\chi_{2}})
=\C\delta_{\chi_{1}^{(\lambda)},\chi_{2}}$.
Hence the dimension of the right hand side in \eqref{xxxx} is less
than or equal to $1$ by Theorem \ref{ehigher-rank-m1}.  We obtain the
following proposition:

\begin{proposition}\label{prop-32}
Let $M$ be an irreducible $\charge{+}$-module containing an
$\Free{+}$-submodule isomorphic to $\Fremo{\lambda}$ 
and let $\epsilon_{1},\,\epsilon_{2}\in\{\pm\}$.
Then the fusion rule of type 
$\fusion{M}{\charge{T_{\chi_{1}},\epsilon_{1}}}{\charge{T_{\chi_{2}},
\epsilon_{2}}}$ is zero if ${\chi_{2}}\neq {\chi_{1}}^{(\lambda)}$ 
and is less than or equal to $1$ if ${\chi_{2}}={\chi_{1}}^{(\lambda)}$.  
\end{proposition}

\begin{corollary}\label{rrrrr}
For any $\lambda\in L^{\circ}$ with $2\lambda\notin L$ and
$\epsilon_{1},\epsilon_{2}\in\{\pm\}$, the fusion rule of type
$\fusion{\charlam{\lambda}}{\charge{T_{\chi_{1}},\epsilon_{1}}}
{\charge{T_{\chi_{2}},\epsilon_{2}}}$
is $\delta_{\chi_{1}^{(\lambda)},\chi_{2}}$.
\end{corollary}
\begin{proof}It is clear from Propositions \ref{twist1} and \ref{prop-32}.
\end{proof}
Finally we prove the following proposition:
\begin{proposition}\label{prop-34}
Let $\lambda\in L^{\circ}$ with $2\lambda\in L$ and 
$\chi$ a central character of $\hat{L}/K$ such that $\chi(\kappa)=-1$.
Then 

(1) the fusion rules of types 
$\fusion{\charlam{\lambda}^{+}}{\charge{T_{\chi},\pm}}
{\charge{T_{\chi}^{(\lambda)},\mp}}$ 
and $\fusion{\charlam{\lambda}^{-}}{\charge{T_{\chi},\pm}}
{\charge{T_{\chi}^{(\lambda)},\pm}}$ are $0$ if $c_{\chi}(\lambda)=1$,

(2) the fusion rules of types 
$\fusion{\charlam{\lambda}^{+}}{\charge{T_{\chi},\pm}}
{\charge{T_{\chi}^{(\lambda)},\pm}}$ 
and $\fusion{\charlam{\lambda}^{-}}{\charge{T_{\chi},\pm}}
{\charge{T_{\chi}^{(\lambda)},\mp}}$ are zero if $c_{\chi}(\lambda)=-1$.
\end{proposition}
\begin{proof} 
We shall only prove that the fusion rule of type 
$\fusion{\charlam{\lambda}^{+}}{\charge{T_{\chi},+}}
{\charge{T_{\chi}^{(\lambda)},+}}$
is $0$ when $c_{\chi}(\lambda)=-1$. 
The others can be dealt similarly.

Let $\mathcal{Y}(\,\cdot\,,z)$ be an intertwining operator of type
$\fusion{\charlam{\lambda}^{+}}{\charge{T_{\chi},+}}
{\charge{T_{\chi}^{(\lambda)},+}}$
and $f_{\lambda}$ the linear map defined as in \eqref{wwww}.  As in
the proof of Proposition \ref{qqqqq}, we take a subset
$S_{\lambda}\subset\lambda+L$ such that
$\charlam{\lambda}^{+}\cong\bigoplus_{\mu\in
S_{\lambda}}\Fremo{\mu}\,(\cong\Free{+}\oplus\bigoplus_{\mu\in
S_{\lambda}}\Fremo{\mu}$ if $\lambda\in L$) as $\Free{+}$-modules.  We
recall the $\Free{+}$-module isomorphism
$\phi:\charlam{\lambda}^{+}[\lambda]\to\Fremo{\lambda}$ for
$\lambda\in S_{\lambda}$.  We may also assume that
$\lambda,3\lambda\in S_{\lambda}$.  By \eqref{jacobi2}, we have for
any $a\in\Fremo{2\lambda},\,u\in\Fremo{\lambda}$ and sufficiently
large integer $k$
\begin{align*}
&(z_{0}+z_{2})^{k}\mathcal{Y}^{\rm tw}_{2\lambda}(a,z_{0}+z_{2})
\mathcal{Y}^{\rm tw}_{\lambda}(u,z_{2})(v\otimes t)\\
&\quad=\frac{1}{2}\sum_{p=0,1}(z_{2}+z_{0})^{k}
\mathcal{Y}^{\rm tw}_{(1+(-1)^{p}2)\lambda}
(\mathcal{Y}_{(-1)^{p}2\lambda,\lambda}
(\theta^{p}(a),z_{0})u,z_{2})(v\otimes t)\\
&\quad=\frac{1}{2}(z_{2}+z_{0})^{k}
\left(\mathcal{Y}^{\rm tw}_{3\lambda}
(\mathcal{Y}_{2\lambda,\lambda}(a,z_{0})u,z_{2})
+\mathcal{Y}^{\rm tw}_{-\lambda}
(\mathcal{Y}_{-2\lambda,\lambda}(\theta(a),z_{0})u,z_{2})\right)(v\otimes t)
\end{align*}
for $v\in M(1)(\theta)^+$ and $t\in T_{\chi}.$ 
This and \eqref{twjacobi1} show that 
for $a=b+\theta(b)\in\charge{+}[2\lambda]$ 
with $b\in\Fremo{2\lambda},\,u\in\charlam{\lambda}^{+}[\lambda]$ 
and $v\in\Fretw{+}$,  
\begin{align}\label{8888}
\begin{split}
&(z_{0}+z_{2})^{k}Y(a,z_{0}+z_{2})\mathcal{Y}(u,z_{2})(v\otimes t)\\
&\,=(z_{0}+z_{2})^{k}\left(p_{+}^{\rm tw}\left(\mathcal{Y}_{2\lambda}^{\rm tw}(b,z_{0}+z_{2})+\mathcal{Y}_{-2\lambda}^{\rm tw}(\theta(b),z_{0}+z_{2})\right)\mathcal{Y}_{\lambda}^{\rm tw}(\phi(u),z_{2})v\right)\otimes e_{2\lambda}f_{\lambda}(t)\\
&\,=\frac{1}{2}(z_{2}+z_{0})^{k}\left(p_{+}^{\rm tw}\mathcal{Y}^{\rm tw}_{3\lambda}(\mathcal{Y}_{2\lambda,\lambda}(b,z_{0})\phi(u),z_{2})v+p_{+}^{\rm tw}\mathcal{Y}^{\rm tw}_{-\lambda}(\mathcal{Y}_{-2\lambda,\lambda}(\theta(b),z_{0})\phi(u),z_{2})v\right.\\
&\left.\qquad+p_{+}^{\rm tw}\mathcal{Y}^{\rm tw}_{-\lambda}(\mathcal{Y}_{-2\lambda,\lambda}(\theta(b),z_{0})\phi(u),z_{2})v+p_{+}^{\rm tw}\mathcal{Y}^{\rm tw}_{3\lambda}(\mathcal{Y}_{2\lambda,\lambda}(b,z_{0})\phi(u),z_{2})v\right)\otimes e_{2\lambda}f_{\lambda}(t)\\
&\,=(z_{2}+z_{0})^{k}\left(p_{+}^{\rm tw}\left(\mathcal{Y}^{\rm tw}_{3\lambda}(\mathcal{Y}_{2\lambda,\lambda}(b,z_{0})\phi(u),z_{2})v+\mathcal{Y}^{\rm tw}_{-\lambda}(\mathcal{Y}_{-2\lambda,\lambda}(\theta(b))\phi(u),z_{2})v\right)\right)\otimes e_{2\lambda}f_{\lambda}(t).
\end{split}
\end{align}

Consider $\phi(Y(a,z)u)$ for any $a\in\charge{+}[2\lambda]$ and $u\in\charlam{\lambda}^{+}[\lambda]$.
Note that $u=\phi(u)+\theta(\phi(u))$. We have  
\begin{align*}
Y(a,z)u=&\epsilon(2\lambda,\lambda)\left(\mathcal{Y}_{2\lambda,\lambda}(a,z)\phi(u)+\mathcal{Y}_{-2\lambda,-\lambda}(\theta(a),z)\theta(\phi(u))\right)\\
&+\epsilon(-2\lambda,\lambda)\left(\mathcal{Y}_{2\lambda,-\lambda}(a,z)\theta(\phi(u))+\mathcal{Y}_{-2\lambda,\lambda}(\theta(a),z)\phi(u)\right).
\end{align*}
Since $3\lambda,\,\lambda\in S_{\lambda}$, we have  
\[
\phi(Y(a,z)u)=\epsilon(2\lambda,\lambda)\mathcal{Y}_{2\lambda,\lambda}(a,z)\phi(u)+\epsilon(-2\lambda,\lambda)\mathcal{Y}_{2\lambda,-\lambda}(a,z)\theta(\phi(u)).
\]
Using \eqref{untwist1} gives 
\begin{align*}
\phi(Y(a,z)u)=\epsilon(2\lambda,\lambda)\mathcal{Y}_{2\lambda,\lambda}(a,z)\phi(u)+\epsilon(-2\lambda,\lambda)\theta\mathcal{Y}_{-2\lambda,\lambda}(\theta(a),z)\phi(u).
\end{align*}

Note that $p_{+}\mathcal{Y}^{\rm tw}_{\nu}(u,z)w=p_{+}\mathcal{Y}^{\rm tw}_{-\nu}(\theta(u),z)w$ for any $u\in\Fremo{\nu}$ and $w\in\Fretw{+}$.
Hence 
\begin{align}\label{hhhh}
\begin{split}
&\mathcal{Y}(Y(a,z_{0})u,z_{2})(v\otimes t)\\
&\quad=\epsilon(2\lambda,\lambda)p_{+}^{\rm tw}\mathcal{Y}^{\rm tw}_{3\lambda}(\mathcal{Y}_{2\lambda,\lambda}(a,z_{0})\phi(u),z_{2})v)\otimes f_{3\lambda}(t)\\
&\qquad+\epsilon(-2\lambda,\lambda)p_{+}^{\rm tw}\mathcal{Y}^{\rm tw}_{-\lambda}(\mathcal{Y}_{-2\lambda,\lambda}(\theta(a),z_{0})\phi(u),z_{2})v)\otimes f_{\lambda}(t).
\end{split}
\end{align}

On the other hand, \eqref{8888} gives 
\begin{align}\label{kkkk}
\begin{split}
&(z_{0}+z_{2})^{k}Y(a,z_{0}+z_{2})\mathcal{Y}(u,z_{2})(v\otimes t)\\
&\quad=(z_{2}+z_{0})^{k}\left(p_{+}^{\rm tw}\mathcal{Y}^{\rm tw}_{3\lambda}(\mathcal{Y}_{2\lambda,\lambda}(a,z_{0})\phi(u),z_{2})v\otimes e_{2\lambda}f_{\lambda}(t)\right.\\
&\qquad+\left.p_{+}^{\rm tw}\mathcal{Y}^{\rm tw}_{-\lambda}(\mathcal{Y}_{-2\lambda,\lambda}(\theta(a),z_{0})\phi(u),z_{2})v\otimes e_{2\lambda}f_{\lambda}(t)\right).
\end{split}
\end{align}
Since $\mathcal{Y}$ is an intertwining operator we have the identity
\begin{align*}
&(z_{0}+z_{2})^{k}Y(a,z_{0}+z_{2})\mathcal{Y}(u,z_{2})(v\otimes t)=(z_{2}+z_{0})^{k}\mathcal{Y}(Y(a,z_{0})u,z_{2})(v\otimes t)
\end{align*}
for sufficiently large integer $k$.
It follows from  \eqref{hhhh} and \eqref{kkkk} that 
\begin{multline*}
p_{+}^{\rm tw}\mathcal{Y}^{\rm tw}_{3\lambda}(\mathcal{Y}_{2\lambda,\lambda}(a,z_{0})\phi(u),z_{2})v)\otimes (\epsilon(2\lambda,\lambda)f_{3\lambda}(t)-e_{2\lambda}f_{\lambda}(t))\\
=p_{+}^{\rm tw}\mathcal{Y}^{\rm tw}_{-\lambda}(\mathcal{Y}_{-2\lambda,\lambda}(\theta(a),z_{0})\phi(u),z_{2})v)\otimes (e_{2\lambda}f_{\lambda}(t)-\epsilon(-2\lambda,\lambda)f_{\lambda}(t)).
\end{multline*}
Since the least powers of $z_{0}$ in 
\[
\mathcal{Y}^{\rm tw}_{3\lambda}(\mathcal{Y}_{2\lambda,\lambda}(e^{2\lambda},z_{0})e^{\lambda},z_{2})v\quad{\rm and}\quad\mathcal{Y}^{\rm tw}_{-\lambda}(\mathcal{Y}_{-2\lambda,\lambda}(e^{-2\lambda},z_{0})e^{\lambda},z_{2})v
\] 
are ${(2\lambda,\lambda)}$ and ${-(2\lambda,\lambda)}$ respectively, 
we see that if $\lambda\neq0$,  
\begin{align}\label{2222}
\chi^{(\lambda)}(e_{2\lambda})f_{\lambda}(t)
=e_{2\lambda}f_{\lambda}(t)=\epsilon(-2\lambda,\lambda)f_{\lambda}(t)
\end{align}
for any $t\in T_{\chi}$. That is,
\begin{align}\label{3333}
c_{\chi}(\lambda)f_{\lambda}(t)=f_{\lambda}(t).
\end{align}
The condition $c_{\chi}(\lambda)=-1$ forces $f_{\lambda}=0$. 
This shows $\mathcal{Y}(\,\cdot\,,z)=0$. 
\end{proof}

By Propositions \ref{duality}, \ref{twist1}--\ref{prop-32}, \ref{prop-34} and 
Corollary \ref{rrrrr} we immediately have:

\begin{proposition}\label{jjjjjj} 
Let $M^i\,(i=1,2,3)$ be irreducible $\charge{+}$-modules and 
assume that one of them is of twisted type.
Then the fusion rule of type $\fusion{M^1}{M^2}{M^3}$ is either $0$ or $1$. 
The fusion rule of type $\fusion{M^1}{M^2}{M^3}$ is $1$ if and only if 
$M^i\,(i=1,2,3)$ satisfy the following conditions;

\noindent
{\rm(i)} $M^1=\charlam{\lambda}$ for $\lambda\in L^{\circ}$ such that 
$2\lambda\notin L$ and $(M^2,\,M^3)$ is one of the following pairs:
\begin{itemize}
\item[] $(\charge{T_{\chi},\pm},\charge{T_{\chi}^{(\lambda)},\pm}),
\,(\charge{T_{\chi},\pm},
\charge{T_{\chi}^{(\lambda)},\mp})$ 
for any irreducible $\hat{L}/K$-module $T_{\chi}$. 
\end{itemize}

\noindent
{\rm(ii)} $M^1=\charlam{\lambda}^+$ for $\lambda\in L^{\circ}$ such that 
$2\lambda\in L$ and $(M^2,\,M^3)$ is one of the following pairs:
\begin{itemize}
\item[] $(\charge{T_{\chi},\pm},\charge{T_{\chi}^{(\lambda)},\pm}),
\,((\charge{T_{\chi}^{(\lambda)},\pm})',(\charge{T_{\chi},\pm})')$ 
for any irreducible $\hat{L}/K$-module $T_{\chi}$ 
such that $c_{\chi}(\lambda)=1$,
\item[] $(\charge{T_{\chi},\pm},\charge{T_{\chi}^{(\lambda)},\mp}),
\,((\charge{T_{\chi}^{(\lambda)},\pm})',(\charge{T_{\chi},\mp})')$ 
for any irreducible $\hat{L}/K$-module $T_{\chi}$ 
such that $c_{\chi}(\lambda)=-1$.
\end{itemize} 

\noindent
{\rm(iii)} $M^1=\charlam{\lambda}^-$ for $\lambda\in L^{\circ}$ 
such that $2\lambda\in L$  and $(M^2,\,M^3)$ is one of the following pairs:
\begin{itemize}
\item[] $(\charge{T_{\chi},\pm},\charge{T_{\chi}^{(\lambda)},\mp}),
\,((\charge{T_{\chi}^{(\lambda)},\mp})',(\charge{T_{\chi},\pm})')$ 
for any irreducible $\hat{L}/K$-module $T_{\chi}$ 
such that $c_{\chi}(\lambda)=1$,
\item[] $(\charge{T_{\chi},\pm},\charge{T_{\chi}^{(\lambda)},\pm}),\,((\charge{T_{\chi}^{(\lambda)},\pm})',(\charge{T_{\chi},\pm})')$ for any irreducible $\hat{L}/K$-module $T_{\chi}$ such that $c_{\chi}(\lambda)=-1$.
\end{itemize} 
 
\noindent
{\rm(iv)} $M^1=\charge{T_{\chi},+}$ for an irreducible $\hat{L}/K$-module $T_{\chi}$ and $(M^2,\,M^3)$ is one of the following pairs:
\begin{itemize}
\item[] $(\charlam{\lambda},\charge{T_{\chi}^{(\lambda)},\pm}),\,((\charge{T_{\chi}^{(\lambda)},\pm})',(\charlam{\lambda})')$ for $\lambda\in L^{\circ}$ such that $2\lambda\notin L$,
\item[] $(\charlam{\lambda}^{\pm},\charge{T_{\chi}^{(\lambda)},\pm}),\,((\charge{T_{\chi}^{(\lambda)},\pm})',(\charlam{\lambda}^{\pm})')$ for $\lambda\in L^{\circ}$ such that $2\lambda\in L$ and that $c_{\chi}(\lambda)=1$,
\item[] $(\charlam{\lambda}^{\pm},\charge{T_{\chi}^{(\lambda)},\mp}),\,((\charge{T_{\chi}^{(\lambda)},\mp})',(\charlam{\lambda}^{\pm})')$ for $\lambda\in L^{\circ}$ such that $2\lambda\in L$ and that $c_{\chi}(\lambda)=-1$.
\end{itemize} 

\noindent
{\rm(v)} $M^1=\charge{T_{\chi},-}$ for an irreducible $\hat{L}/K$-module $T_{\chi}$ and $(M^2,\,M^3)$ is one of the following pairs:
\begin{itemize}
\item[] $(\charlam{\lambda},\charge{T_{\chi}^{(\lambda)},\pm}),\,((\charge{T_{\chi}^{(\lambda)},\pm})',(\charlam{\lambda})')$ for $\lambda\in L$ such that $2\lambda\notin L^{\circ}$,
\item[] $(\charlam{\lambda}^{\pm},\charge{T_{\chi}^{(\lambda)},\mp}),\,((\charge{T_{\chi}^{(\lambda)},\pm})',(\charlam{\lambda}^{\mp})')$ for $\lambda\in L^{\circ}$ such that $2\lambda\in L$ and that $c_{\chi}(\lambda)=1$,
\item[] $(\charlam{\lambda}^{\pm},\charge{T_{\chi}^{(\lambda)},\pm}),\,((\charge{T_{\chi}^{(\lambda)},\mp})',(\charlam{\lambda}^{\mp})')$ for $\lambda\in L^{\circ}$ such that $2\lambda\in L$ and that $c_{\chi}(\lambda)=-1$.
\end{itemize} 
\end{proposition}

Now Theorem \ref{fusioncharge} follows from
Propositions \ref{fffffff} and \ref{jjjjjj}. 

\subsection{Fusion product for $\charge{+}$}
Let $V$ be a vertex operator algebra and $\{W_{i}\}_{i\in I}$ 
the set of all the equivalence classes of irreducible $V$-modules.
For any representatives $M^{i}$ of $W_i\,(i\in I)$, 
we write $N_{ij}^{k}$ for the fusion rule of 
type $\fusion{M^{i}}{M^{j}}{M^{k}}$ for $i,j,k\in I$.
The fusion rules $N_{ij}^{k}$ are independent of a choice of representatives. 
Here and further we assume that $I$ is a finite set and 
that all the fusion rules are finite. 

Set $\mathcal{R}(V)=\bigoplus_{i\in I}\C W_{i}$ 
a vector space over $\C$ with basis $\{W_{i}\}_{i\in I}$. 
Then the product of $\mathcal{R}(V)$ is defined by
\[
W_{i}\times W_{j}=\sum_{k}N_{ij}^{k}\,W_{k}
\]
for any $i,j\in I$. 
By Proposition \ref{duality} the product is commutative. 
%but in general is not associative. 
The commutative algebra $\mathcal{R}(V)$ is called 
the {\rm fusion algebra} of $V$. 

Denote by $W_{i}'$ the equivalence class of the contragredient module of
a representative of $W_{i}$.
Then for any $i\in I$ there exists uniquely $i'\in I$ 
such that $W_{i}'=W_{i'}$.
By Proposition \ref{duality}, we have $N_{ij}^{k}=N_{ik'}^{j'}$ 
for any $i,j,k\in I$. 

We now describe the fusion products for $\charge{+}$.  For simplicity,
we introduce notations of equivalence classes of irreducible
$\charge{+}$-modules.  For $\lambda\in L^{\circ}$, we set $[\lambda]$
be the equivalent class of irreducible $\charge{+}$-modules isomorphic
to $\charlam{\lambda}$.  When $2\lambda\in L$, we denote by
$[\lambda]^{\pm}$ the equivalent class of irreducible
$\charge{+}$-modules isomorphic to $\charlam{\lambda}^{\pm}$.  By
abuse of notations we set $[\lambda]=[\lambda]^++[\lambda]^-$ for
$\lambda\in L^{\circ}$ with $2\lambda\in L$.  We then have that
$[\lambda]=[-\lambda]$ and $[\lambda+\alpha]=[\lambda]$ for any
$\lambda\in L^{\circ}$ and $\alpha\in L$.  This implies that
$[\lambda+\mu]=[\lambda-\mu]$ for $\lambda,\mu\in L^{\circ}$ if
$2\lambda\in L$ or $2\mu\in L$.  For a central character $\chi$ of
$\hat{L}/K$ with $\chi(\kappa)=-1$, we write $[\chi]^{\pm}$ for the
equivalence classes of irreducible $\charge{+}$-modules
$\charge{T_{\chi},\pm}$, respectively.
 
%From Proposition \ref{eeee} we get the following identification:
%\begin{align}
%[\lambda]'&=[-\lambda]\quad\hbox{for $\lambda\in L^{\circ}$ such that $2\lambda\notin L$},\\
%([\lambda]^{\pm})'&=[\lambda]^{\pm}\quad\hbox{for $\lambda\in L^{\circ}$ such that $2\lambda\in L$ and $2(\lambda,\lambda)$ is even},\\
%([\lambda]^{\pm})'&=[\lambda]^{\mp}\quad\hbox{for $\lambda\in L^{\circ}$ such that $2\lambda\in L$ and $2(\lambda,\lambda)$ is odd},\\
%([\chi]^{\pm})'&=[\chi']^{\pm}\quad\hbox{for a central character $\chi$ such that $\chi(\kappa)=-1$}.
%\end{align}

We set $S_{0}=\{\lambda\in L^{\circ}|2\lambda\in L\}$ and $S_{1}=\{\lambda\in L^{\circ}|2\lambda\notin L\}$.
By Theorem \ref{fusioncharge} we have the following fusion products: 
\begin{align}
[\lambda]\times[\mu]&=[\lambda+\mu]+[\lambda-\mu]\quad\hbox{for $\lambda,\mu\in S_{1}$}\\
[\lambda]^{\pm}\times[\mu]&=[\lambda+\mu]\quad\hbox{for $\lambda\in S_{0},\,\mu\in S_{1}$},\\
[\lambda]^{+}\times[\mu]^{\pm}&=[\lambda+\mu]^{\pm}\quad\hbox{for $\lambda,\,\mu\in S_{0}$ such that $\pi_{\lambda,2\mu}=1$},\label{ooooooo}\\
[\lambda]^{+}\times[\mu]^{\pm}&=[\lambda+\mu]^{\mp}\quad\hbox{for $\lambda,\,\mu\in S_{0}$ such that $\pi_{\lambda,2\mu}=-1$},\label{pppppp}\\
[\lambda]^{-}\times[\mu]^{-}&=[\lambda+\mu]^{+}\quad\hbox{for $\lambda,\,\mu\in S_{0}$ such that $\pi_{\lambda,2\mu}=1$},\\
[\lambda]^{-}\times[\mu]^{-}&=[\lambda+\mu]^{-}\quad\hbox{for $\lambda,\,\mu\in S_{0}$ such that $\pi_{\lambda,2\mu}=-1$},\\
[\lambda]\times[\chi]^{\pm}&=[\chi^{(\lambda)}]^{+}+[\chi^{(\lambda)}]^{-}\quad\hbox{for $\lambda\in S_{1}$},\\
[\lambda]^{+}\times[\chi]^{\pm}&=[\chi^{(\lambda)}]^{\pm}\quad\hbox{for $\lambda\in S_{0}$ such that $c_{\chi}(\lambda)=1$},\\
[\lambda]^{+}\times[\chi]^{\pm}&=[\chi^{(\lambda)}]^{\mp}\quad\hbox{for $\lambda\in S_{0}$ such that $c_{\chi}(\lambda)=-1$},\label{qwer}\\
[\lambda]^{-}\times[\chi]^{-}&=[\chi^{(\lambda)}]^{+}\quad\hbox{for $\lambda\in S_{0}$ such that $c_{\chi}(\lambda)=1$},\\
[\lambda]^{-}\times[\chi]^{-}&=[\chi^{(\lambda)}]^{-}\quad\hbox{for $\lambda\in S_{0}$ such that $c_{\chi}(\lambda)=-1$}.
\end{align}
The other products are derived from these products with the symmetries of fusion rules in Proposition \ref{duality}. 
For example, the product of $[\chi_{1}]^{+}$ and $[\chi_{2}]^{+}$ is given by
\[
[\chi_{1}]^{+}\times[\chi_{2}]^{+}=\sum[\lambda]+\sum[\mu]^{+}+\sum[\nu]^{-},
\]
where $\lambda$ runs through $S_{1}$ such that $\chi_{1}^{(\lambda)}=\chi_{2}'$, $\mu$ runs through $S_{0}$ such that $\chi_{1}^{(\mu)}=\chi_{2}'$ and that $c_{\chi_{1}}(\mu)(-1)^{2(\mu,\mu)}=1$, and $\nu$ runs through $S_{0}$ such that $\chi_{1}^{(\nu)}=\chi_{2}'$ and that $c_{\chi_{1}}(\mu)(-1)^{2(\mu,\mu)}=-1$.

\begin{theorem}\label{ppppp}
The fusion algebra $\mathcal{R}(\charge{+})$ is a commutative associative algebra.
\end{theorem}
\begin{proof}
For any equivalence classes $W_{1},W_{2}$ and $W_{3}$ of irreducible $\charge{+}$-modules, we have to prove that $W_{1}\times(W_{2}\times W_{3})=(W_{1}\times W_{2})\times W_{3}$.
We can do this case by case.
For example we shall prove 
\begin{align}\label{9ki}
[\lambda]^{+}\times([\mu]^{+}\times[\nu]^{-})=([\lambda]^{+}\times[\mu]^{+})\times[\nu]^{-}
\end{align} 
for $\lambda,\mu,\nu\in S_{0}$ such that $\pi_{\lambda,2\mu}=1$ and $\pi_{\mu,2\nu}=-1$ and 
\begin{align}\label{yht}
[\lambda]^{+}\times([\mu]^{+}\times[\chi]^{-})=([\lambda]^{+}\times[\mu]^{+})\times[\chi]^{-}
\end{align} 
for $\lambda,\mu\in S_{0}$ and a central character $\chi$ of $\hat{L}/K$ with $\chi(\kappa)=-1$ such that $\pi_{\lambda,2\mu}=1$ and $c_{\chi}(\mu)=-1$.
We first show \eqref{9ki}. 
By using \eqref{ooooooo} and \eqref{pppppp}, we have 
\[
[\mu]^{+}\times[\nu]^{-}=[\mu+\nu]^{+},\quad [\lambda]^{+}\times[\mu]^{+}=[\lambda+\mu]^{+}.
\]
Since 
\begin{align*}
&\pi_{\lambda,2\mu+2\nu}=e^{(\lambda,2\mu+2\nu)\pi i}\w_{q}^{c_{0}(2\mu+2\nu,\lambda)}=\pi_{\lambda,2\mu}\pi_{\lambda,2\nu}=\pi_{\lambda,2\nu},\\
&\pi_{\lambda+\mu,2\nu}=e^{(\lambda+\mu,2\nu)\pi i}\w_{q}^{c_{0}(2\nu,\lambda+\mu)}=\pi_{\lambda,2\nu}\pi_{\mu,2\nu}=-\pi_{\lambda,2\nu},
\end{align*}
we see that 
\begin{align*}
&[\lambda]^{+}\times([\mu]^{+}\times[\nu]^{-})=[\lambda]^{+}\times[\mu+\nu]^{+}=[\lambda+\mu+\nu]^{\pm},\\
&([\lambda]^{+}\times[\mu]^{+})\times[\nu]^{-}=[\lambda+\mu]^{+}\times[\nu]^{-}=[\lambda+\mu+\nu]^{\pm}
\end{align*}
if $\pi_{\lambda,2\nu}=\pm1$ respectively.
Thus \eqref{9ki} holds. 

Next we show \eqref{yht}.
Equation \eqref{qwer} imply 
\[
[\mu]^{+}\times[\chi]^{-}=[\chi^{(\mu)}]^{+}.
\]
Then we see that 
\begin{align*}
&[\lambda]^{+}\times([\mu]^{+}\times[\chi]^{-})=[\lambda]^{+}\times[\chi^{(\mu)}]^{-}=[\left(\chi^{(\mu)}\right)^{(\lambda)}]^{\pm}
\end{align*}
if $c_{\chi}(\lambda)=\pm1$ respectively.
On the other hand, 
\[
([\lambda]^{+}\times[\mu]^{+})\times[\chi]^{-}=[\lambda+\mu]^{+}\times[\chi]^{-}=[\chi^{(\lambda+\mu)}]^{\mp}
\]
if $c_{\chi}(\lambda+\mu)=-c_{\chi}(\lambda)=\pm1$ respectively.
Since $\left(\chi^{(\mu)}\right)^{(\lambda)}(a)=\chi^{(\mu)}(\sigma_{\lambda}(a))=\chi(\sigma_{\mu}\sigma_{\lambda}(a))=\chi(\sigma_{\lambda+\mu}(a))=\chi^{(\lambda+\mu)}(a)$ for any $a\in Z(\hat{L}/K)$, we have $\left(\chi^{(\mu)}\right)^{(\lambda)}=\chi^{(\lambda+\mu)}$. 
Therefore, \eqref{yht} holds.  
\end{proof}

\subsection{Application}

In this section we apply the results on fusion rules for $V_L^+$-modules
to study orbifold vertex operator algebras constructed from $V_L$ 
and automorphism $\theta$ when $L$ is unimodular.

Let $L$ be a positive-definite even unimodular lattice.  That is,
$L=L^{\circ}.$ Then $V_L$ is a holomorphic vertex operator algebra in
the sense that $V_L$ is rational and $V_L$ is the only irreducible
$V_L$-module up to isomorphism (see \cite{D1} and
\cite{DLM2}). Moreover, $V_L$ has a unique irreducible
$\theta$-twisted module $V_L^T$ up to isomorphism where $T$ is the
unique simple module for $\hat L/K$ such that $\kappa K$ acts as $-1$
(see \cite{FLM} and \cite{D2}). Recall that $V_L^T=M(1)(\theta)\otimes
T$ and $d$ is the rank of $L.$ The weight gradation of $V_L^T$ is
given by
\begin{align}\label{grad2}
V_L^T=\sum_{n\in\frac{1}{2}\Z_{\geq 0}}(V_L^T)_{n+\frac{d}{16}}
\end{align}
(see \cite{DL2}).  Since $L$ is unimodular, $d$ is divisible by $8$.
Hence the $L(0)$-weights of either $\charge{T,+}$ or $\charge{T,-}$
are integers (half integers for the other).  We denote by
$\charge{T,e}$ (resp, $\charge{T,o}$) the irreducible
$\charge{+}$-submodules of $\charge{T}$ with integral (half integral)
$L(0)$-weights.  It is clear that $\charge{T,e}=\charge{T,+}$ if $d/8$
is even and $\charge{T,e}=\charge{T,-}$ if $d/8$ is odd.  By Theorems
\ref{chargeclass} (also see \cite{AD}) and \ref{fusioncharge}, we
have:

\begin{proposition}\label{special} 
Let $L$ be a positive-definite even unimodular lattice.

\noindent
{\rm(i)} The vertex operator algebra 
$V_L^+$ has exactly $4$ irreducible modules 
$V_L^{\pm}, V_L^{T,\pm}$ up to isomorphism.

\noindent
{\rm(ii)} The fusion rules among modules are
\begin{align*} 
&V_L^+\times W=W\times V_L^+=W,\quad V_L^{-}\times V_L^-=V_L^+,\\
&V_L^-\times \charge{T,e}= \charge{T,e}\times V_L^-=V_L^{T,o},\quad V_L^-\times V_L^{T,o}= V_L^{T,o}\times V_L^-=\charge{T,e},\\
&\charge{T,e}\times \charge{T,e}=V_L^{T,o}\times V_L^{T,o}=V_L^{+},\quad \charge{T,e}\times V_L^{T,o}=V_L^{T,o}\times \charge{T,e}=V_L^{-},
\end{align*}
where $W$ is any irreducible $V_L^+$-module.
\end{proposition}
\begin{remark}
If $L$ is the Leech lattice, the irreducible modules for $V_L^+$ has
been classified previously in \cite{D3} by using the representation theory
for the Virasoro algebra of central charge $1/2.$ 
\end{remark}

The main result in this subsection is the following:
\begin{proposition}\label{holomorph} 
Let $L$ be a positive-definite even unimodular lattice.
Assume that $V_L^+$ is rational and $V=V_L^++V_L^{T,e}$ is
a vertex operator algebra. Then $V$ is a holomorphic vertex operator
algebra and $C_2$-cofinite.
\end{proposition}

\begin{proof} 
It is known that $\charge{+}$ is $C_2$-cofinite (see \cite{Yam} and
\cite{ABD}).  Since $V$ is
$C_2$-cofinite as a $\charge{+}$-module by \cite{Bu}, it is also
$C_2$-cofinite as a vertex operator algebra.

We assume that $V=V_L^++V_L^{T,-}.$ The case
that $V=V_L^++V_L^{T,+}$ can be proved similarly.
We first prove that $V$ is the only irreducible $V$-module
up to isomorphism. Let $W$ be an irreducible $V$-module. Then $W$ is
a completely reducible $V_L^+$-module. Let $M$ be an irreducible
$V_L^+$-submodule of $W.$ If $M=V_L^+$ or $V_L^{T,-}$ using the fusion
rule given in Proposition \ref{special} shows that $V_L^+$ 
is always contained in $W$ as a $V_L^+$-submodule. So
$W$ contains a vacuum like vector and thus isomorphic to $V$ 
(see \cite{liform}). 

If $M=V_L^-,$ then $V_L^{T,-}\times V_L^-=V_L^{T,+}$ is a $V_L^+$-submodule
of $W.$ Note that $V_L^{T,-}$ has integral weight and  $V_L^{T,+}$ has 
strictly half integral weights. So  $W$ has both integral weights
from $V_L^-$ and half integral weights from $V_L^{T,+}.$ 
But this is impossible as $W$ is irreducible. Similarly, $M$ cannot
be $V_L^{T,+}.$ 

We now prove that $V$ is rational. That is, any admissible $V$-module
is completely reducible. Let $W$ be an admissible $V$-module
and $M$ be the maximal semisimple admissible submodule. Then
$V=M\oplus X$ for a $V_L^+$-submodule of $W$ as $V_L^+$ is
rational. If $X\ne 0$ then $W/M$ is a $V$-module. So as $V_L^+$-module
$W/M=X$ contains a $V_L^+$-submodule isomorphic to $V_L^+.$ 
This shows that $X$ contains a vacuum-like vector $x$ and the 
$V$-submodule $Z$ of $W$ generated by $u$ is isomorphic to $V.$
Clearly, $M\cap Z=0$ and $M\oplus Z$ is a semisimple 
admissible $V$-submodule of $W$ and strictly contains 
$M.$ This contradiction shows that $W=M.$
\end{proof}

Again, if $L$ is the Leech lattice, this result has been given in
\cite{D3} before.

\end{document}